\documentclass[12pt,a4paper]{article}
\usepackage[latin1]{inputenc}
\usepackage{amsmath,amsfonts,amsthm,enumitem,aliascnt,graphicx,multirow,natbib,hyperref}

\numberwithin{equation}{section}

\theoremstyle{plain}
\newtheorem{theorem}{Theorem}

\newaliascnt{corollary}{theorem}

\aliascntresetthe{corollary}

\newaliascnt{lemma}{theorem}
\newtheorem{lemma}[lemma]{Lemma}
\aliascntresetthe{lemma}

\newaliascnt{proposition}{theorem}
\newtheorem{proposition}[proposition]{Proposition}
\aliascntresetthe{proposition}

\theoremstyle{definition}
\newtheorem{example}{Example}

\newcommand{\midil}{\,|\,}
\allowdisplaybreaks[4]

\makeatletter
\def\blfootnote{\xdef\@thefnmark{}\@footnotetext}
\makeatother

\author{Joydeep Chowdhury and Probal Chaudhuri\\
Indian Statistical Institute, Kolkata}
\title{Convergence rates for kernel regression in infinite dimensional spaces}
\date{}

\begin{document}
\maketitle

\begin{abstract}
We consider a nonparametric regression setup, where the covariate is a random element in a complete separable metric space, and the parameter of interest associated with the conditional distribution of the response lies in a separable Banach space. We derive the optimum convergence rate for the kernel estimate of the parameter in this setup. The small ball probability in the covariate space plays a critical role in determining the asymptotic variance of kernel estimates. Unlike the case of finite dimensional covariates, we show that the asymptotic orders of the bias and the variance of the estimate achieving the optimum convergence rate may be different for infinite dimensional covariates. Also, the bandwidth, which balances the bias and the variance, may lead to an estimate with suboptimal mean square error for infinite dimensional covariates. We describe a data-driven adaptive choice of the bandwidth, and derive the asymptotic behavior of the adaptive estimate.
\\

\textit{Keywords:}
Adaptive estimate, Bias variance decomposition, Gaussian process, Maximum likelihood regression, Mean square error, Optimal bandwidth, Small ball probability, \textit{t} process
\end{abstract}

\blfootnote{Accepted in AISM \url{http://www.ism.ac.jp/editsec/aism/}}

\section{Introduction} \label{sec:1}
Suppose that we have a nonparametric regression problem, where the covariate $ \mathbf{X} $ is a random element in a complete separable metric space, and the response $ \mathbf{Y} $ lies in some arbitrary measure space. Our parameter of interest, which is denoted as $ \Theta( \mathbf{x} ) $, is a parameter associated with the conditional distribution of $ \mathbf{Y} $ given $ \mathbf{X} = \mathbf{x} $.
Let $ ( \mathbf{X}_1, \mathbf{Y}_1 ), \ldots, ( \mathbf{X}_n, \mathbf{Y}_n ) $ be a sample of \textit{i.i.d.} observations from the joint distribution of $ ( \mathbf{X}, \mathbf{Y} ) $, and our objective is to estimate $ \Theta( \mathbf{x} ) $ based on this sample.
In the particular case, where the response $ \mathbf{Y} $ is a real random variable, the covariate space is the $ q $-dimensional Euclidean space $ \mathbb{R}^q $ and $ \Theta( \mathbf{x} ) = \mathbb{E}[ \mathbf{Y} \midil \mathbf{X} = \mathbf{x} ] $, \cite{stone1980optimal} proved that the optimal convergence rate of a nonparametric estimate $ \widehat{\Theta}_n( \mathbf{x} ) $ of $ \Theta( \mathbf{x} ) $ is $ n^{- ( \beta / (2 \beta + q) ) } $.
Here, $ \beta $ is a positive constant such that $ | \Theta( \mathbf{z} ) - \Theta( \mathbf{x} ) | = O( \| \mathbf{z} - \mathbf{x} \|^\beta ) $ as $ \mathbf{z} \longrightarrow \mathbf{x} $, with $ \| \cdot \| $ being the Euclidean norm in $ \mathbb{R}^q $.
The optimum achievable convergence rate for nonparametric regression with finite dimensional covariate was further investigated in \cite{stone1982optimal}, \cite{ibragimov1980nonparametric}, \cite{yatracos1988lower}, \cite{donoho1991geometrizingII,donoho1991geometrizingIII}, etc.
However, when the dimension of the covariate space is infinite, the expressions of the optimum rate of convergence derived by these authors are no longer valid.

Recently nonparametric regression with functional covariates has been studied in \cite{masry2005nonparametric}, \cite{ferraty2007nonparametric}, \cite{rachdi2007nonparametric}, etc. These authors investigated nonparametric estimation of the conditional mean when the covariate is functional, and the response is real-valued. They investigated the consistency and the asymptotic normality of kernel estimates as well as data-driven selection of the bandwidth.
In \cite{ferraty2012regression} and \cite{lian2012convergence}, the problem of nonparametric regression when both the response and the covariate are functions is investigated, where the parameter of interest is the conditional mean of the response given the covariate. In \cite{ferraty2012regression}, asymptotic normality of the estimate of the conditional mean is derived and a bootstrap implementation is described. In \cite{lian2012convergence}, an upper bound of the convergence rate of the estimate of the conditional mean is established.

The problem of optimum convergence rate of a nonparametric regression estimate was explored in \cite{mas2012lower} and \cite{chagny2014adaptive} when the covariate is infinite dimensional.
In \cite{mas2012lower}, the usual mean regression problem with a real-valued response was considered, and a lower bound for the rate of convergence of the minimax risk was established (see Theorem 3 in \cite{mas2012lower}).
In \cite{chagny2014adaptive}, the optimum convergence rate was derived for the estimate of the conditional distribution function of a real-valued response given an infinite dimensional covariate.
In both these cases, the methodology developed is restricted to the conditional mean of some real-valued response, and cannot be applied when the response is infinite dimensional, or the parameter of interest is not the conditional mean.

In most of the existing literature on nonparametric regression with functional data, the authors considered real or multivariate responses and functional covariates. However, regression problems, where the response itself may be infinite dimensional in nature, are also common. Authors who investigated regression with functional responses and covariates, like \cite{ferraty2012regression} and \cite{lian2012convergence}, considered the conditional mean as their parameter of interest. But, one may also be interested in various parameters of the conditional distribution of the response other than the conditional mean, like the conditional variance and covariance, the conditional coefficient of variation, the conditional correlation, etc.
In our study, the parameter of interest $ \Theta( \mathbf{x} ) $ lies in some separable Banach space, and it covers a large class of parameters of interest including those stated above. We shall investigate the convergence rate of a large class of kernel estimates in this setup and derive the optimal convergence rate.

In \autoref{sec:2}, our regression setup and the kernel estimates are described in detail. In \autoref{sec:3}, we discuss an asymptotic bias--variance decomposition of our kernel estimate, and study the asymptotic behavior of the bias and the variance terms. We show that the asymptotic behavior of the variance term critically depends on the small ball probability in the covariate space. The main convergence results for the estimate $ \widehat{\Theta}_n( \mathbf{x} ) $ are presented in \autoref{sec:4}. A data-driven method of bandwidth selection along with the asymptotic behavior of the adaptive estimate is presented in \autoref{sec:adaptive}. In the same section, we demonstrate the adaptive estimates in simulated datasets from several regression models. \autoref{sec:5} contains concluding remarks and discussion. The proofs and related mathematical details are provided in \autoref{sec:6}.

\section{Kernel estimates} \label{sec:2}
We assume that the covariate $ \mathbf{X} $ is a random element in some complete separable metric space $ ( \mathcal{C}, d ) $ with $ d $ being the metric, and the response $ \mathbf{Y} $ is a random element in some measurable space $ \mathcal{R} $ equipped with some appropriate $ \sigma $-field and probability measure. Denote the conditional probability measure of $ \mathbf{Y} $ given $ \mathbf{X} = \mathbf{x} $ as $ \mu( \cdot \midil \mathbf{x} ) $. We want to estimate a parameter $ \Theta( \mathbf{x} ) $ associated with $ \mu( \cdot \midil \mathbf{x} ) $ for a fixed $ \mathbf{x} \in \mathcal{C} $. We employ the nonparametric kernel regression method developed by \cite{nadaraya1964estimating} and \cite{watson1964smooth}. Let $ K( \cdot ) $ be a suitable kernel function with associated bandwidth $ h > 0 $. To estimate $ \Theta( \mathbf{x} ) $, we first construct the weighted empirical probability measure $ \mu_n( \cdot \midil \mathbf{x} ) $ that assigns probability mass
\begin{align*}
W_{i,n} = \frac{K( h^{-1} d( \mathbf{x}, \mathbf{X}_i ) )}{\sum_{i=1}^{n} K( h^{-1} d( \mathbf{x}, \mathbf{X}_i ) )}
\end{align*}
to the data point $ \mathbf{Y}_i $ for $ i = 1, \ldots, n $. The kernel estimate $ \widehat{\Theta}_n( \mathbf{x} ) $ of $ \Theta( \mathbf{x} ) $ is the corresponding parameter associated with $ \mu_n( \cdot \midil \mathbf{x} ) $.

We require the concept of a type 2 Banach spaces in our subsequent discussion. A separable Banach space is called type 2 if there is a positive constant $ c $ such that for any $ n \ge 1 $ and independent zero-mean random elements $ \mathbf{Z}_1, \ldots, \mathbf{Z}_n $ in that Banach space with $ \mathbb{E}\| \mathbf{Z}_i \|^2 < \infty $ for $ i = 1, \ldots, n $, we have $ \mathbb{E}\| \mathbf{Z}_1 + \cdots + \mathbf{Z}_n \|^2 \le c ( \mathbb{E}\| \mathbf{Z}_1 \|^2 + \cdots + \mathbb{E}\| \mathbf{Z}_n \|^2 ) $ \citep[p.~158]{araujo1980central}.
Also, a Banach space is said to have a Schauder basis $ \{ e_n \} $ if for every element $ v $ in that space, $ v = \sum_{n=1}^{\infty} v_n e_n $ for some sequence of real numbers $ \{ v_n \} $.
Separable Hilbert spaces and $ L_p[ a, b ] $ spaces with $ p \ge 2 $ and $ -\infty \le a < b \le \infty $ are well-known examples of type 2 Banach spaces with Schauder bases.
Henceforth, $ \mathbb{I}( \cdot ) $ will denote the usual indicator function.

We now discuss some examples. These examples demonstrate the use of kernel estimates in a diverse class of statistical models. In all our subsequent discussion, the expectation of a random element in a separable Banach space is defined in the sense of Bochner \citep[p.~100]{araujo1980central}.

\begin{example}[\textbf{Mean regression}] \label{ex:2.1}
: Consider $ \Theta( \mathbf{x} ) = \mathbb{E}[ \Psi( \mathbf{Y} ) \midil \mathbf{X} = \mathbf{x} ] \in \mathcal{B} $, where $ \mathcal{B} $ is a type 2 Banach space, and $ \Psi( \cdot ) $ is a function from $ \mathcal{R} $ to $ \mathcal{B} $. The estimate $ \widehat{\Theta}_n( \mathbf{x} ) $ of $ \Theta( \mathbf{x} ) $ is
\begin{align*}
\widehat{\Theta}_n( \mathbf{x} ) = \frac{\sum_{i=1}^{n} \Psi( \mathbf{Y}_i ) K( h^{-1} d( \mathbf{x}, \mathbf{X}_i ) )}{\sum_{i=1}^{n} K( h^{-1} d( \mathbf{x}, \mathbf{X}_i ) )} .
\end{align*}

Some examples of $ \Psi( \cdot ) $ and the resulting $ \Theta( \mathbf{x} ) $ are the following. Let $ \mathbf{Y} \in \mathbb{R} $, and $ \Psi( \mathbf{Y} ) = \mathbb{I}( \mathbf{Y} \le y ) $, where $ y \in \mathbb{R} $. Then, $ \Theta( \mathbf{x} ) $ is the conditional distribution of $ \mathbf{Y} $ given $ \mathbf{X} = \mathbf{x} $ at $ y $ (see \cite{ferraty2006estimating}, \cite{chagny2014adaptive}, etc.). Alternatively, if $ \Psi( \mathbf{Y} ) = \mathbf{Y}^r $, $ \Theta( \mathbf{x} ) $ is the conditional $ r $th moment of $ \mathbf{Y} $ given $ \mathbf{X} = \mathbf{x} $. Next, let $ \mathbf{Y} $ be a random vector in $ \mathbb{R}^q $. For $ \mathbf{u}, \mathbf{v} \in \mathbb{R}^q $ with $ \mathbf{u} = [ u_1, \ldots, u_q ] $ and $ \mathbf{v} = [ v_1, \ldots, v_q ] $, let $ \mathbf{u} \le \mathbf{v} $ denote $ u_i \le v_i $ for $ i = 1, \ldots, q $. Then, for $ \Psi( \mathbf{Y} ) = \mathbb{I}( \mathbf{Y} \le \mathbf{y} ) $, where $ \mathbf{y} \in \mathbb{R}^q $, $ \Theta( \mathbf{x} ) $ becomes the conditional multivariate distribution of $ \mathbf{Y} $ at $ \mathbf{y} $ given $ \mathbf{X} = \mathbf{x} $. When $ \mathbf{Y} $ is a univariate or multivariate random variable, the choice $ \Psi( \mathbf{Y} ) = \mathbf{Y} $ corresponds to the conditional mean of a univariate or multivariate response given $ \mathbf{X} = \mathbf{x} $ (see, e.g., \cite{ferraty2006nonparametric}, \cite{ferraty2007nonparametric}, etc.). Similarly, when $ \mathbf{Y} \in \mathcal{B} $ and $ \mathcal{B} $ is a separable Hilbert space, the choices $ \Psi( \mathbf{Y} ) = \mathbf{Y} $ and $ \Psi( \mathbf{Y} ) = \mathbf{Y} \otimes \mathbf{Y} $ (the outer product of $ \mathbf{Y} $ with itself) correspond to the conditional mean and the second conditional moment of $ \mathbf{Y} $ given $ \mathbf{X} = \mathbf{x} $, respectively (see \cite{ferraty2012regression}).
Note that when $ \mathcal{B} = \mathbb{R}^q $, $ \mathbf{Y} \otimes \mathbf{Y} $ becomes the $ q \times q $ matrix $ \mathbf{Y} \mathbf{Y}^t $.
\end{example}

\begin{example}[\textbf{Functions of conditional mean}] \label{ex:2.2}
: Let $ \mathcal{B}_1, \mathcal{B}_2 $ be two separable Banach spaces, $ \mathcal{U} $ be an open subset of $ \mathcal{B}_1 $, and $ \Psi( \cdot ) : \mathcal{R} \longrightarrow \mathcal{B}_1 $ be such that $ \mathbb{E}[ \Psi( \mathbf{Y} ) \midil \mathbf{X} = \mathbf{x} ] \in \mathcal{U} $. For $ \Gamma( \cdot ) : \mathcal{U} \longrightarrow \mathcal{B}_2 $, we consider $ \Theta( \mathbf{x} ) = \Gamma( \mathbb{E}[ \Psi( \mathbf{Y} ) \midil \mathbf{X} = \mathbf{x} ] ) $. Here, the kernel regression estimate $ \widehat{\Theta}_n( \mathbf{x} ) $ is
\begin{align*}
\widehat{\Theta}_n( \mathbf{x} ) = \Gamma\left( \frac{\sum_{i=1}^{n} \Psi( \mathbf{Y}_i ) K( h^{-1} d( \mathbf{x}, \mathbf{X}_i ) ) }{\sum_{i=1}^{n} K( h^{-1} d( \mathbf{x}, \mathbf{X}_i ) )} \right) .
\end{align*}

As a special case, let $ \mathbf{Y} $ be a real random variable. Let $ \Psi( \mathbf{Y} ) = ( \mathbf{Y}^2, \mathbf{Y} ) $ and $ \mathcal{U} = \{ ( u, v ) \in \mathbb{R}^2 \midil u > v^2,\, v > 0 \} $. Let $ \Gamma( \cdot ) : \mathcal{U} \longrightarrow \mathbb{R} $ be defined by $ \Gamma( u, v ) = v^{-1} \sqrt{u - v^2} $. Then, $ \Theta( \mathbf{x} ) = \Gamma( \mathbb{E}[ \Psi( \mathbf{Y} ) \midil \mathbf{X} = \mathbf{x} ] ) $ is the conditional coefficient of variation of $ \mathbf{Y} $ given $ \mathbf{X} = \mathbf{x} $ (see, e.g., \cite{dette2009testing}, \cite{dette2012testing}).

As another special case, let $ \mathbf{Y} = ( Y_1, Y_2 ) $ be a bivariate random variable. Let $ \Psi( \mathbf{Y} ) = \Psi( Y_1, Y_2 ) = ( Y_1 Y_2, Y_1, Y_2, Y_1^2, Y_2^2 ) $ and $ \mathcal{U} = \{ ( s, t, u, v, w ) \in \mathbb{R}^5 \midil v > t^2, w > u^2 \} $. Let $ \Gamma( \cdot ) : \mathcal{U} \longrightarrow \mathbb{R} $ be defined by
\begin{align*}
\Gamma( s, t, u, v, w ) = \frac{s - t u}{\sqrt{v - t^2} \sqrt{w - u^2}} .
\end{align*}
Then, $ \Theta( \mathbf{x} ) = \Gamma( \mathbb{E}[ \Psi( \mathbf{Y} ) \midil \mathbf{X} = \mathbf{x} ] ) $ is the conditional correlation coefficient of $ Y_1 $ and $ Y_2 $ given $ \mathbf{X} = \mathbf{x} $ (see, e.g., \citet[p.~13]{klemela2014multivariate}).

As the third special case, let the response space $ \mathcal{R} $ be a separable Hilbert space, and $ \mathcal{B}_2 $ denote the space of Hilbert--Schmidt operators on $ \mathcal{R} $. Note that $ \mathcal{B}_2 $ is a Hilbert space \citep[p.~195]{bhatia2009notes}. Also, the space of finite rank operators is dense in the space of Hilbert--Schmidt operators \citep[p.~196]{bhatia2009notes}, and the space of finite rank operators on a separable Hilbert space is itself separable. Consequently, $ \mathcal{B}_2 $ is a separable Hilbert space. Set $ \mathcal{B}_1 = \mathcal{B}_2 \times \mathcal{R} $. Define $ \Psi( \mathbf{Y} ) = ( \mathbf{Y} \otimes \mathbf{Y}, \mathbf{Y} ) $, $ \mathcal{U} = \mathcal{B}_1 $ and $ \Gamma( \mathbf{u}, \mathbf{v} ) = \mathbf{u} - \mathbf{v} \otimes \mathbf{v} $. Then, $ \Theta( \mathbf{x} ) = \mathbb{COV}[ \mathbf{Y} \midil \mathbf{X} = \mathbf{x} ] $, which is the conditional covariance operator of $ \mathbf{Y} $ given $ \mathbf{X} = \mathbf{x} $ (see, e.g., \cite{ferraty2012regression}). Note that when $ \mathbf{Y} $ is real random variable, this choice of $ \Psi( \mathbf{Y} ) $ and $ \Gamma( \cdot, \cdot ) $ corresponds to the conditional variance of $ \mathbf{Y} $ given $ \mathbf{X} = \mathbf{x} $.
\end{example}

\begin{example}[\textbf{Maximum likelihood regression}] \label{ex:2.3}
: Nonparametric estimation in a likelihood based regression problem with finite dimensional covariate was investigated in \cite{staniswalis1989kernel}, \cite{chaudhuri1995likelihood} and \cite{aerts1997local}. Let the covariate $ \mathbf{X} $ and the response $ \mathbf{Y} $ be random elements in the complete separable metric spaces $ \mathcal{C} $ and $ \mathcal{R} $, respectively. Suppose $ \mathbf{Y} $ given $ \mathbf{X} $ has a conditional density with respect to some sigma-finite measure in $ \mathcal{R} $, and it is given by $ f( \cdot \midil \Theta( \mathbf{x} ) ) $ for $ \mathbf{X} = \mathbf{x} $, where $ \Theta( \cdot ) : \mathcal{C} \longrightarrow \mathbb{R}^q $. We assume that the form of the function $ f( \cdot \midil \cdot ) $ is known, but $ \Theta( \cdot ) $ is unknown. We are interested in estimating $ \Theta( \mathbf{x} ) $ using maximum weighted likelihood procedure, where $ \mathbf{x} \in \mathcal{C} $ is fixed. The kernel estimate $ \widehat{\Theta}_n( \mathbf{x} ) $ of $ \Theta( \mathbf{x} ) $ is given by
\begin{align}
\widehat{\Theta}_n( \mathbf{x} ) = \arg \max_{ \mathbf{t} \in \mathbb{R}^q } \prod\limits_{i=1}^{n} [ f( \mathbf{Y}_i \midil \mathbf{t} ) ]^{W_{i,n}( \mathbf{x} )} , \text{ where } W_{i,n}( \mathbf{x} ) = \frac{K( h^{-1} d( \mathbf{x}, \mathbf{X}_i ) )}{\sum_{i=1}^{n} K( h^{-1} d( \mathbf{x}, \mathbf{X}_i ) )} .
\label{ex2.3eq1}
\end{align}
So, when $ f( \mathbf{y} \midil \mathbf{t} ) $ is a differentiable function of $ \mathbf{t} \in \mathbb{R}^q $, $ \widehat{\Theta}_n( \mathbf{x} ) $ is the solution (in $ \mathbf{t} $) of the likelihood equation
\begin{align}
\sum_{i=1}^{n} \left[ \nabla( \log f( \mathbf{Y}_i \midil \mathbf{t} ) ) \right] W_{i,n}( \mathbf{x} ) = \mathbf{0} .
\label{ex2.3eq1a}
\end{align}
Here $ \nabla $ denotes the gradient vector of first partial derivatives with respect to $ \mathbf{t} $.
\end{example}

It is well known that when the covariate $ \mathbf{X} $ is finite dimensional, say $ \mathbf{X} \in \mathbb{R}^q $, and $ \mathbf{X} $ has a continuous positive density at $ \mathbf{x} $, one needs to have a sequence of bandwidths $ \{ h_n \} $ such that $ h_n \longrightarrow 0 $ and $ n h_n^q \longrightarrow \infty $ as $ n \longrightarrow \infty $ to ensure the consistency of the kernel regression estimate $ \widehat{\Theta}_n( \mathbf{x} ) $ (see, e.g., chapter 3 in \cite{hardle1990applied}). To deal with covariates, which are not necessarily finite dimensional, define $ \phi( \mathbf{z}, h ) = \mathbb{P}[ d( \mathbf{z}, \mathbf{X} ) \le h ] $. The function $ \phi( \mathbf{z}, h ) $ is known as the small ball probability function, and it plays an important role in the asymptotic properties of nonparametric regression estimates. We make the following assumptions on the kernel and the sequence of bandwidths, which are required to establish the consistency of the estimates and derive their convergence rates.
\begin{enumerate}[label=A(\roman*), ref=A(\roman*)]
\item The kernel $ K( \cdot ) $ is supported on $ [ 0, 1 ] $ with $ K( u ) $ being bounded and bounded away from 0 for $ 0 \le u \le 1 $, i.e., there are constants $ 0 < l \le L $ such that $ l \le K( u ) \le L $ for all $ 0 \le u \le 1 $.
\label{assume:a1}

\item The bandwidth $ h_n \longrightarrow 0 $, and $ n \phi( \mathbf{x}, h_n ) \longrightarrow \infty $ as $ n \longrightarrow \infty $.
\label{assume:a2}
\end{enumerate}
The choice of the kernel $ K( \cdot ) $ described in assumption \ref{assume:a1} is equivalent to the type I kernel described in \citet[p.~42]{ferraty2006nonparametric}. This is a popular choice of kernel in the literature on nonparametric regression involving functional covariates (see, e.g., \cite{ferraty2006estimating}, \cite{burba2009k}, \cite{chagny2014adaptive,chagny2016adaptive}, etc.).
Note that for $ \mathbf{X} \in \mathbb{R}^q $ having a continuous positive density at $ \mathbf{x} $, the condition $ n \phi( \mathbf{x}, h_n ) \longrightarrow \infty $ as $ n \longrightarrow \infty $ in assumption \ref{assume:a2} is equivalent to $ n h_n^q \longrightarrow \infty $ as $ n \longrightarrow \infty $. Assumption \ref{assume:a2} is required to ensure the consistency of the kernel estimates involving an infinite dimensional covariate, and is also used in earlier works (see, e.g., \cite{ferraty2007nonparametric}).

\section{Bias-variance decomposition} \label{sec:3}
Let $ \mathcal{B} $ be a separable type 2 Banach space with a Schauder basis, and $ \Theta( \cdot ) : \mathcal{C} \longrightarrow \mathcal{B} $. For $ \mathbf{x} \in \mathcal{C} $, we consider the class of kernel regression estimates, which satisfy
\begin{align} \label{eq:main}
\widehat{\Theta}_n( \mathbf{x} ) - \Theta( \mathbf{x} ) = B_n( \mathbf{x} ) + V_n( \mathbf{x} ) + R_n( \mathbf{x} ) ,
\end{align}
where
\begin{align}
& B_n( \mathbf{x} ) = \mathbb{L}_\mathbf{x}\left( \frac{ \sum_{i=1}^{n} F( \mathbf{X}_i ) K( h_n^{-1} d( \mathbf{x}, \mathbf{X}_i ) ) }{ \sum_{i=1}^{n} K( h_n^{-1} d( \mathbf{x}, \mathbf{X}_i ) ) } - F( \mathbf{x} ) \right) ,
\label{eq:bias} \\
& V_n( \mathbf{x} ) = \mathbb{L}_\mathbf{x}\left( \frac{ \sum_{i=1}^{n} [ G( \mathbf{Y}_i ) - \mathbb{E}[ G( \mathbf{Y}_i ) \midil \mathbf{X}_i ] ] K( h_n^{-1} d( \mathbf{x}, \mathbf{X}_i ) ) }{ \sum_{i=1}^{n} K( h_n^{-1} d( \mathbf{x}, \mathbf{X}_i ) ) } \right) 
\label{eq:variance} .
\end{align}
Here, $ F( \cdot ) : \mathcal{C} \longrightarrow \mathcal{G} $, $ G( \cdot ) : \mathcal{R} \longrightarrow \mathcal{G} $, $ \mathbb{L}_\mathbf{x}( \cdot ) : \mathcal{G} \longrightarrow \mathcal{B} $ and $ \mathcal{G} $ is a separable Banach space. $ \mathbb{L}_\mathbf{x}( \cdot ) $ is a continuous linear map. The functions $ F( \cdot ) $, $ G( \cdot ) $ and the remainder term $ R_n( \mathbf{x} ) $ are assumed to satisfy the following conditions.
\begin{enumerate}[label=B(\roman*), ref=B(\roman*)]
\item \label{cond:1}
Let $ \beta > 0 $ be a constant. Then, $ F( \cdot ) \in \mathcal{F}( \mathbf{x}, \beta, \mathcal{G} ) $. Here, $ \mathcal{F}( \mathbf{x}, \beta, \mathcal{G} ) $ is a class of functions $ F( \cdot ) : \mathcal{C} \longrightarrow \mathcal{G} $ such that for some constant $ b_F > 0 $, $ \| F( \mathbf{z} ) - F( \mathbf{x} ) \| \le b_F d( \mathbf{x}, \mathbf{z} )^{\beta} $ for all $ \mathbf{z} $ lying in a neighborhood of $ \mathbf{x} $.

\item \label{cond:3}
$ G( \cdot ) $ is such that for some $ \nu > 2 $, $ \mathbb{E}[ \| G( \mathbf{Y} ) - \mathbb{E}[ G( \mathbf{Y} ) \midil \mathbf{X} = \mathbf{z} ] \|^\nu \midil \mathbf{X} = \mathbf{z} ] $ is uniformly bounded for $ \mathbf{z} $ lying in a neighborhood of $ \mathbf{x} $.

\item \label{cond:4}
$ R_n( \mathbf{x} ) = o_\mathbb{P}( \delta_n ) $ as $ n \longrightarrow \infty $, where $ \delta_n = \max\big\{ h_n^\beta, \big[ n \phi( \mathbf{x}, h_n ) \big]^{-1/2} \big\} $.
\end{enumerate}

Note that $ \mathbb{E}[ V_n( \mathbf{x} ) \midil \mathbf{X}_1, \ldots, \mathbf{X}_n ] = \mathbf{0} $. We can view $ B_n( \mathbf{x} ) $ as the bias term and $ V_n( \mathbf{x} ) $ as the variance term in kernel regression. Note that condition \ref{cond:1} is related to the smoothness of the regression function. Condition \ref{cond:3} imposes a bound on the variability of the residual of the regression. Condition \ref{cond:4} essentially states that the remainder term in our bias--variance decomposition is asymptotically negligible. We shall now verify the validity of the above bias--variance decomposition in Examples \ref{ex:2.1}, \ref{ex:2.2} and \ref{ex:2.3} in \autoref{sec:2}.

\textit{Example 1 (continued).}
Recall \autoref{ex:2.1} considered in \autoref{sec:2}. We can set
\begin{align*}
& B_n( \mathbf{x} ) = \frac{ \sum_{i=1}^{n} \mathbb{E}[ \Psi( \mathbf{Y}_i ) \midil \mathbf{X}_i ] K( h_n^{-1} d( \mathbf{x}, \mathbf{X}_i ) ) }{ \sum_{i=1}^{n} K( h_n^{-1} d( \mathbf{x}, \mathbf{X}_i ) ) } - \mathbb{E}[ \Psi( \mathbf{Y} ) \midil \mathbf{X} = \mathbf{x} ] , \\
& V_n( \mathbf{x} ) = \frac{ \sum_{i=1}^{n} [ \Psi( \mathbf{Y}_i ) - \mathbb{E}[ \Psi( \mathbf{Y}_i ) \midil \mathbf{X}_i ] ] K( h_n^{-1} d( \mathbf{x}, \mathbf{X}_i ) ) }{ \sum_{i=1}^{n} K( h_n^{-1} d( \mathbf{x}, \mathbf{X}_i ) ) } , \\	
& \text{and } R_n( \mathbf{x} ) = \mathbf{0} .
\end{align*}
Then, setting $ G( \mathbf{Y} ) = \Psi( \mathbf{Y} ) $, $ F( \mathbf{X} ) = \mathbb{E}[ \Psi( \mathbf{Y} ) \midil \mathbf{X} ] $ and $ \mathbb{L}_\mathbf{x}( \cdot ) $ to be the identity map on $ \mathcal{B} $, \eqref{eq:main} holds for any kernel satisfying \ref{assume:a1} and any sequence of bandwidths $ \{ h_n \} $ satisfying \ref{assume:a2}. Here, condition \ref{cond:4} is trivially satisfied. Note that in this case, $ F( \mathbf{z} ) = \Theta( \mathbf{z} ) $, and so condition \ref{cond:1} is satisfied when $ \Theta( \mathbf{z} ) $ is Holder continuous at $ \mathbf{x} $ with exponent $ \beta $, and the class $ \mathcal{F}( \mathbf{x}, \beta, \mathcal{B} ) $ can be taken as the class of all Holder continuous functions at $ \mathbf{x} $.
Condition \ref{cond:3} is satisfied when for some $ \nu > 2 $, $ \mathbb{E}[ \| \Psi( \mathbf{Y} ) - \mathbb{E}[ \Psi( \mathbf{Y} ) \midil \mathbf{X} = \mathbf{z} ] \|^\nu \midil \mathbf{X} = \mathbf{z} ] $ is uniformly bounded for $ \mathbf{z} $ lying in a neighborhood of $ \mathbf{x} $. In particular, \ref{cond:3} holds for the location-scale type model $ \Psi( \mathbf{Y} ) = l( \mathbf{X} ) + s( \mathbf{X} ) \mathbf{U} $, where $ l( \cdot ) : \mathcal{C} \longrightarrow \mathcal{B} $ and $ s( \cdot ) : \mathcal{C} \longrightarrow ( 0, \infty ) $ are continuous functions, and $ \mathbf{U} $ is a zero-mean random element in $ \mathcal{B} $, which is independent of $ \mathbf{X} $ with $ \mathbb{E}[ \| \mathbf{U} \|^\nu ] < \infty $ for some $ \nu > 2 $.

\textit{Example 2 (continued).}
Consider again the class of estimators described in \autoref{ex:2.2}. The following proposition asserts that the bias--variance decomposition \eqref{eq:main} holds for those estimators.
\begin{theorem} \label{ex:prop1}
In \autoref{ex:2.2} considered in \autoref{sec:2}, let $ \mathcal{B}_2 $ be a type 2 Banach space. Let the kernel function $ K( \cdot ) $ satisfy \ref{assume:a1} and the bandwidths $ \{ h_n \} $ satisfy \ref{assume:a2}.
Assume that $ \Gamma( \cdot ) $ is Fr\'{e}chet differentiable with derivative $ \Gamma'( \cdot ) $. Let $ \mathbb{L}_\mathbf{x}( \cdot ) = \Gamma'\left( \mathbb{E}[ \Psi( \mathbf{Y} ) \midil \mathbf{X} = \mathbf{x} ] \right)( \cdot ) $, $ G( \mathbf{Y} ) = \Psi( \mathbf{Y} ) $, $ F( \mathbf{z} ) = \mathbb{E}[ G( \mathbf{Y} ) \midil \mathbf{X} = \mathbf{z} ] $, and conditions \ref{cond:1} and \ref{cond:3} hold. Then, \ref{cond:4} is also satisfied, and consequently the bias--variance decomposition in \eqref{eq:main} holds.
\end{theorem}
Note that in all the specific cases discussed in \autoref{ex:2.2}, namely, the coefficient of variation, the correlation coefficient and the covariance operator, the function $ \Gamma( \cdot ) $ satisfies the differentiability condition stated in \autoref{ex:prop1}, and its derivative can be computed in a straight forward way.

When $ \mathbf{Y} $ is a real-valued random variable and $ \Theta( \mathbf{x} ) $ is the conditional coefficient of variation of $ \mathbf{Y} $ given $ \mathbf{X} = \mathbf{x} $ as mentioned in \autoref{ex:2.2}, the assumption \ref{cond:1} is satisfied if $ \mathbb{E}[ \mathbf{Y}^2 \midil \mathbf{X} = \mathbf{z} ] $ and $ \mathbb{E}[ \mathbf{Y} \midil \mathbf{X} = \mathbf{z} ] $ are both Holder continuous at $ \mathbf{x} $ with exponent $ \beta $. Further, the assumption \ref{cond:3} is satisfied if $ \mathbb{E}[ \mathbf{Y}^{4+\alpha} \midil \mathbf{X} = \mathbf{z} ] $ is uniformly bounded for $ \mathbf{z} $ lying in a neighborhood of $ \mathbf{x} $ for some $ \alpha > 0 $. Note that conditions $ \mathbb{E}[ \mathbf{Y} \midil \mathbf{X} = \mathbf{x} ] > 0 $ and $ \mathbb{V}[ \mathbf{Y} \midil \mathbf{X} = \mathbf{x} ] > 0 $ ensure that $ \mathbb{E}[ \Psi( \mathbf{Y} ) \midil \mathbf{X} = \mathbf{x} ] $ lies in the domain of $ \Gamma( \cdot ) $.

When $ \mathbf{Y} = ( Y_1, Y_2 ) $ is a bivariate random variable, and $ \Theta( \mathbf{x} ) $ is the conditional correlation between $ Y_1 $ and $ Y_2 $ given $ \mathbf{X} = \mathbf{x} $ as mentioned in \autoref{ex:2.2}, assumption \ref{cond:1} is satisfied if each of $ \mathbb{E}[ Y_1 Y_2 \midil \mathbf{X} = \mathbf{z} ] $, $ \mathbb{E}[ Y_1 \midil \mathbf{X} = \mathbf{z} ] $, $ \mathbb{E}[ Y_2 \midil \mathbf{X} = \mathbf{z} ] $, $ \mathbb{E}[ Y_1^2 \midil \mathbf{X} = \mathbf{z} ] $ and $ \mathbb{E}[ Y_2^2 \midil \mathbf{X} = \mathbf{z} ] $ is Holder continuous at $ \mathbf{x} $ with exponent $ \beta $. Assumption \ref{cond:3} is satisfied if $ \mathbb{E}[ \| \mathbf{Y} \|^{4+\alpha} \midil \mathbf{X} = \mathbf{z} ] $ is uniformly bounded for $ \mathbf{z} $ lying in a neighborhood of $ \mathbf{x} $ for some $ \alpha > 0 $. Further, $ \mathbb{V}[ Y_1 \midil \mathbf{X} = \mathbf{x} ] , \mathbb{V}[ Y_2 \midil \mathbf{X} = \mathbf{x} ] > 0 $ ensure that $ \mathbb{E}[ \Psi( \mathbf{Y} ) \midil \mathbf{X} = \mathbf{x} ] $ lies in the domain of $ \Gamma( \cdot ) $.

One can verify that when $ \Theta( \mathbf{x} ) $ is the conditional covariance of $ \mathbf{Y} $ given $ \mathbf{X} = \mathbf{x} $, assumption \ref{cond:1} is satisfied if $ \mathbb{E}[ \mathbf{Y} \otimes \mathbf{Y} \midil \mathbf{X} = \mathbf{z} ] $ and $ \mathbb{E}[ \mathbf{Y} \midil \mathbf{X} = \mathbf{z} ] $ are both Holder continuous at $ \mathbf{x} $ with exponent $ \beta $. Assumption \ref{cond:3} holds if $ \mathbb{E}[ \| \mathbf{Y} \|^{4+\alpha} \midil \mathbf{X} = \mathbf{z} ] $ is uniformly bounded for $ \mathbf{z} $ lying in a neighborhood of $ \mathbf{x} $ for some $ \alpha > 0 $.

\textit{Example 3 (continued).}
In the case of \autoref{ex:2.3} considered in \autoref{sec:2}, define $ g( \mathbf{y} \midil \mathbf{t} ) = \log f( \mathbf{y} \midil \mathbf{t} ) $, where $ \mathbf{t} \in \mathbb{R}^q $. Let $ \mathcal{T} $ be an open ball in $ \mathbb{R}^q $ containing the range of $ \Theta( \cdot ) $. We now assume some Cramer-type regularity conditions on the log-likelihood $ g( \mathbf{y} \midil \mathbf{t} ) $ that are required for asymptotic analysis of weighted maximum likelihood estimates (see, e.g., \cite{chaudhuri1995likelihood}). The support of $ f( \mathbf{y} \midil \mathbf{t} ) $ is assumed to be same for all $ \mathbf{t} \in \mathcal{T} $, and $ g( \mathbf{y} \midil \mathbf{t} ) $ is assumed to be thrice continuously differentiable with respect to $ \mathbf{t} $ for $ \mathbf{t} \in \mathcal{T} $. Denote the Hessian matrix of all second order partial derivatives of $ g( \mathbf{y} \midil \mathbf{t} ) $ with respect to $ \mathbf{t} $ as $ \Delta_2( g( \mathbf{y} \midil \mathbf{t} ) ) $, and the array of all third order partial derivatives of $ g( \mathbf{y} \midil \mathbf{t} ) $ with respect to $ \mathbf{t} $ as $ \Delta_3( g( \mathbf{y} \midil \mathbf{t} ) ) $. Define $ \mathbf{I}( \Theta( \mathbf{z} ) ) = - \mathbb{E}[ \Delta_2( g( \mathbf{Y} \midil \Theta( \mathbf{z} ) ) ) \midil \mathbf{X} = \mathbf{z} ] $, and assume that $ \mathbf{I}( \Theta( \mathbf{z} ) ) $ is finite, positive definite and continuous for $ \mathbf{z} $ lying in a neighborhood of $ \mathbf{x} $. Also, assume that for $ \mathbf{t} \in \mathcal{T} $, there exist two nonnegative random variables $ \mathbf{D}_1( \mathbf{Y} \midil \mathbf{t} ) $ and $ \mathbf{D}_2( \mathbf{Y} \midil \mathbf{t} ) $ such that $ \mathbb{E}[ \mathbf{D}_1( \mathbf{Y} \midil \mathbf{t} ) ]^2 < \infty $, $ \mathbb{E}[ \mathbf{D}_2( \mathbf{Y} \midil \mathbf{t} ) ] < \infty $, and $ \| \Delta_2( g( \mathbf{Y} \midil \mathbf{s} ) ) \| \le \mathbf{D}_1( \mathbf{Y} \midil \mathbf{t} ) $, $ \| \Delta_3( g( \mathbf{Y} \midil \mathbf{s} ) ) \| \le \mathbf{D}_2( \mathbf{Y} \midil \mathbf{t} ) $ for any $ \mathbf{s} $ in some neighborhood of $ \mathbf{t} $ contained in $ \mathcal{T} $.

In the next proposition, we see that the decomposition \eqref{eq:main} along with conditions \ref{cond:1}--\ref{cond:4} is satisfied for the weighted maximum likelihood estimate $ \widehat{\Theta}_n( \mathbf{x} ) $ defined in \eqref{ex2.3eq1}.
\begin{theorem} \label{ex:prop2}
In \autoref{ex:2.3} considered in \autoref{sec:2},
assume that $ \Theta( \cdot ) \in \mathcal{F}( \mathbf{x}, \beta, \mathbb{R}^q ) $ for some $ \beta > 0 $, where $ \mathcal{F}( \mathbf{x}, \beta, \mathbb{R}^q ) $ is as defined in \ref{cond:1}. Let the kernel function $ K( \cdot ) $ satisfy \ref{assume:a1} and the bandwidths $ \{ h_n \} $ satisfy \ref{assume:a2}. Then, under the Cramer type regularity conditions stated above, the decomposition \eqref{eq:main} along with conditions \ref{cond:1}--\ref{cond:4} will hold for $ \widehat{\Theta}_n( \mathbf{x} ) $ in \eqref{ex2.3eq1} if we choose $ \mathbb{L}_\mathbf{x}( \cdot ) = [ \mathbf{I}( \Theta( \mathbf{x} ) ) ]^{-1}( \cdot ) $, $ G( \mathbf{Y} ) = \nabla g( \mathbf{Y} \midil \Theta( \mathbf{X} ) ) $ and $ F( \mathbf{X} ) = \mathbf{I}( \Theta( \mathbf{x} ) )( \Theta( \mathbf{X} ) ) $, where $ g( \mathbf{y} \midil \mathbf{t} ) = \log f( \mathbf{y} \midil \mathbf{t} ) $.
\end{theorem}

\subsection{Asymptotic behavior of the bias and the variance} \label{subsec:3_1}
In this subsection, the orders of convergence of the bias term $ B_n( \mathbf{x} ) $ and the variance term $ V_n( \mathbf{x} ) $ in \eqref{eq:main} are investigated. It follows from assumptions \ref{assume:a2} and \ref{cond:1} that $ \| B_n( \mathbf{x} ) \| \le \| \mathbb{L}_\mathbf{x} \| b_F h_n^\beta $ for all sufficiently large $ n $. So, for all choices of bandwidths $ \{ h_n \} $ with $ h_n \longrightarrow 0^+ $ as $ n \longrightarrow \infty $,
\begin{align} \label{eq:b2}
\mathbb{E}\left[ \| B_n( \mathbf{x} ) \|^2 \right] \le ( \| \mathbb{L}_\mathbf{x} \| b_F )^2 h_n^{2 \beta}
\end{align}
for all sufficiently large $ n $.
The inequality \eqref{eq:b2} leads to an upper bound of the rate of convergence of the bias term, and it will be used later to study the asymptotic properties of the estimate $ \widehat{\Theta}_n( \mathbf{x} ) $.

We next discuss the asymptotic behavior of the variance term $ V_n( \mathbf{x} ) $. We derive an upper bound of the convergence rate of $ \mathbb{E}[ \| V_n( \mathbf{x} ) \|^2 ] $ in the theorem below.
\begin{theorem} \label{thm:1}
Under \ref{assume:a1}, \ref{assume:a2} and \ref{cond:3}, $ n \phi( \mathbf{x}, h_n ) \mathbb{E}[ \| V_n( \mathbf{x} ) \|^2 ] = O( 1 ) $ as $ n \longrightarrow \infty $.
\end{theorem}
The following condition is needed to derive the asymptotic distribution of $ V_n( \mathbf{x} ) $.
\begin{enumerate}[label=B(\roman*), ref=B(\roman*)]
\setcounter{enumi}{3}
\item \label{cond:5}
$ \mathcal{B} $ is a separable Hilbert space, and $ G( \cdot ) $ in \eqref{eq:variance} is such that the covariance operator $ \mathbb{D}( \cdot, \cdot \mid \mathbf{z} ) : \mathcal{B} \times \mathcal{B} \longrightarrow \mathbb{R} $ defined by
$ \mathbb{D}( \mathbf{u}, \mathbf{v} \midil \mathbf{z} ) 
= \mathbb{E}[ 
\langle \mathbf{u}, \mathbb{L}_\mathbf{x}( G( \mathbf{Y} ) - \mathbb{E}[ G( \mathbf{Y} ) \midil \mathbf{X} = \mathbf{z} ] ) \rangle
\langle \mathbf{v}, \mathbb{L}_\mathbf{x}( G( \mathbf{Y} ) - \mathbb{E}[ G( \mathbf{Y} ) \midil \mathbf{X} = \mathbf{z} ] ) \rangle
\midil \mathbf{X} = \mathbf{z} ] $, where $ \mathbf{u}, \mathbf{v} \in \mathcal{B} $, converges to $ \mathbb{D}( \cdot, \cdot \mid \mathbf{x} ) $ in the trace norm as $ \mathbf{z} \longrightarrow \mathbf{x} $, and $ \mathbb{D}( \cdot, \cdot \mid \mathbf{x} ) $ is a bounded positive definite operator.
\end{enumerate}
Condition \ref{cond:5} is related to the smoothness of the conditional distribution of the residual in the regression given the covariate, and it holds in many common models. For example, consider the location-scale type model $ \mathbb{L}_\mathbf{x}( G( \mathbf{Y} ) ) = l( \mathbf{X} ) + s( \mathbf{X} ) \mathbf{U} $, where $ l( \cdot ) : \mathcal{C} \longrightarrow \mathcal{B} $ and $ s( \cdot ) : \mathcal{C} \longrightarrow ( 0, \infty ) $ are continuous functions, and $ \mathbf{U} $ is a zero-mean random element in $ \mathcal{B} $, which is independent of $ \mathbf{X} $, having a bounded positive definite covariance operator.

From assumption \ref{assume:a1}, it follows that
$ l^j \phi( \mathbf{x}, h ) \le \mathbb{E}\left[ K^j( h^{-1} d( \mathbf{x}, \mathbf{X} ) ) \right] \le L^j \phi( \mathbf{x}, h ) $
for any positive integer $ j $ and any bandwidth $ h > 0 $. Define
\begin{align}
E_n^{(j)}( \mathbf{x} ) = \left[ \phi( \mathbf{x}, h_n ) \right]^{-1} \mathbb{E}\left[ K^j( h_n^{-1} d( \mathbf{x}, \mathbf{X} ) ) \right]
\label{eq:assume:a1}
\end{align}
for all positive integer $ j $. Note that
\begin{align}
0 < L^{-1} l \le \left[ E_n^{(2)}( \mathbf{x} ) \right]^{-1/2} E_n^{(1)}( \mathbf{x} ) \allowbreak \le l^{-1} L < \infty
\label{eq:assume:a2}
\end{align}
for all $ n $. In the next theorem, we establish the asymptotic Gaussianity of $ V_n( \mathbf{x} ) $.
\begin{theorem} \label{thm:2}
Let the kernel function $ K( \cdot ) $ satisfy \ref{assume:a1}, and the sequence of bandwidths $ \{ h_n \} $ satisfy \ref{assume:a2}. Then, under conditions \ref{cond:3} and \ref{cond:5},
\begin{align*}
\left[ n \phi( \mathbf{x}, h_n ) \right]^{1/2} \left[ E_n^{(2)}( \mathbf{x} ) \right]^{-1/2} \allowbreak E_n^{(1)}( \mathbf{x} ) V_n( \mathbf{x} ) \longrightarrow \mathbf{W}
\end{align*}
\textit{in distribution} as $ n \longrightarrow \infty $, where $ \mathbf{W} $ is a zero mean Gaussian random element in $ \mathcal{B} $ with covariance operator $ \mathbb{D}( \cdot, \cdot \midil \mathbf{x} ) $.
\end{theorem}
Recall that $ R_n( \mathbf{x} ) = \mathbf{0} $ for the mean-type regression problems described in \autoref{ex:2.1}. So, for these class of regression problems, from \autoref{thm:2} and \eqref{eq:main} we get that
\begin{align*}
\left[ n \phi( \mathbf{x}, h_n ) \right]^{1/2} \left[ E_n^{(2)}( \mathbf{x} ) \right]^{-1/2} \allowbreak E_n^{(1)}( \mathbf{x} ) \left[ \widehat{\Theta}_n( \mathbf{x} ) - \Theta( \mathbf{x} ) - B_n( \mathbf{x} ) \right] \longrightarrow \mathbf{W}
\end{align*}
\textit{in distribution} as $ n \longrightarrow \infty $.
Define
\begin{align*}
e_n[ G( \mathbf{Y} ) \midil \mathbf{x} ] = \frac{ \sum_{i=1}^{n} G( \mathbf{Y}_i ) K( h_n^{-1} d( \mathbf{x}, \mathbf{X}_i ) ) }{ \sum_{i=1}^{n} K( h_n^{-1} d( \mathbf{x}, \mathbf{X}_i ) ) } .
\end{align*}
The covariance operator $ \mathbb{D}( \cdot, \cdot \midil \mathbf{x} ) $ of $ \mathbf{W} $ may be estimated by
\begin{align*}
& \widehat{\mathbb{D}}_n( \mathbf{u}, \mathbf{v} \midil \mathbf{x} ) \\
& = \frac{ \sum_{i=1}^{n} \left[ \left\langle \mathbf{u}, \mathbb{L}_\mathbf{x}( G( \mathbf{Y}_i ) - e_n[ G( \mathbf{Y} ) \midil \mathbf{x} ] ) \right\rangle \left\langle \mathbf{v}, \mathbb{L}_\mathbf{x}( G( \mathbf{Y}_i ) - e_n[ G( \mathbf{Y} ) \midil \mathbf{x} ] ) \right\rangle \right] K\left( \frac{d( \mathbf{x}, \mathbf{X}_i )}{h_n} \right) }{ \sum_{i=1}^{n} K\left( \frac{d( \mathbf{x}, \mathbf{X}_i )}{h_n} \right) } .
\end{align*}

The function $ \phi( \mathbf{x}, h ) $ plays a central role in determining the convergence rate and the asymptotic distribution of $ V_n( \mathbf{x} ) $, and we discuss it in detail in the next subsection.

\subsection{The small ball probability function} \label{subsec:3_2}
When the covariate $ \mathbf{X} $ is finite dimensional, say $ \mathbf{X} \in \mathbb{R}^q $, and it has a continuous positive density at $ \mathbf{x} $, it follows that $ \phi( \mathbf{x}, h ) \sim h^q $ as $ h \longrightarrow 0^+ $. But if $ \mathbf{X} $ is a random element in an infinite dimensional space, getting the asymptotic order of $ \phi( \mathbf{x}, h ) $ as $ h \longrightarrow 0^+ $ is much more difficult, and the available results in this area are mostly for the case where $ \mathbf{X} $ is a Gaussian process (see, e.g., \cite{lifshits2013gaussian}). In the literature, the popular approach has been to first derive the limiting behavior of $ \log \phi( \mathbf{0}, h ) $ as $ h \longrightarrow 0^+ $, when $ \mathbf{X} $ is a Gaussian random element centered at $ \mathbf{0} $. Then, one makes a connection between $ \phi( \mathbf{x}, h ) $ and $ \phi( \mathbf{0}, h ) $ for suitable $ \mathbf{x} $ and sufficiently small $ h $.

The asymptotic behavior of $ \log \phi( \mathbf{0}, h ) $ was investigated in \cite{li2001small} for real-valued centered Gaussian Markov processes on $ [ 0, 1 ] $ under the $ L_p $-norm, where $ 1 \le p \le \infty $. It was shown there that in such a case, $ h^2 \log \phi( \mathbf{0}, h ) \longrightarrow - c_1 $ as $ h \longrightarrow 0^+ $, where $ c_1 > 0 $ is a constant depending on $ p $. For $ \mathbf{X} $ being a fractional Brownian motion on $ [ 0, 1 ] $ with Hurst index $ \gamma \in ( 0, 1 ) $, it was shown in Theorem 4.6 in \cite{li2001gaussian} that under the $ L_\infty $-norm, $ -c_2 h^{-1 / \gamma} \le \log \phi( \mathbf{0}, h ) \le -c_3 h^{-1 / \gamma} $ for all $ 0 < h \le 1 $. Here, $ c_2 $ and $ c_3 $ are positive constants depending on $ \gamma $. For $ \mathbf{X} $ being an integrated fractional Brownian motion with Hurst index $ \gamma \in ( 0, 1 ) $, it was established in Theorem 4.10 of \cite{li2001gaussian} that under the $ L_\infty $-norm, $ -c_4 h^{-1 / (1 + \gamma)} \le \log \phi( \mathbf{0}, h ) \le -c_5 h^{-1 / (1 + \gamma)} $ for all $ 0 < h \le 1 $, where $ c_4 $ and $ c_5 $ are positive constants depending on $ \gamma $.

For the L\'{e}vy fractional Brownian motion on $ [ 0, 1 ]^q $ with Hurst index $ \gamma \in ( 0, 1 ) $, it was proved in Theorem 5.1 in \cite{li2001gaussian} that under the $ L_\infty $-norm, $ -c_6 h^{- q / \gamma} \le \log \phi( \mathbf{0}, h ) \le -c_7 h^{- q / \gamma} $ for all $ 0 < h \le 1 $. Here, $ c_6 $ and $ c_7 $ are positive constants depending on $ \gamma $ and $ q $. For a Brownian sheet on $ [ 0, 1 ]^q $, it follows from Theorem 5.3 in \cite{li2001gaussian} that under the $ L_2 $-norm, $ -c_8 h^{-2} ( \log (1 / h) )^{(2 q - 2)} \allowbreak \le \log \phi( \mathbf{0}, h ) \le -c_9 h^{-2} ( \log (1 / h) )^{(2 q - 2)} $ as $ h \longrightarrow 0^+ $, where $ c_8, c_9 > 0 $ are constants depending on $ q $. It was shown in Theorem 5.4 in \cite{li2001gaussian} that if $ \mathbf{X} $ is a Brownian sheet on $ [ 0, 1 ]^2 $, we have $ -c_{10} h^{-2} ( \log (1 / h) )^3 \le \log \phi( \mathbf{0}, h ) \le -c_{11} h^{-2} ( \log (1 / h) )^3 $ under the $ L_\infty $-norm, where $ c_{10}, c_{11} > 0 $ are constants.

\subsection{Shifted small ball probability} \label{subsec:3_3}
As we have already mentioned, the asymptotic behavior of $ \log \phi( \mathbf{x}, h ) $ is derived by establishing some relationship between $ \phi( \mathbf{x}, h ) $ and $ \phi( \mathbf{0}, h ) $. As described in subsection 1.2 in \cite{mas2012lower}, one can establish a relation between $ \phi( \mathbf{x}, h ) $ and $ \phi( \mathbf{0}, h ) $ if the probability measure of $ \mathbf{X} - \mathbf{x} $ is absolutely continuous with respect to the probability measure of $ \mathbf{X} $, and the density of the measure of $ \mathbf{X} - \mathbf{x} $ with respect to the measure of $ \mathbf{X} $ is suitably smooth. This approach is motivated from the Cameron--Martin Theorem describing the Radon--Nikodym derivative of a Weiner measure translated by $ \mathbf{x} $ with respect to the centered Weiner measure, where $ \mathbf{x} $ is an element of the reproducing kernel Hilbert space associated with the centered Weiner measure (see \cite{cameron1944transformations}). When $ \mathbf{X} $ is a centered Gaussian random element in a separable Banach space, and $ \mathbf{x} $ is an element of the associated reproducing kernel Hilbert space, from Theorem 3.1 in \cite{li2001gaussian} we get that $ \exp[ - ( 1 / 2 ) \| \mathbf{x} \|_\mu^2 ] \phi( \mathbf{0}, h ) \le \phi( \mathbf{x}, h ) \le \phi( \mathbf{0}, h ) $ for all $ h > 0 $, where $ \| \cdot \|_\mu $ is the norm in the reproducing kernel Hilbert space. But this result is not very useful for our purpose since the probability of the event that an infinite dimensional Gaussian random element lies in its reproducing kernel Hilbert space is zero (see Corollary 7.1 in \cite{lukic2001stochastic}). Fortunately, it follows from Remark 2.2 in \cite{dereich2005probabilities} that when $ \mathbf{X} $ is a centered Gaussian random element in a separable Banach space, then for almost all $ \mathbf{x} $, $ ( \phi( \mathbf{0}, h/2 ) )^2 \le \phi( \mathbf{x}, h ) \le \phi( \mathbf{0}, h ) $ for all sufficiently small $ h $. How small $ h $ needs to be depends on that particular $ \mathbf{x} $. On the other hand, it follows from Theorem 2.1 in \cite{hoffmann1979lower} that for $ \mathbf{X} $ being a centered Gaussian random element in a separable infinite dimensional Hilbert space, we have $ \exp[ - ( 1 / 2 ) \| \mathbf{x} \|^2 ] \phi( \mathbf{0}, h ) \le \phi( \mathbf{x}, h ) \le \phi( \mathbf{0}, h ) $ for all $ h > 0 $.

Let $ \mathbf{X} $ be a centered Gaussian random element in a separable Hilbert space. The Karhunen--Loeve expansion of $ \mathbf{X} $ is $ \mathbf{X} = \sum_{j=1}^{\infty} \sqrt{\lambda_j} Z_j \boldsymbol{\psi}_j $, where $ \{ Z_j \} $ is a collection of independent normal random variables with mean 0 and variance 1, $ \{ \lambda_j \} $ is the sequence of decreasing eigenvalues of the covariance operator of $ \mathbf{X} $, and $ \{ \boldsymbol{\psi}_j \} $ is an orthonormal basis of the Hilbert space. Here, the small ball probability $ \phi( \mathbf{x}, h ) $ can be related to the rate of decrease of the sequence $ \{ \lambda_j \} $. As discussed in subsection 4.1 in \cite{chagny2014adaptive}, for certain rates of decrease for $ \{ \lambda_j \} $, e.g., if for some $ \alpha > 1 $, $ j^\alpha \lambda_j $ is bounded and bounded away from 0 for all $ j $, we may have $ c_{12} h^{p_1} \exp( -c_{13} h^{-q_1} ) \le \phi( \mathbf{x}, h ) \le c_{14} h^{p_2} \exp( -c_{15} h^{-q_1} ) $ for positive constants $ c_{12} $, $ c_{13} $, $ c_{14} $, $ c_{15} $, $ p_1 $, $ p_2 $ and $ q_1 $. Alternatively, for some other rates, e.g., if $ j \exp[ 2 j ] \lambda_j $ is bounded and bounded away from 0 for all $ j $, we may have $ c_{16} h^{p_3} \exp[ -c_{17} ( \log ( 1 / h ) )^{q_2} ] \le \phi( \mathbf{x}, h ) \le c_{18} h^{p_4} \exp[ -c_{17} ( \log ( 1 / h ) )^{q_2} ] $ for positive constants $ c_{16} $, $ c_{17} $, $ c_{18} $, $ p_3 $, $ p_4 $ and $ q_2 > 1 $ (see subsection 4.1 in \cite{chagny2014adaptive}). See also Theorem 4.4, Examples 4.5, 4.6 and 4.7 in \cite{hoffmann1979lower} for a discussion on the relation between the small ball probability $ \phi( \mathbf{x}, h ) $ and the rate of decrease of $ \{ \lambda_j \} $.

From the discussion on the small ball probability functions above, it is now clear that in a diverse collection of cases, we have
\begin{align}
C_1 h^{t_1} \exp\left[ -C_2 \left( \frac{1}{h} \right)^{t_2} \left( \log \frac{1}{h} \right)^{t_3} \right]
\le \phi( \mathbf{x}, h ) 
\le C_3 h^{t_4} \exp\left[ -C_4 \left( \frac{1}{h} \right)^{t_2} \left( \log \frac{1}{h} \right)^{t_3} \right]
\label{eq:1}
\end{align}
as $ h \longrightarrow 0^+ $. Here, $ C_1, C_2, C_3, C_4 > 0 $ and $ t_1, t_2, t_3, t_4 \ge 0 $ are appropriate constants, all of which, except $ C_1 $, are independent of $ \mathbf{x} $. $ C_1 $ may or may not depend on $ \mathbf{x} $, but if it depends on $ \mathbf{x} $ then $ C_1 = C_1' \exp[ - ( 1 / 2 ) \| \mathbf{x} \|^2 ] $ for some positive constant $ C_1' $. For infinite dimensional covariates, either $ t_2 > 0 $, or $ t_3 > 1 $ is an integer with $ C_2 = C_4 $.
Define
\begin{align}
m( h ) = C_2 ( 1 / h )^{t_2} ( \log ( 1 / h ) )^{t_3}
\label{eq:mh}
\end{align}
for $ 0 < h < 1 $. We shall derive the optimum convergence rates of the estimates in terms of $ m( h ) $.

The previous discussion of small ball probabilities are concerned with only Gaussian random elements. We next consider small ball probabilities of some infinite dimensional non-Gaussian distributions.
Let $ \mathcal{B}_1 $ and $ \mathcal{B}_2 $ be separable Banach spaces, and $ f( \cdot ) : \mathcal{B}_2 \longrightarrow \mathcal{B}_1 $ be a function such that for any $ \mathbf{u} \in \mathcal{B}_2 $, there exist constants $ r, s > 0 $, which may depend on $ \mathbf{u} $, such that for any $ \mathbf{v} \in \mathcal{B}_2 $ sufficiently close to $ \mathbf{u} $, we have
$ r \| \mathbf{v} - \mathbf{u} \|
\le \| f( \mathbf{v} ) - f( \mathbf{u} ) \| 
\le s \| \mathbf{v} - \mathbf{u} \| $. Any Frechet differentiable function $ f( \cdot ) $ with a Frechet differentiable inverse satisfies such a condition.
If $ \mathbf{T} $ and $ \mathbf{G} $ are random elements with $ \mathbf{T} = f( \mathbf{G} ) $, and the small ball probability of $ \mathbf{G} $ satisfies the bounds described in \eqref{eq:1}, then similar bounds also hold for $ \mathbf{T} $ (see \autoref{prop:3}).
An example of such a non-Gaussian process $ \mathbf{T} $ is the geometric Brownian motion in an $ L_2 $ space, where $ f( \cdot ) $ is the pointwise exponential map \citep[p.~67]{oksendal2003stochastic}.

Next, let $ \mathbf{G} $ be a Gaussian process whose small ball probability satisfies the bounds in \eqref{eq:1}, and $ \mathbf{T} = \mathbf{G} / \mathbf{U} $, where $ \mathbf{U} $ is a bounded positive random variable independent of $ \mathbf{G} $. Then, \eqref{eq:1} will also hold for the small ball probabilities of $ \mathbf{T} $ (see \autoref{prop:1}). Also, bounds similar to \eqref{eq:1} can be established for the small ball probabilities of an infinite dimensional \textit{t}-process, whose corresponding Gaussian process has small ball probabilities satisfying \eqref{eq:1} (see \autoref{prop:2}).

The bounds in \eqref{eq:1} were considered in \citet[p.~209]{ferraty2006nonparametric} with $ C_2 = C_4 $, $ t_1 = t_4 = 0 $, and they called it the small ball probability function of an exponential-type process. For $ t_2 = 0 $, $ t_3 = 1 $ and appropriate values of the parameters $ C_1 $, $ C_2 $, $ C_3 $ and $ C_4 $, \eqref{eq:1} yields the case of a finite dimensional covariate $ \mathbf{X} $ with a continuous positive density at $ \mathbf{x} $, or a fractal-type process as defined in \citet[p.~207]{ferraty2006nonparametric}.

\section{Convergence rate} \label{sec:4}
We now derive the optimum achievable convergence rate for kernel estimates satisfying the bias--variance decomposition \eqref{eq:main}. As we shall see, the function $ m( h ) $ defined in \eqref{eq:mh} plays a central role in determining the convergence rate of the estimate $ \widehat{\Theta}_n( \mathbf{x} ) $. We shall consider the covariate space to be infinite dimensional. The case of finite dimensional covariates is extensively discussed in the past literature (see, e.g., \cite{stone1980optimal,stone1982optimal}, \cite{ibragimov1980nonparametric}, \cite{yatracos1988lower}, \cite{donoho1991geometrizingII,donoho1991geometrizingIII}). 
In order to consider only infinite dimensional covariates, we assume that in \eqref{eq:1}, either $ t_2 > 0 $, or $ t_3 > 1 $ with $ C_2 = C_4 $ in all subsequent discussions.
In that case, $ m( h ) $ is a strictly decreasing positive function, and $ m^{-1}( \cdot ) $, which is the inverse function of $ m( \cdot ) $, is well-defined.
In the next theorem, we see that $ \left( m^{-1}( \log n ) \right)^{\beta} $ is an attainable rate of convergence of $ \widehat{\Theta}_n( \mathbf{x} ) $. Also, under certain additional conditions, $ \left( m^{-1}( \log n ) \right)^{2 \beta} $ is an attainable rate of convergence of the mean square error of $ \widehat{\Theta}_n( \mathbf{x} ) $.
\begin{theorem} \label{thm:up}
Suppose that in \eqref{eq:1}, we have either $ t_2 > 0 $, or $ t_3 > 1 $ with $ C_2 = C_4 $. Then, for any kernel $ K( \cdot ) $ satisfying \ref{assume:a1} and $ \Theta( \mathbf{x} ) $ satisfying \eqref{eq:main} along with conditions \ref{cond:1}--\ref{cond:4}, there is a sequence of bandwidths $ \{ h_n \} $ satisfying \ref{assume:a2} such that $ \big\| \widehat{\Theta}_n( \mathbf{x} ) - \Theta( \mathbf{x} ) \big\| = O_\mathbb{P}\big( \left( m^{-1}( \log n ) \right)^\beta \big) $ as $ n \longrightarrow \infty $, where $ m( h ) $ is as defined in \eqref{eq:mh}.
Further, if $ \mathbb{E}[ \| R_n( \mathbf{x} ) \|^2 ] = o\left( \delta_n^2 \right) $ as $ n \longrightarrow \infty $, where $ \delta_n $ is as defined in \ref{cond:4},
$ \mathbb{E}\big\| \widehat{\Theta}_n( \mathbf{x} ) - \Theta( \mathbf{x} ) \big\|^2 = O\big( \left( m^{-1}( \log n ) \right)^{2 \beta} \big) $ as $ n \longrightarrow \infty $ for the aforementioned sequence of bandwidths $ \{ h_n \} $.
\end{theorem}
Recall that when the parameter of interest is a conditional mean type function as described in \autoref{ex:2.1} in \autoref{sec:2}, $ R_n( \mathbf{x} ) = \mathbf{0} $. So, in that case the condition $ \mathbb{E}[ \| R_n( \mathbf{x} ) \|^2 ] = o( \delta_n^2 ) $ assumed in the second part of the above theorem is trivially satisfied.

\subsection{Lower bound on the convergence rate} \label{subsec:4_1}
We now proceed to investigate the lower bound of the convergence rate of $ \widehat{\Theta}_n( \mathbf{x} ) $.
In the next proposition, we establish an asymptotic lower bound of the sequence of bandwidths $ \{ h_n \} $ that leads to consistent kernel regression estimates. This result will be needed while deriving the lower bound of the convergence rate of a kernel estimate.
\begin{theorem} \label{thm:1a}
Suppose that in the upper and the lower bounds in the shifted small ball probability in \eqref{eq:1}, we have either $ t_2 > 0 $, or $ t_3 > 1 $ with $ C_2 = C_4 $. Then, for any sequence of bandwidths $ \{ h_n \} $, which satisfies assumption \ref{assume:a2}, we have $ h_n / m^{-1}( \log n ) $ bounded away from $ 0 $ as $ n \longrightarrow \infty $, where $ m( h ) $ is as defined in \eqref{eq:mh}.
\end{theorem}

Define
\begin{align*}
\tilde{B}_n( \mathbf{x} ) = \mathbb{L}_\mathbf{x}\left( \mathbb{E}\left[ \left( F( \mathbf{X} ) - F( \mathbf{x} ) \right) \allowbreak \frac{K( h_n^{-1} d( \mathbf{x}, \mathbf{X} ) )}{E_n^{(1)}( \mathbf{x} ) \phi( \mathbf{x}, h_n )} \right] \right) ,
\end{align*}
where $ \mathbb{L}_\mathbf{x}( \cdot ) $ and $ F( \cdot ) $ are as defined after \eqref{eq:variance}, and $ E_n^{(1)}( \mathbf{x} ) $ is as defined in \eqref{eq:assume:a1}.
Also, let $ \{ \mathbf{e}_1, \mathbf{e}_2, \ldots \} $ be a Schauder basis of $ \mathcal{B} $, such that for any $ \mathbf{v} \in \mathcal{B} $, $ \mathbf{v} = \sum_{n=1}^{\infty} v_n \mathbf{e}_n $ for a sequence of real numbers $ \{ v_n \} $. Let $ \tilde{\phi}_i \in \mathcal{B}^* $ be the projection functional corresponding to $ \mathbf{e}_i $, i.e., $ \mathbf{v} = \sum_{i=1}^{\infty} \tilde{\phi}_i( \mathbf{v} ) \mathbf{e}_i $ for all $ \mathbf{v} \in \mathcal{B} $.
Consider the following assumptions.
\begin{enumerate}[label=C(\roman*), ref=C(\roman*)]
\item \label{cond:6}
There is $ \Theta( \cdot ) : \mathcal{C} \longrightarrow \mathcal{B} $ with the corresponding $ \mathbb{L}_\mathbf{x}( \cdot ) $ and $ F( \cdot ) $ such that for any sequence of bandwidths $ \{ h_n \} $ satisfying \ref{assume:a2},
\begin{align} \label{eq:b3}
h_n^{-\beta} \left\| \tilde{B}_n( \mathbf{x} ) \right\| > b_1 > 0
\end{align}
for all sufficiently large $ n $.

\item \label{cond:7}
Let $ G( \cdot ) $ be as defined after \eqref{eq:variance}. For some positive integer $ i_0 $, the conditional variance function $ \mathbb{V}( \mathbf{z} ) : \mathcal{C} \longrightarrow \mathbb{R} $ defined by $ \mathbb{V}( \mathbf{z} ) = \mathbb{E}\big[ \big( \tilde{\phi}_{i_0}\big( \mathbb{L}_\mathbf{x}( G( \mathbf{Y} ) - \mathbb{E}[ G( \mathbf{Y} ) \midil \mathbf{X} = \mathbf{z} ] ) \big) \big)^2 \,\big|\, \mathbf{X} = \mathbf{z} \big] $ converges to $ \mathbb{V}( \mathbf{x} ) $ as $ \mathbf{z} \longrightarrow \mathbf{x} $, and $ \mathbb{V}( \mathbf{x} ) > 0 $.
\end{enumerate}
Condition \ref{cond:7}, like condition \ref{cond:5}, is related to the smoothness of the conditional distribution of the residual in the regression. In fact, condition \ref{cond:7} holds in the same location-scale type models, which we described after condition \ref{cond:5}. Condition \ref{cond:6} gives a lower bound on the rate of convergence of the bias part of the estimate. Inequality \eqref{eq:b2} and condition \ref{cond:6} together imply that the rate of convergence of the bias part is same as $ h_n^\beta $ as $ n \longrightarrow \infty $. The following two conditions are sufficient to ensure that \ref{cond:6} holds.
\begin{enumerate}[label=(\alph*), ref=(\alph*)]
\item \label{cond:6a}
There is a constant $ 0 < s < 1 $ such that $ \phi( \mathbf{x}, s h ) / \phi( \mathbf{x}, h ) $ is bounded away from 1 for all sufficiently small $ h > 0 $.

\item \label{cond:6b}
Let $ \mathbb{L}_\mathbf{x}( \mathcal{F}(\mathbf{x}, \beta, \mathcal{G} ) ) $ be the class of all functions defined by the composition $ \mathbb{L}_\mathbf{x} \circ H $, where $ H \in \mathcal{F}(\mathbf{x}, \beta, \mathcal{G} ) $ and $ \mathcal{F}(\mathbf{x}, \beta, \mathcal{G} ) $ is as defined in \ref{cond:1}. Then $ \mathbb{L}_\mathbf{x}( \mathcal{F}(\mathbf{x}, \beta, \mathcal{G} ) ) $ contains the function $ \mathbf{z} \mapsto d( \mathbf{x}, \mathbf{z} )^\beta \mathbf{v} $, where $ \mathbf{v} \in \mathcal{B} $, and $ \mathbf{z} $ lies in a neighborhood of $ \mathbf{x} $.
\end{enumerate}
Condition \ref{cond:6a} is satisfied when in \eqref{eq:1}, $ t_2 > 0 $, or $ t_1 < t_4 $, or $ C_2 = C_4 $. We observe that at least one of these is true in the examples that we have described in \autoref{subsec:3_2}.

Now, we derive the lower bound of the order of convergence of the bias term $ B_n( \mathbf{x} ) $ in \eqref{eq:main} under \ref{cond:1} and \ref{cond:6}. Note that $ B_n( \mathbf{x} ) = \tilde{B}_n( \mathbf{x} ) + \tilde{R}_n( \mathbf{x} ) $, where
\begin{align*}
& \tilde{R}_n( \mathbf{x} ) \\
&= \mathbb{L}_\mathbf{x}\left( \frac{ \sum_{i=1}^{n} ( F( \mathbf{X}_i )- F( \mathbf{x} ) ) K( h_n^{-1} d( \mathbf{x}, \mathbf{X}_i ) ) }{ \sum_{i=1}^{n} K( h_n^{-1} d( \mathbf{x}, \mathbf{X}_i ) ) } - \tilde{B}_n( \mathbf{x} ) \right) \\
&= \mathbb{L}_\mathbf{x}\left( \frac{ \sum_{i=1}^{n} ( F( \mathbf{X}_i ) - F( \mathbf{x} ) ) K\left( \frac{d( \mathbf{x}, \mathbf{X}_i )}{h_n} \right) }{ \sum_{i=1}^{n} K\left( \frac{d( \mathbf{x}, \mathbf{X}_i )}{h_n} \right) } 
- \frac{1}{n} \sum_{i=1}^{n} \frac{( F( \mathbf{X}_i )- F( \mathbf{x} ) ) K\left( \frac{d( \mathbf{x}, \mathbf{X}_i )}{h_n} \right)}{E_n^{(1)}( \mathbf{x} ) \phi( \mathbf{x}, h_n )} \right) \\
& \quad + \mathbb{L}_\mathbf{x}\left( \frac{1}{n} \sum_{i=1}^{n} \frac{( F( \mathbf{X}_i ) - F( \mathbf{x} ) ) K\left( \frac{d( \mathbf{x}, \mathbf{X}_i )}{h_n} \right)}{E_n^{(1)}( \mathbf{x} ) \phi( \mathbf{x}, h_n )}
- \mathbb{E}\left[ \frac{( F( \mathbf{X} ) - F( \mathbf{x} ) ) K\left( \frac{d( \mathbf{x}, \mathbf{X} )}{h_n} \right)}{E_n^{(1)}( \mathbf{x} ) \phi( \mathbf{x}, h_n )} \right] \right) .
\end{align*}
It follows form condition \ref{assume:a1} and Markov inequality that
\begin{align*}
\left( \frac{1}{n} \sum_{i=1}^{n} \allowbreak \frac{K( h_n^{-1} d( \mathbf{x}, \mathbf{X}_i ) )}{E_n^{(1)}( \mathbf{x} ) \phi( \mathbf{x}, h_n )} - 1 \right) = O_\mathbb{P}\left( \big[ n \phi( \mathbf{x}, h_n ) \big]^{-1/2} \right)
\end{align*}
as $ n \longrightarrow \infty $. Hence, from conditions \ref{assume:a1}, \ref{assume:a2} and \ref{cond:1}, we have
\begin{align*}
& \mathbb{L}_\mathbf{x}\left( \frac{ \sum_{i=1}^{n} ( F( \mathbf{X}_i ) - F( \mathbf{x} ) ) K\left( \frac{d( \mathbf{x}, \mathbf{X}_i )}{h_n} \right) }{ \sum_{i=1}^{n} K\left( \frac{d( \mathbf{x}, \mathbf{X}_i )}{h_n} \right) } 
- \frac{1}{n} \sum_{i=1}^{n} \frac{( F( \mathbf{X}_i )- F( \mathbf{x} ) ) K\left( \frac{d( \mathbf{x}, \mathbf{X}_i )}{h_n} \right)}{E_n^{(1)}( \mathbf{x} ) \phi( \mathbf{x}, h_n )} \right) \\
& \quad = o_\mathbb{P}( h_n^\beta )
\end{align*}
as $ n \longrightarrow \infty $. Also, from assumptions \ref{assume:a1}, \ref{assume:a2}, \ref{cond:1} and Markov inequality, it follows that
\begin{align*}
& \mathbb{L}_\mathbf{x}\left( \frac{1}{n} \sum_{i=1}^{n} \frac{( F( \mathbf{X}_i ) - F( \mathbf{x} ) ) K\left( \frac{d( \mathbf{x}, \mathbf{X}_i )}{h_n} \right)}{E_n^{(1)}( \mathbf{x} ) \phi( \mathbf{x}, h_n )}
- \mathbb{E}\left[ \frac{( F( \mathbf{X} ) - F( \mathbf{x} ) ) K\left( \frac{d( \mathbf{x}, \mathbf{X} )}{h_n} \right)}{E_n^{(1)}( \mathbf{x} ) \phi( \mathbf{x}, h_n )} \right] \right) \\
& \quad = o_\mathbb{P}( h_n^\beta )
\end{align*}
as $ n \longrightarrow \infty $. Hence,
\begin{align} \label{eq:b1}
\tilde{R}_n( \mathbf{x} ) = o_\mathbb{P}( h_n^\beta ) \text{ as } n \longrightarrow \infty .
\end{align}
Note that inequality \eqref{eq:b3} provides a lower bound of the convergence rate for the bias term $ B_n( \mathbf{x} ) $ in view of \eqref{eq:b1}, and this will be used to determine a lower bound of the rate of convergence of $ \widehat{\Theta}_n( \mathbf{x} ) $.
Also note that $ \tilde{\phi}_{i_0}\big( \mathbb{L}_\mathbf{x}( G( \mathbf{Y} ) ) \big) $ in condition \ref{cond:7} is a real-valued random variable. So, the convergence condition in \ref{cond:7} of the conditional variance may be viewed as a special case of condition \ref{cond:5}.
We now state the theorem on the lower bound of the convergence rate of $ \big\| \widehat{\Theta}_n( \mathbf{x} ) - \Theta( \mathbf{x} ) \big\| $.
\begin{theorem} \label{thm:low}
Suppose that in \eqref{eq:1}, we have either $ t_2 > 0 $, or $ t_3 > 1 $ with $ C_2 = C_4 $, the kernel $ K( \cdot ) $ satisfies \ref{assume:a1}, the sequence of bandwidths $ \{ h_n \} $ satisfies \ref{assume:a2}, and the decomposition \eqref{eq:main} along with conditions \ref{cond:1}--\ref{cond:4}, \ref{cond:6} and \ref{cond:7} hold. Then,
\begin{align*}
\liminf_{n \longrightarrow \infty} \mathbb{P}\left[ \left( m^{-1}( \log n ) \right)^{-\beta} \allowbreak \left\| \widehat{\Theta}_n( \mathbf{x} ) - \Theta( \mathbf{x} ) \right\| > c \right] > 0
\end{align*}
for some constant $ c > 0 $ depending on $ \Theta( \mathbf{x} ) $, where $ m( h ) $ is as defined in \eqref{eq:mh}.
\end{theorem}

\autoref{thm:low} implies that we cannot get a faster rate of convergence than $ \left( m^{-1}( \log n ) \right)^{\beta} $, since $ \left( m^{-1}( \log n ) \right)^{-\beta} \allowbreak \big\| \widehat{\Theta}_n( \mathbf{x} ) - \Theta( \mathbf{x} ) \big\| $ does not converge to $ 0 $ \textit{in probability} as $ n \longrightarrow \infty $.
Further, from \autoref{thm:low} it follows that $ \left( m^{-1}( \log n ) \right)^{2 \beta} $ is a lower bound for the rate of convergence of the mean square error $ \mathbb{E}\big\| \widehat{\Theta}_n( \mathbf{x} ) - \Theta( \mathbf{x} ) \big\|^2 $.
Hence, combining \autoref{thm:up} and \autoref{thm:low}, we get that $ \left( m^{-1}( \log n ) \right)^{\beta} $ and $ \left( m^{-1}( \log n ) \right)^{2 \beta} $ are the optimum rates of convergence of $ \widehat{\Theta}_n( \mathbf{x} ) $ and its mean square error, respectively, when all the conditions of the two theorems are satisfied. We now deduce simplified expressions of the optimum rates for the specific infinite dimensional covariate distributions considered in \autoref{subsec:3_2}.

For $ \mathbf{X} $ being a real-valued continuous Gaussian Markov process on $ [ 0, 1 ] $, under the $ L_p $-norm, we have $ \left( m^{-1}( \log n ) \right)^{\beta} = O( ( \log n )^{- \beta / 2} ) $ as $ n \longrightarrow \infty $. For fractional Brownian motion on $ [ 0, 1 ] $ with Hurst index $ \gamma \in ( 0, 1 ) $, under the $ L_\infty $-norm, we have $ t_2 = 1 / \gamma $, and consequently $ \left( m^{-1}( \log n ) \right)^{\beta} = O( ( \log n )^{- \gamma \beta} ) $ as $ n \longrightarrow \infty $. On the other hand, for an integrated fractional Brownian motion with Hurst index $ \gamma $ and under the $ L_\infty $-norm, we have $ t_2 = 1 / ( 1 + \gamma ) $ and $ \left( m^{-1}( \log n ) \right)^{\beta} = O( ( \log n )^{- ( 1 + \gamma ) \beta} ) $ as $ n \longrightarrow \infty $. When $ \mathbf{X} $ is a L\'{e}vy fractional Brownian motion on $ [ 0, 1 ]^q $ with Hurst index $ \gamma $, $ t_2 = q / \gamma $ and $ \left( m^{-1}( \log n ) \right)^{\beta} = O( ( \log n )^{- \gamma \beta / q} ) $ as $ n \longrightarrow \infty $.

In the class of processes $ H_{X,L} $ considered in subsection 4.1 of \cite{chagny2014adaptive}, $ t_2 > 0 $ and $ t_3 = 0 $, and we have $ \left( m^{-1}( \log n ) \right)^{2 \beta} = O( ( \log n )^{- 2 \beta / t_2} ) $ as $ n \longrightarrow \infty $. On the other hand, for the class of processes $ H_{X,M} $ considered by these authors, we have $ \left( m^{-1}( \log n ) \right)^{2 \beta} = O\big( \exp[ - 2 \beta C_2^{- 1 / t_3} \allowbreak ( \log n )^{1 / t_3} ] \big) $ as $ n \longrightarrow \infty $.
Note that these rates coincide with the optimal rates of convergence of the mean square error described in \citet[Table 1, p.~2363]{chagny2014adaptive}, which were derived when the response is real-valued and the parameter of interest is the conditional distribution function of the response. We have covered this particular case of the parameter of interest in \autoref{ex:2.1}.

\subsection{Asymptotic dominance of bias over variance} \label{subsec:4_2}
Recall that in the case of finite dimensional covariates, the bias and the variance terms in nonparametric regression have the same rate of convergence (see, e.g., \cite{stone1980optimal}, \cite{ibragimov1980nonparametric}).
In fact, \cite{mas2012lower} chose the bandwidth of the kernel estimate by balancing the asymptotic orders of the bias and the variance (see Lemma 1 and the preceding discussion in \cite{mas2012lower}) even when the covariate is infinite dimensional. However, as we shall show now, the optimum choice of the bandwidth in a kernel estimate, as described in the proof of \autoref{thm:up}, leads to different asymptotic orders of the bias and the variance when the covariate is infinite dimensional in nature, i.e., when we have either $ t_2 > 0 $ or $ t_3 > 1 $ with $ C_2 = C_4 $ in \eqref{eq:1}.

\begin{theorem} \label{thm:3}
Suppose that either $ t_2 > 0 $ or $ t_3 > 1 $ with $ C_2 = C_4 $ in the bounds in \eqref{eq:1}. Also, let the kernel $ K( \cdot ) $ satisfy \ref{assume:a1}, and the decomposition \eqref{eq:main} along with conditions \ref{cond:1}--\ref{cond:4} hold.
Then, for any $ \Theta( \mathbf{x} ) $ satisfying \ref{cond:6}, the ratio $ \| V_n( \mathbf{x} ) \| / \| B_n( \mathbf{x} ) \| = o_\mathbb{P}( 1 ) $ as $ n \longrightarrow \infty $ for the optimum choice of bandwidth $ \{ h_n \} $ described in the proof of \autoref{thm:up}.
\end{theorem}

\autoref{thm:3} illustrates that our optimum bandwidth, which minimizes \eqref{eq:thm:up} in the proof of \autoref{thm:up}, does not balance the convergence rates of the variance and the bias in kernel regression if the covariate is infinite dimensional. Instead, the ratio of the variance to the bias for our optimal choice of bandwidth tends to zero as the sample size increases. This phenomenon is due to the exponential decay of the small ball probability function in infinite dimensional spaces. When the covariate is infinite dimensional, we may have very small number of observations in a neighborhood in the covariate space. To cope with this problem, one has to use relatively larger bandwidths than what is required for finite dimensional covariates. This results in an `over-smoothed' estimate with its bias asymptotically larger than its variance.
It will be appropriate to note here that the optimum convergence rate derived in \autoref{thm:up} and \autoref{thm:low} is same as the one derived in \cite{mas2012lower} for estimation of the conditional mean of a real-valued response, where the chosen bandwidth balances the bias and the variance \citep[p.~1760]{mas2012lower}. However, our optimum bandwidth, which does not try to balance the bias and the variance in the decomposition \eqref{eq:main}, will often lead to an estimate with higher statistical precision compared to an estimate based on a bandwidth that balances the bias and the variance. In several cases, the statistical error will be substantially lower when our optimum bandwidth is used as demonstrated in \autoref{coro:1}.

\begin{theorem} \label{coro:1}
Suppose  assumptions \ref{assume:a1}, \ref{assume:a2}, \ref{cond:1}--\ref{cond:4} and \ref{cond:7} are satisfied.
Let $ \widehat{\Theta}_n^{(b)}( \mathbf{x} ) $ be an estimate of $ \Theta( \mathbf{x} ) $ constructed using bandwidth $ h_n^{(b)} $, which satisfies \ref{assume:a2} and balances the bias and the variance so that $ ( h_n^{(b)} )^{2 \beta} n \phi( \mathbf{x}, h_n^{(b)} ) $ is bounded and bounded away from 0 as $ n \longrightarrow \infty $. Also, let $ \widehat{\Theta}_n^{(op)}( \mathbf{x} ) $ be an estimate of $ \Theta( \mathbf{x} ) $ constructed using our optimum bandwidth minimizing \eqref{eq:thm:up} in the proof of \autoref{thm:up}.
Assume that $ t_2 > 0 $ in the bounds in \eqref{eq:1}.
Then, for any $ \beta_1 > \beta $ and any $ \Theta( \cdot ) $ for which the corresponding $ F( \cdot ) \in \mathcal{F}( \mathbf{x}, \beta_1, \mathcal{G} ) \subseteq \mathcal{F}( \mathbf{x}, \beta, \mathcal{G} ) $,
\begin{align*}
\frac{\left\| \widehat{\Theta}_n^{(op)}( \mathbf{x} ) - \Theta( \mathbf{x} ) \right\|}{\left\| \widehat{\Theta}_n^{(b)}( \mathbf{x} ) - \Theta( \mathbf{x} ) \right\|}
= o_\mathbb{P}( 1 )
\quad
\text{as }
n \longrightarrow \infty .
\end{align*}
Further, if $ \mathbb{E}[ \| R_n( \mathbf{x} ) \|^2 ] = o\left( \delta_n^2 \right) $ as $ n \longrightarrow \infty $, where $ \delta_n $ is as defined in \ref{cond:4}, then
\begin{align*}
\frac{\mathbb{E}\left\| \widehat{\Theta}_n^{(op)}( \mathbf{x} ) - \Theta( \mathbf{x} ) \right\|^2}{\mathbb{E}\left\| \widehat{\Theta}_n^{(b)}( \mathbf{x} ) - \Theta( \mathbf{x} ) \right\|^2}
= o( 1 )
\quad
\text{as }
n \longrightarrow \infty .
\end{align*}
\end{theorem}

Recall that for conditional mean type functions described in \autoref{ex:2.1} in \autoref{sec:2}, $ R_n( \mathbf{x} ) = \mathbf{0} $, and the condition $ \mathbb{E}[ \| R_n( \mathbf{x} ) \|^2 ] = o( \delta_n^2 ) $ assumed in the second part of the above theorem is trivially satisfied.

\section{Adaptive selection of bandwidths} \label{sec:adaptive}
In practice, one has to choose the bandwidth $ h $ by some data-driven adaptive procedure. Such adaptive choice of bandwidth, when the covariate is functional, has been investigated in \cite{chagny2014adaptive,chagny2016adaptive} for the kernel estimates of the conditional distribution and the conditional mean of a real-valued response. Their data-based bandwidth selection procedure can be suitably extended for more general regression problems considered in this paper.

Let $ \mathbb{H}_n $ be a finite collection of bandwidths with cardinality less than or equal to $ n $ such that for any $ h \in \mathbb{H}_n $, $ \phi( \mathbf{x}, h ) \le 2 ( \log n )^{-1} $ and $ \phi( \mathbf{x}, h ) \ge n^{-1} ( \log n )^2 $. Since $ \phi( \mathbf{x}, h ) $ is a monotone increasing function of $ h $, if a sequence of bandwidths $ \{ h_n \} $ is such that $ h_n \in \mathbb{H}_n $ for all $ n $, then $ \{ h_n \} $ satisfies condition \ref{assume:a2}.
In this section, we shall write $ \widehat{\Theta}_n( \mathbf{x}, h ) $, $ B_n( \mathbf{x}, h ) $, $ V_n( \mathbf{x}, h ) $ and $ R_n( \mathbf{x}, h ) $ for $ \widehat{\Theta}_n( \mathbf{x} ) $, $ B_n( \mathbf{x} ) $, $ V_n( \mathbf{x} ) $ and $ R_n( \mathbf{x} ) $, respectively, to indicate the dependence of $ \widehat{\Theta}_n( \mathbf{x} ) $, $ B_n( \mathbf{x} ) $, $ V_n( \mathbf{x} ) $ and $ R_n( \mathbf{x} ) $ on the bandwidth $ h $.
We assume the following:
\begin{enumerate}[label=D(\roman*), ref=D(\roman*)]
\item \label{cond:adaptive2}
There is a constant $ \sigma > 0 $ such that
for any $ \mathbf{z} $ in a certain neighborhood of $ \mathbf{x} $ and every integer $ k \ge 2 $,
\begin{align*}
\mathbb{E}\left[ \| \mathbb{L}_\mathbf{x}\left( G( \mathbf{Y} ) - \mathbb{E}[ G( \mathbf{Y} ) \midil \mathbf{X} = \mathbf{z} ] \right) \|^k \middle\arrowvert \mathbf{X} = \mathbf{z} \right] \le \frac{k!}{2} \sigma^k .
\end{align*}

\item \label{cond:adaptive3}
There are constants $ \epsilon_1 > 0 $, $ \epsilon_2 > 0 $ and $ M > 0 $ such that whenever $ h \le \epsilon_1 $ and $ \| V_n( \mathbf{x}, h ) \| \le \epsilon_2 $, we have
$ \| R_n( \mathbf{x}, h ) \|^2 \le M h^{2 \beta} + \| V_n( \mathbf{x}, h ) \|^2 $.
\end{enumerate}
Condition \ref{cond:adaptive2} is similar to assumption ($ H_\epsilon $) used in \citet[p.~108]{chagny2016adaptive}, which was used to derive the convergence rate of the adaptive estimate of the conditional mean for a real-valued response.
Condition D(ii) describes a bound on the remainder of our bias--variance type decomposition in terms of the bound on the bias part and the variance part.
Condition \ref{cond:adaptive3} is trivially satisfied for the conditional mean type estimates described in \autoref{ex:2.1} in \autoref{sec:2}.
It is also satisfied in the class of regression problems described in \autoref{ex:2.2} in \autoref{sec:2}.
Define the empirical shifted small ball probability $ \widehat{\phi}( \mathbf{x}, h ) = ( 1/n ) \sum_{i=1}^{n} \mathbb{I}( d( \mathbf{x}, \mathbf{X}_i ) \le h ) $. Define
\begin{align*}
D_n( \mathbf{x}, h ) = 
\sigma^2 \zeta_n \frac{\log n}{n \hat{\phi}( \mathbf{x}, h )} \mathbb{I}\left( \hat{\phi}( \mathbf{x}, h ) > 0 \right)
+ \sigma^2 \zeta_n n \mathbb{I}\left( \hat{\phi}( \mathbf{x}, h ) = 0 \right) ,
\end{align*}
where $ \{ \zeta_n \} $ is a sequence of positive constants independent of $ h $, such that $ \zeta_n \longrightarrow \zeta_0 > 0 $ as $ n \longrightarrow \infty $.
The constant $ \zeta_0 $ is described in the proof of \autoref{lemma:adapt4}.
Also define
\begin{align*}
C_n( \mathbf{x}, h ) = \max_{ h' \in \mathbb{H}_n } \left( \left\| \widehat{\Theta}_n( \mathbf{x}, h' ) - \widehat{\Theta}_n( \mathbf{x}, \max\{ h, h' \} ) \right\|^2 - D_n( \mathbf{x}, h' ) \right)_+ .
\end{align*}
$ D_n( \mathbf{x}, h ) $ approximates the upper bound of the variance term and $ C_n( \mathbf{x}, h ) $ approximates the bias term. The data-driven choice of bandwidth is defined as
\begin{align*}
h_n^* = \arg \min_{h \in \mathbb{H}_n} \allowbreak \left[ C_n( \mathbf{x}, h ) + D_n( \mathbf{x}, h ) \right] .
\end{align*}
The following theorem gives an upper bound on the convergence rate of the adaptive estimate $ \widehat{\Theta}_n( \mathbf{x}, h_n^* ) $.
\begin{theorem} \label{thm:5}
Define
\begin{align*}
\lambda_n
= \min\limits_{h \in \mathbb{H}_n}\left[ h^{2 \beta} + \frac{\log n}{n \phi( \mathbf{x}, h )} \right] .
\end{align*}
Let conditions \ref{assume:a1}, \ref{cond:1}, \ref{cond:adaptive2} and \ref{cond:adaptive3} be satisfied. Then,
\begin{align}
\left\| \widehat{\Theta}_n( \mathbf{x}, h_n^* ) - \Theta( \mathbf{x} ) \right\|^2
= O_\mathbb{P}\left( \lambda_n \right)
\quad
\text{as } n \longrightarrow \infty .
\label{eq:8}
\end{align}
Further, for the conditional mean-type functions described in \autoref{ex:2.1} in \autoref{sec:2}, we have
\begin{align}
\mathbb{E}\left\| \widehat{\Theta}_n( \mathbf{x}, h_n^* ) - \Theta( \mathbf{x} ) \right\|^2 
= O\left( \lambda_n \right)
\quad
\text{as } n \longrightarrow \infty .
\label{eq:9}
\end{align}
\end{theorem}

Equation \eqref{eq:8} gives an upper bound for the asymptotic convergence rate of the adaptive estimate. In \cite{chagny2016adaptive}, the adaptive estimate and its convergence rate were derived for the estimation of the conditional mean of a real-valued response in a homoscedastic model. Our setup includes heteroscedastic regression models where the parameter to be estimated is an element of a type 2 Banach space, and it is not necessarily the conditional mean.

\subsection{Numerical demonstration}
\label{subsec:5.1}
We next demonstrate the adaptive estimate $ \widehat{\Theta}_n( \mathbf{x}, h_n^* ) $ in several regression models. In all the examples, we consider the covariate $ \mathbf{X} $ to be a random element in $ L_2[ 0, 1 ] $. The usual norm in $ L_2[ 0, 1 ] $ is denoted as $ \| \cdot \|_2 $. We denote the adaptive choice of the bandwidth as $ h_n^* $. We take $ \zeta_n = \min\{ \sqrt{n}, 1500 \} $ in our computation, the validity of which is ensured from \eqref{eq:adaptivezetaboundfinal}. We substitute $ \phi( \mathbf{x}, h ) $ by $ \hat{\phi}( \mathbf{x}, h ) $ in the construction of the collection of bandwidths $ \mathbb{H}_n $, as done in \cite{chagny2016adaptive}. The parameter $ \sigma^2 $ used to define $ D_n( \mathbf{x}, h ) $ and mentioned in \ref{cond:adaptive2} also needs to be estimated. This is done based on the regression model.
Note that for $ \sigma^2 $, which satisfies 
\begin{align}
\sigma^2 \ge
\| \mathbb{L}_\mathbf{x} \|^2 \sup_{\mathbf{z}} \mathbb{E}\left[ \| G( \mathbf{Y} ) - \mathbb{E}[ G( \mathbf{Y} ) \midil \mathbf{X} = \mathbf{z} ] \|^2 \middle\arrowvert \mathbf{X} = \mathbf{z} \right]
\label{eq:adapt.1}
\end{align}
for $ \mathbf{z} $ lying in some neighborhood of $ \mathbf{x} $, condition \ref{cond:adaptive2} will hold.
Since by construction $ \max \mathbb{H}_n \longrightarrow 0 $ as $ n \longrightarrow \infty $, and the kernel $ K( u ) = I( 0 \le u \le 1 ) $ satisfies condition \ref{assume:a1}, it is enough to consider the supremum in \eqref{eq:adapt.1} over the $ \mathbf{X}_i $s such that $ d( \mathbf{x}, \mathbf{X}_i ) \le \max \mathbb{H}_n $ for estimating $ \sigma^2 $.
So, if $ \hat{\sigma}_1^2 $ is an estimated upper bound of $ \| \mathbb{L}_\mathbf{x} \|^2 $, and $ \hat{\sigma}_2^2( \mathbf{X}_i ) $ is an estimated upper bound of $ \mathbb{E}\left[ \| G( \mathbf{Y}_i ) - \mathbb{E}[ G( \mathbf{Y}_i ) \midil \mathbf{X}_i ] \|^2 \middle\arrowvert \mathbf{X}_i \right] $, then we can take
\begin{align*}
\hat{\sigma}^2 = \hat{\sigma}_1^2 \max\{ \hat{\sigma}_2^2( \mathbf{X}_i ) \midil d( \mathbf{x}, \mathbf{X}_i ) \le \max \mathbb{H}_n \}
\end{align*}
as an estimate of $ \sigma^2 $.
Let $ h_{n,1} = \min \mathbb{H}_n $ and $ h_{n,2} = \max \mathbb{H}_n $.
Denote
\begin{align*}
W_{i,n}^{(1)}( \mathbf{z} ) = \frac{K( {h_{n,1}}^{-1} d( \mathbf{z}, \mathbf{X}_i ) )}{\sum_{i=1}^{n} K( {h_{n,1}}^{-1} d( \mathbf{z}, \mathbf{X}_i ) )} .
\end{align*}
In the case of the mean regression model as described in \autoref{ex:2.1}, an estimate of $ \sigma^2 $ is
\begin{align*}
\hat{\sigma}^2 = \max\left\{ \sum_{j=1}^{n} \left\| \Psi( \mathbf{Y}_j ) - \left( \sum_{k=1}^{n} \Psi( \mathbf{Y}_k ) W_{k,n}^{(1)}( \mathbf{X}_j ) \right) \right\|^2 W_{j,n}^{(1)}( \mathbf{X}_i ) \middle\arrowvert d( \mathbf{x}, \mathbf{X}_i ) \le h_{n,2} \right\} .
\end{align*}
The function $ \Psi( \cdot ) $ is as described in \autoref{ex:2.1}.
The rationale for using the weights $ W_{i,n}^{(1)}( \mathbf{z} ) $ is the same as that described in subsection 4.1.2 in \cite{chagny2016adaptive}.
In case the parameter to be estimated is a function of the conditional mean as discussed in \autoref{ex:2.2}, or in the case of a maximum likelihood regression model as described in \autoref{ex:2.3}, we need to additionally estimate an upper bound of the term $ \| \mathbb{L}_\mathbf{x} \|^2 $ in \eqref{eq:adapt.1}.
For a function of conditional mean type estimate, we have seen in \autoref{ex:prop1} that $ \mathbb{L}_\mathbf{x}( \cdot ) = \Gamma'\left( \mathbb{E}[ \Psi( \mathbf{Y} ) \midil \mathbf{X} = \mathbf{x} ] \right)( \cdot ) $, where $ \Gamma( \cdot ) $ and $ \Psi( \cdot ) $ are as described in \autoref{ex:2.2}.
Hence, we can take $ \hat{\sigma}_1^2 = \left\| \Gamma'\left( \sum_{i=1}^{n} \Psi( \mathbf{Y}_i ) W_{i,n}^{(1)}( \mathbf{x} ) \right) \right\|^2 $, and
\begin{align*}
\hat{\sigma}^2 = \hat{\sigma}_1^2 \max\left\{ \sum_{j=1}^{n} \left\| \Psi( \mathbf{Y}_j ) - \left( \sum_{k=1}^{n} \Psi( \mathbf{Y}_k ) W_{k,n}^{(1)}( \mathbf{X}_j ) \right) \right\|^2 W_{j,n}^{(1)}( \mathbf{X}_i ) \middle\arrowvert d( \mathbf{x}, \mathbf{X}_i ) \le h_{n,2} \right\} .
\end{align*}
Here, the function $ \Psi( \cdot ) $ is as described in \autoref{ex:2.2}.
In a maximum likelihood regression model (\autoref{ex:2.3}), we have seen from \autoref{ex:prop2} that $ \mathbb{L}_\mathbf{x}( \cdot ) = [ \mathbf{I}( \Theta( \mathbf{x} ) ) ]^{-1}( \cdot ) $. So, we need to estimate $ \mathbf{I}( \Theta( \mathbf{x} ) ) = - \mathbb{E}[ \Delta_2( g( \mathbf{Y} \midil \Theta( \mathbf{x} ) ) ) \midil \mathbf{X} = \mathbf{x} ] $, which we estimate by
$ \hat{\mathbf{I}}
= - \sum_{i=1}^{n} \Delta_2( g( \mathbf{Y}_i \mid \widehat{\Theta}_n^{(1)}( \mathbf{x} ) ) ) W_{i,n}^{(1)}( \mathbf{x} ) $,
where $ \widehat{\Theta}_n^{(1)}( \mathbf{x} ) $ is defined as the solution of the likelihood equation \eqref{ex2.3eq1a} with the bandwidth being $ h_{n,1} $.
So, we can take
$ \hat{\sigma}_1^2 = \big\| \hat{\mathbf{I}} \big\|^{-2} $. In this case, $ G( \mathbf{Y} ) = \nabla g( \mathbf{Y} \midil \Theta( \mathbf{X} ) ) $, so that $ \mathbb{E}[ G( \mathbf{Y}_i ) \midil \mathbf{X}_i ] = \mathbf{0} $ for all $ i $. Since $ \Theta( \mathbf{X}_i ) $ is unknown, we estimate $ G( \mathbf{Y}_i ) $ by $ \nabla g( \mathbf{Y}_i \midil \widehat{\Theta}_n^{(1)}( \mathbf{X}_i ) ) $, where $ \widehat{\Theta}_n^{(1)}( \mathbf{X}_i ) $ is the solution of the likelihood equation \eqref{ex2.3eq1a} with $ \mathbf{x} $ replaced by $ \mathbf{X}_i $ and the bandwidth being $ h_{n,1} $. Consequently, here we have
\begin{align*}
\hat{\sigma}^2 = \hat{\sigma}_1^2 \max\left\{ \sum_{j=1}^{n} \left\| \nabla g( \mathbf{Y}_j \midil \widehat{\Theta}_n^{(1)}( \mathbf{X}_j ) ) \right\|^2 W_{j,n}^{(1)}( \mathbf{X}_i ) \middle\arrowvert d( \mathbf{x}, \mathbf{X}_i ) \le h_{n,2} \right\} .
\end{align*}

\begin{figure}
\centering
\includegraphics[width=1\linewidth]{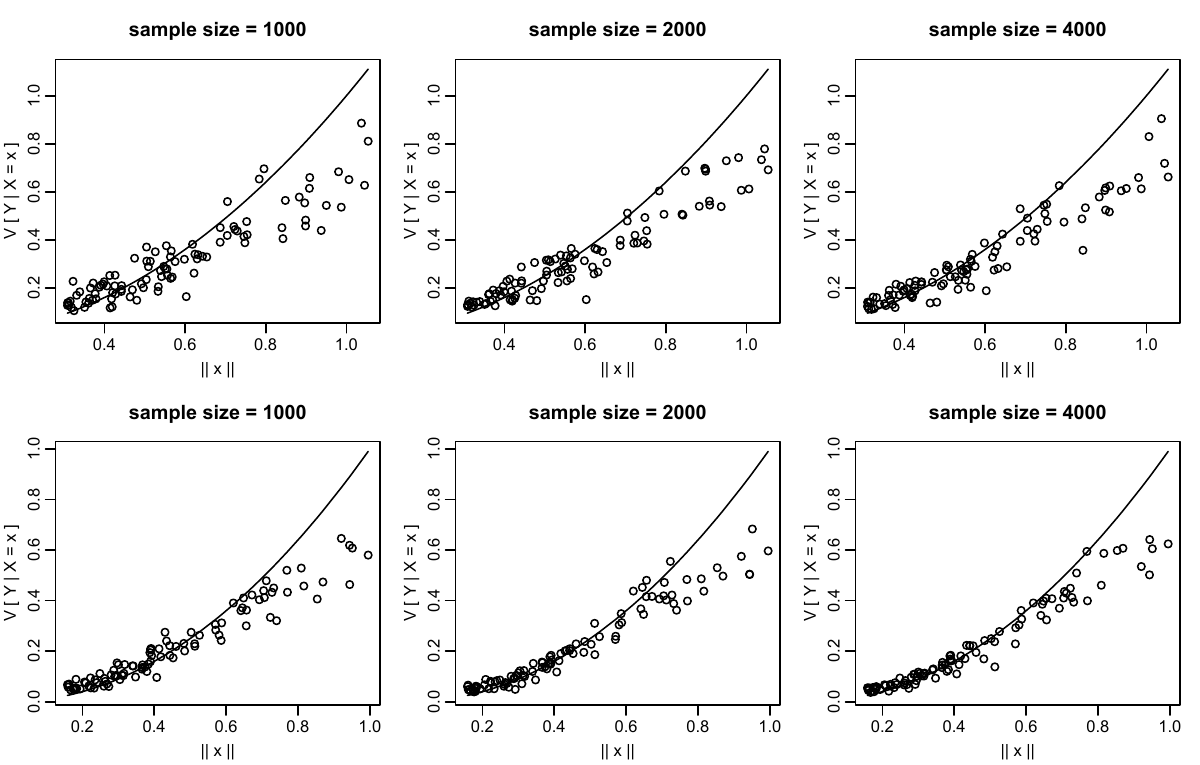}
\caption{Plots of actual conditional variance $ \mathbb{V}[ \mathbf{Y} \midil \mathbf{X} = \mathbf{x} ] $ (line) and its estimates $ \widehat{\mathbb{V}}_n^{(1)}[ \mathbf{Y} \midil \mathbf{X} = \mathbf{x} ] $ (points, first row) and $ \widehat{\mathbb{V}}_n^{(2)}[ \mathbf{Y} \midil \mathbf{X} = \mathbf{x} ] $ (points, second row) for different sample sizes. The covariate is a standard Brownian motion.}
\label{fig:condvarmlevar}
\end{figure}

\begin{figure}
\centering
\includegraphics[width=1\linewidth]{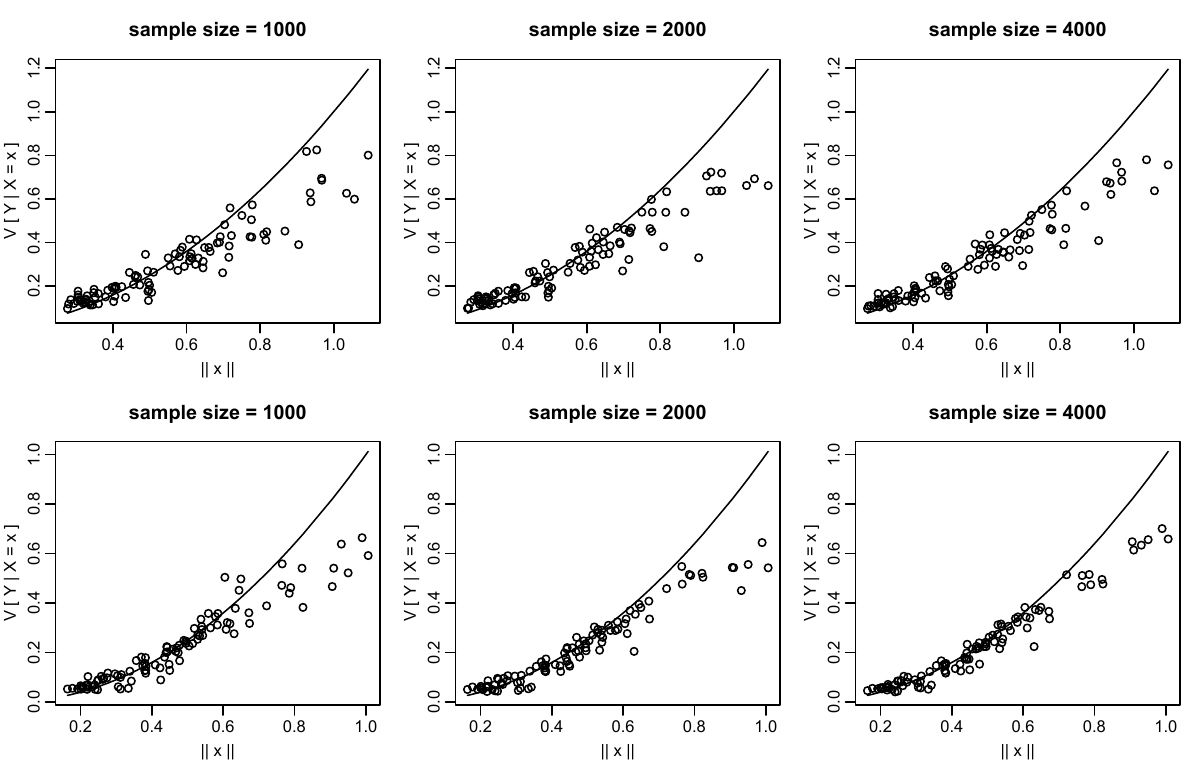}
\caption{Plots of actual conditional variance $ \mathbb{V}[ \mathbf{Y} \midil \mathbf{X} = \mathbf{x} ] $ (line) and its estimates $ \widehat{\mathbb{V}}_n^{(1)}[ \mathbf{Y} \midil \mathbf{X} = \mathbf{x} ] $ (points, first row) and $ \widehat{\mathbb{V}}_n^{(2)}[ \mathbf{Y} \midil \mathbf{X} = \mathbf{x} ] $ (points, second row) for different sample sizes. The covariate is a fractional Brownian motion with Hurst index 0.8.}
\label{fig:condvarmlevarfrac}
\end{figure}

As our first example, we consider $ \mathbf{Y} $ following a normal distribution with mean zero and variance $ \| \mathbf{X} \|_2^2 $. We consider two distributions for $ \mathbf{X} $, namely, the standard Brownian motion and the fractional Brownian motion with Hurst index 0.8. We want to estimate the conditional variance $ \mathbb{V}[ \mathbf{Y} \midil \mathbf{X} = \mathbf{x} ] $, which we do in two ways. In the first case, we estimate $ \mathbb{V}[ \mathbf{Y} \midil \mathbf{X} = \mathbf{x} ] $ by
\begin{align*}
& \widehat{\mathbb{V}}_n^{(1)}[ \mathbf{Y} \midil \mathbf{X} = \mathbf{x} ]
= \sum_{i=1}^{n} \left[ \mathbf{Y}_i - \left( \sum_{i=1}^{n} \mathbf{Y}_i W_{i,n}( \mathbf{x} ) \right) \right]^2 W_{i,n}( \mathbf{x} ) , \\
\intertext{where}
& W_{i,n}( \mathbf{x} )
= \frac{K\left( (h_n^*)^{-1} d( \mathbf{x}, \mathbf{X}_i ) \right)}{\sum_{i=1}^{n} K\left( (h_n^*)^{-1} d( \mathbf{x}, \mathbf{X}_i ) \right)} .
\end{align*}
So, this estimate belongs to the class of estimates described in \autoref{ex:2.2}.
In the second case, we estimate $ \mathbb{V}[ \mathbf{Y} \midil \mathbf{X} = \mathbf{x} ] $ using the weighted maximum likelihood procedure described in \autoref{ex:2.3}, with the conditional density of $ \mathbf{Y} $ given $ \mathbf{X} $ being the density of the normal random variable with mean zero and variance $ \| \mathbf{X} \|_2^2 $. In this case, our estimate turns out to be
\begin{align*}
\widehat{\mathbb{V}}_n^{(2)}[ \mathbf{Y} \midil \mathbf{X} = \mathbf{x} ]
= \sum_{i=1}^{n} \mathbf{Y}_i^2 W_{i,n}( \mathbf{x} ) .
\end{align*}
We randomly generate 100 values of $ \mathbf{x} $ from the distribution of $ \mathbf{X} $ and compute the estimates $ \widehat{\mathbb{V}}_n^{(1)}[ \mathbf{Y} \midil \mathbf{X} = \mathbf{x} ] $ and $ \widehat{\mathbb{V}}_n^{(2)}[ \mathbf{Y} \midil \mathbf{X} = \mathbf{x} ] $ for each of them.
We plot the estimates $ \widehat{\mathbb{V}}_n^{(1)}[ \mathbf{Y} \midil \mathbf{X} = \mathbf{x} ] $ and $ \widehat{\mathbb{V}}_n^{(2)}[ \mathbf{Y} \midil \mathbf{X} = \mathbf{x} ] $ along with the actual $ \mathbb{V}[ \mathbf{Y} \midil \mathbf{X} = \mathbf{x} ] $ values for different sample sizes against the values of $ \| \mathbf{x} \| $ in \autoref{fig:condvarmlevar}, where $ \mathbf{X} $ is a standard Brownian motion. When $ \mathbf{X} $ is a fractional Brownian motion with Hurst index 0.8, we plot the estimates along with the actual values in \autoref{fig:condvarmlevarfrac}.
from the plots, we observe that the two estimates $ \widehat{\mathbb{V}}_n^{(1)}[ \mathbf{Y} \midil \mathbf{X} = \mathbf{x} ] $ and $ \widehat{\mathbb{V}}_n^{(2)}[ \mathbf{Y} \midil \mathbf{X} = \mathbf{x} ] $ have no noticeable differences.
Also, there appears to be some underestimation when the value of $ \mathbb{V}[ \mathbf{Y} \midil \mathbf{X} = \mathbf{x} ] $ is high. This is due to the fact that the $ \mathbb{V}[ \mathbf{Y} \midil \mathbf{X} = \mathbf{X}_i ] $ values for $ \mathbf{X}_i $ lying in a neighborhood of $ \mathbf{x} $ tend to be smaller than $ \mathbb{V}[ \mathbf{Y} \midil \mathbf{X} = \mathbf{x} ] $, and the kernel estimate is based on those $ \mathbf{X}_i $ and their corresponding $ \mathbf{Y}_i $ values.
We also observe that the deviations of the estimated values from the actual values are marginally less when the covariate is a fractional Brownian motion with Hurst index 0.8, compared to the case where the covariate is a standard Brownian motion.
This may be due to the fact that the distribution of $ \| \mathbf{X} \|_2 $ is more concentrated at lower values when $ \mathbf{X} $ is a fractional Brownian motion with Hurst index 0.8 compared to the distribution of the same when $ \mathbf{X} $ is a standard Brownian motion.

\begin{figure}
\centering
\includegraphics[width=1\linewidth]{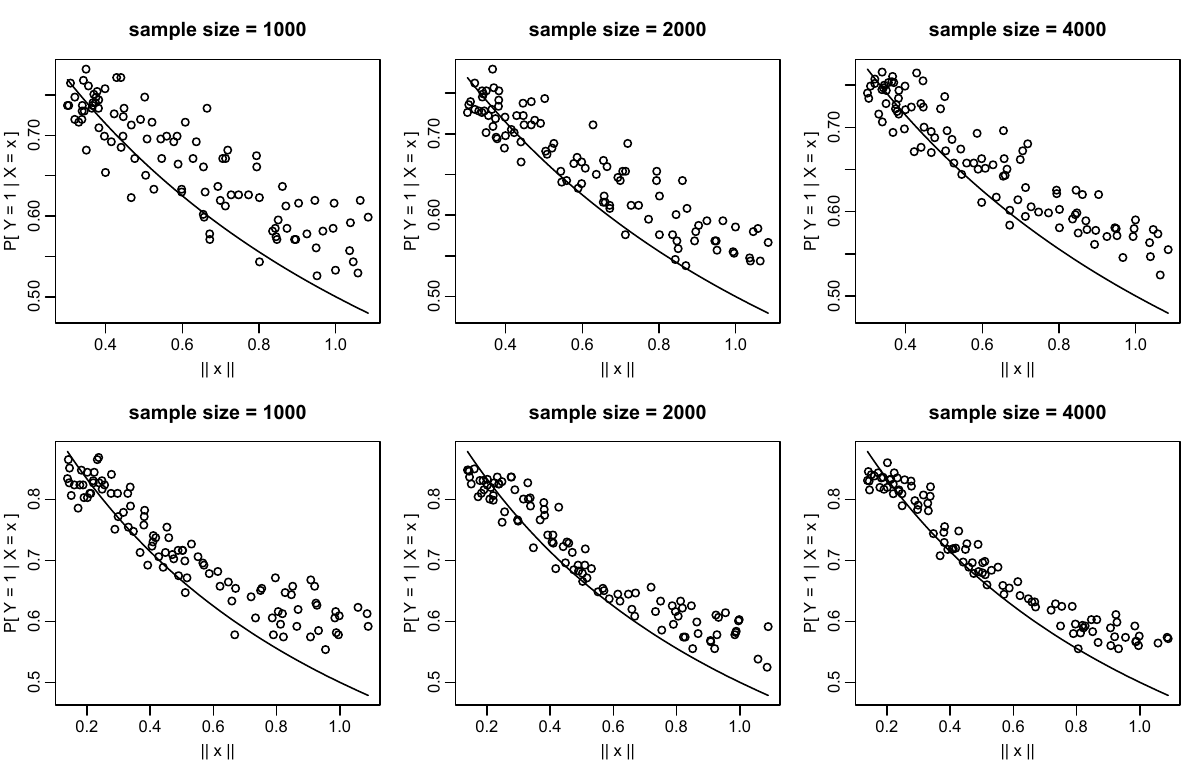}
\caption{Plots of estimated probability $ \widehat{\mathbb{P}}_n[ \mathbf{Y} = 1 \midil \mathbf{X} = \mathbf{x} ] $ (points) and actual probability $ \mathbb{P}[ \mathbf{Y} = 1 \midil \mathbf{X} = \mathbf{x} ] $ (line). The covariate is the standard Brownian motion in the first row, and the fractional Brownian motion with Hurst index 0.8 in the second row.}
\label{fig:mleprob}
\end{figure}

In the second example, we take $ \mathbf{Y} $ to be a Bernoulli random variable, with $ \mathbb{P}[ \mathbf{Y} = 1 \midil \mathbf{X} = \mathbf{x} ] = 1 - \mathbb{P}[ \mathbf{Y} = 0 \midil \mathbf{X} = \mathbf{x} ] = 1 / ( 1 + \| \mathbf{X} \|_2 ) $. Here our parameter of interest is $ \mathbb{P}[ \mathbf{Y} = 1 \midil \mathbf{X} = \mathbf{x} ] $, and we estimate it using the weighted maximum likelihood procedure described in \autoref{ex:2.3}, while employing the adaptive choice of the bandwidth. We again consider two distributions of $ \mathbf{X} $, namely the standard Brownian motion and the fractional Brownian motion with Hurst index 0.8, randomly generate 100 values of $ \mathbf{x} $ from the distribution of $ \mathbf{X} $ and plot the estimated values and the actual probabilities against the values of $ \| \mathbf{x} \| $ in \autoref{fig:mleprob} for three different sample sizes.
The improvement in accuracy of the estimate over the sample sizes is noticeable.
We also observe that there appears to be some overestimation for small values of $ \mathbb{P}[ \mathbf{Y} = 1 \midil \mathbf{X} = \mathbf{x} ] $, which is due to the fact that values of $ \mathbb{P}[ \mathbf{Y} = 1 \midil \mathbf{X} = \mathbf{X}_i ] $ for $ \mathbf{X}_i $ lying in a neighborhood of $ \mathbf{x} $ tend to be larger than $ \mathbb{P}[ \mathbf{Y} = 1 \midil \mathbf{X} = \mathbf{x} ] $ in such a case.
Further, like in the first example, we observe that the deviations of the estimated values from the actual values are less when the covariate is a fractional Brownian motion with Hurst index 0.8, compared to the case where the covariate is a standard Brownian motion.

\begin{figure}
\centering
\includegraphics[width=1\linewidth]{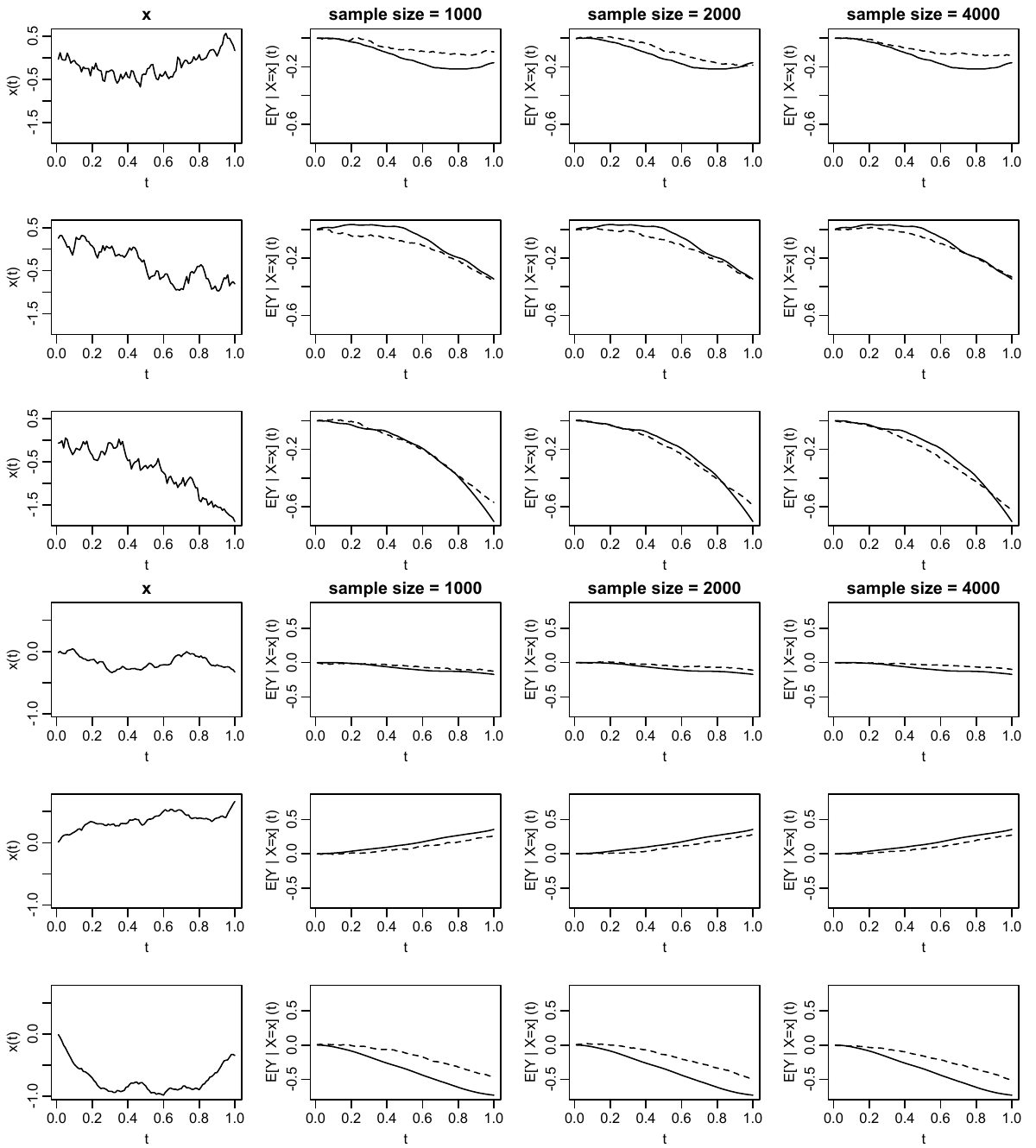}
\caption{Plots of estimated conditional mean curves $ \widehat{\mathbb{E}}_n[ \mathbf{Y} \midil \mathbf{X} = \mathbf{x} ] $ (dashed line) and actual conditional mean curves $ \mathbb{E}[ \mathbf{Y} \midil \mathbf{X} = \mathbf{x} ] $ (solid line). The covariate is a standard Brownian motion in the first three rows, and a fractional Brownian motion with Hurst index 0.8 in the last three rows.}
\label{fig:condmean}
\end{figure}

In the third example, we consider a functional response $ \mathbf{Y} $, defined by $ \mathbf{Y}( t ) = \int_{0}^{t} \mathbf{X}( t ) dt + \mathbf{E}( t ) $, where $ \mathbf{E}( \cdot ) $ is a Brownian motion independent of $ \mathbf{X}( \cdot ) $ with the covariance operator $ \mathbb{COV}( \mathbf{E}( s ), \mathbf{E}( t ) ) = 0.25 \min\{ s, t \} $. We want to estimate the conditional mean curve $ \mathbb{E}[ \mathbf{Y} \midil \mathbf{X} = \mathbf{x} ] $, for some fixed value $ \mathbf{x} $ of the covariate. We again consider two distributions of $ \mathbf{X} $, namely the standard Brownian motion and the fractional Brownian motion with Hurst index 0.8. In each case, we generate 3 random curves as values of $ \mathbf{x} $ and plot the adaptive estimates of the corresponding conditional means for different sample sizes in \autoref{fig:condmean}. In the first column, we have plotted the curves chosen as values of $ \mathbf{x} $. The first three rows in \autoref{fig:condmean} present the estimated conditional mean curves and the actual conditional mean curves for different sample sizes corresponding to the respective values of $ \mathbf{x} $ in the particular rows when the covariate is a standard Brownian motion. The last three rows in \autoref{fig:condmean} present the estimated conditional mean curves and the actual conditional mean curves for different sample sizes when the covariate is a fractional Brownian motion with Hurst index 0.8.
We observe that in all the cases, the estimates follow the actual curves closely.

\subsection{Demonstration in real data}
We now demonstrate the adaptive estimates of several regression parameters in the Tecator data. The Tecator data is a popular dataset available in the R package `caret'. This dataset contains the percentage values of moisture, fat and protein contents of 215 meat samples along with their absorbance spectra in the wavelength range 850--1050 nm measured by a Tecator spectroscope. The chemical contents of the meat samples are measured by analytical chemistry, which is expensive. The spectra of the samples are measured using a Tecator spectroscope, which is relatively inexpensive compared to the analytical chemistry method. So, it is economically beneficial to be able to predict the chemical contents of a sample from its spectra. Hence, we consider the fat and the protein content values as the response and the curve of the absorbance spectra as the covariate. We denote the percentage values of the fat and the protein contents as $ \mathbf{Y}_1 $ and $ \mathbf{Y}_2 $, respectively, and curve of the absorbance spectra as $ \mathbf{X} $. So, the covariate $ \mathbf{X} $ is a random function here, which we consider as a random element in the $ L_2 $ space. We consider 5 regression parameters of interest, namely $ \mathbb{E}[ \mathbf{Y}_1 \midil \mathbf{X} = \mathbf{x} ] $, $ \mathbb{E}[ \mathbf{Y}_2 \midil \mathbf{X} = \mathbf{x} ] $, $ \mathbb{VAR}[ \mathbf{Y}_1 \midil \mathbf{X} = \mathbf{x} ] $, $ \mathbb{VAR}[ \mathbf{Y}_2 \midil \mathbf{X} = \mathbf{x} ] $ and $ \mathbb{COR}[ \mathbf{Y}_1, \mathbf{Y}_2 \midil \mathbf{X} = \mathbf{x} ] $.
\begin{figure}
\centering
\includegraphics[width=1\linewidth]{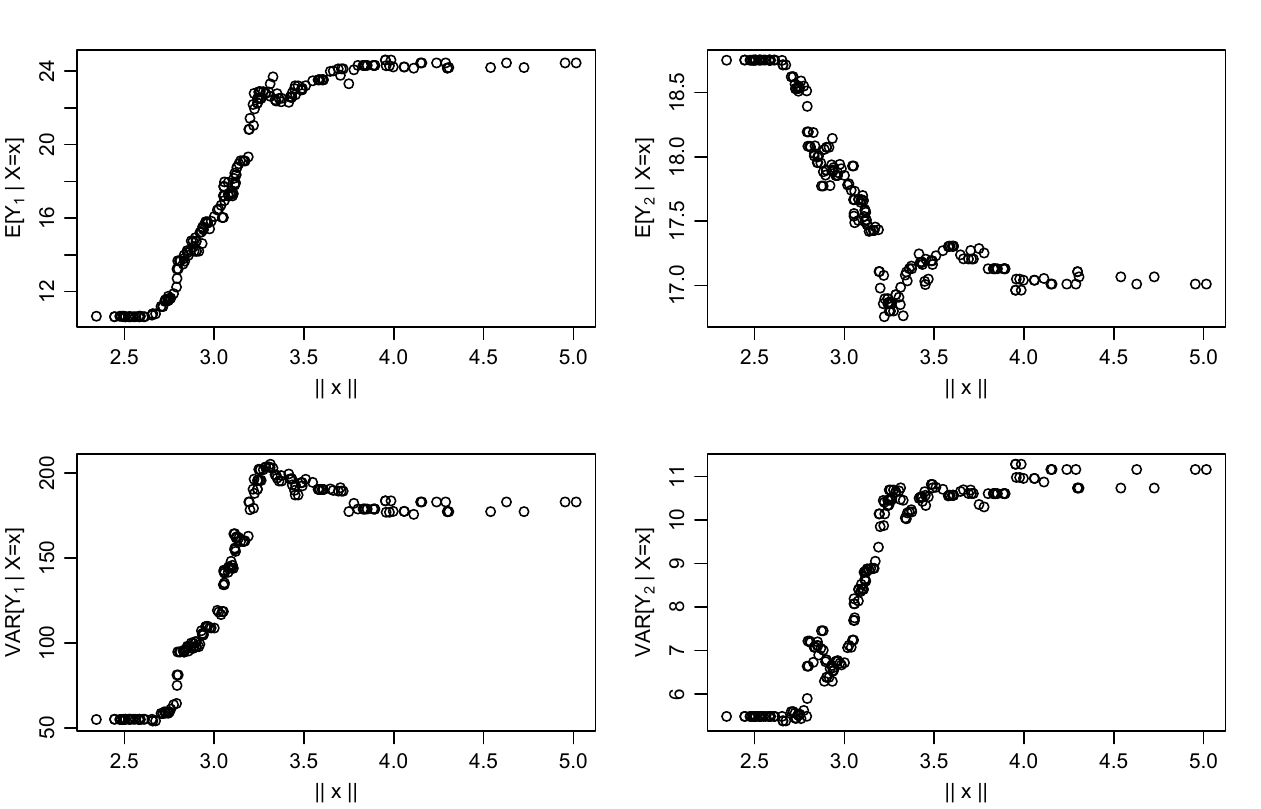}
\caption{Plots of adaptive estimates of $ \mathbb{E}[ \mathbf{Y}_1 \midil \mathbf{X} = \mathbf{x} ] $, $ \mathbb{E}[ \mathbf{Y}_2 \midil \mathbf{X} = \mathbf{x} ] $, $ \mathbb{VAR}[ \mathbf{Y}_1 \midil \mathbf{X} = \mathbf{x} ] $ and $ \mathbb{VAR}[ \mathbf{Y}_2 \midil \mathbf{X} = \mathbf{x} ] $ against the $ L_2 $-norm of $ \mathbf{x} $ in the Tecator data.}
\label{fig:adapttecator1}
\end{figure}
\begin{figure}
\centering
\includegraphics[width=1\linewidth]{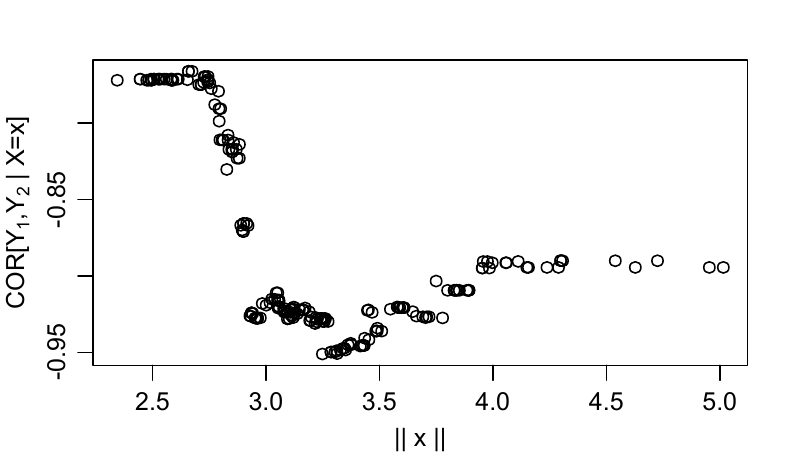}
\caption{Plot of adaptive estimates of $ \mathbb{COR}[ \mathbf{Y}_1, \mathbf{Y}_2 \midil \mathbf{X} = \mathbf{x} ] $ against the $ L_2 $-norm of $ x $ in the Tecator data.}
\label{fig:adapttecator2}
\end{figure}
We compute the adaptive estimates of all this parameters, where $ \mathbf{x} $ varies over all the sample curves of the absorbance spectra. We plot the adaptive estimates of $ \mathbb{E}[ \mathbf{Y}_1 \midil \mathbf{X} = \mathbf{x} ] $, $ \mathbb{E}[ \mathbf{Y}_2 \midil \mathbf{X} = \mathbf{x} ] $, $ \mathbb{VAR}[ \mathbf{Y}_1 \midil \mathbf{X} = \mathbf{x} ] $ and $ \mathbb{VAR}[ \mathbf{Y}_2 \midil \mathbf{X} = \mathbf{x} ] $ against the $ L_2 $ norm of $ \mathbf{x} $ in \autoref{fig:adapttecator1}, and the adaptive estimate of $ \mathbb{COR}[ \mathbf{Y}_1, \mathbf{Y}_2 \midil \mathbf{X} = \mathbf{x} ] $ against the $ L_2 $ norm of $ \mathbf{x} $ in \autoref{fig:adapttecator2}. The clear patterns of variation of the regression parameters over the covariate values are noticeable in each of the plots.

\begin{figure}
\centering
\includegraphics[width=1\linewidth]{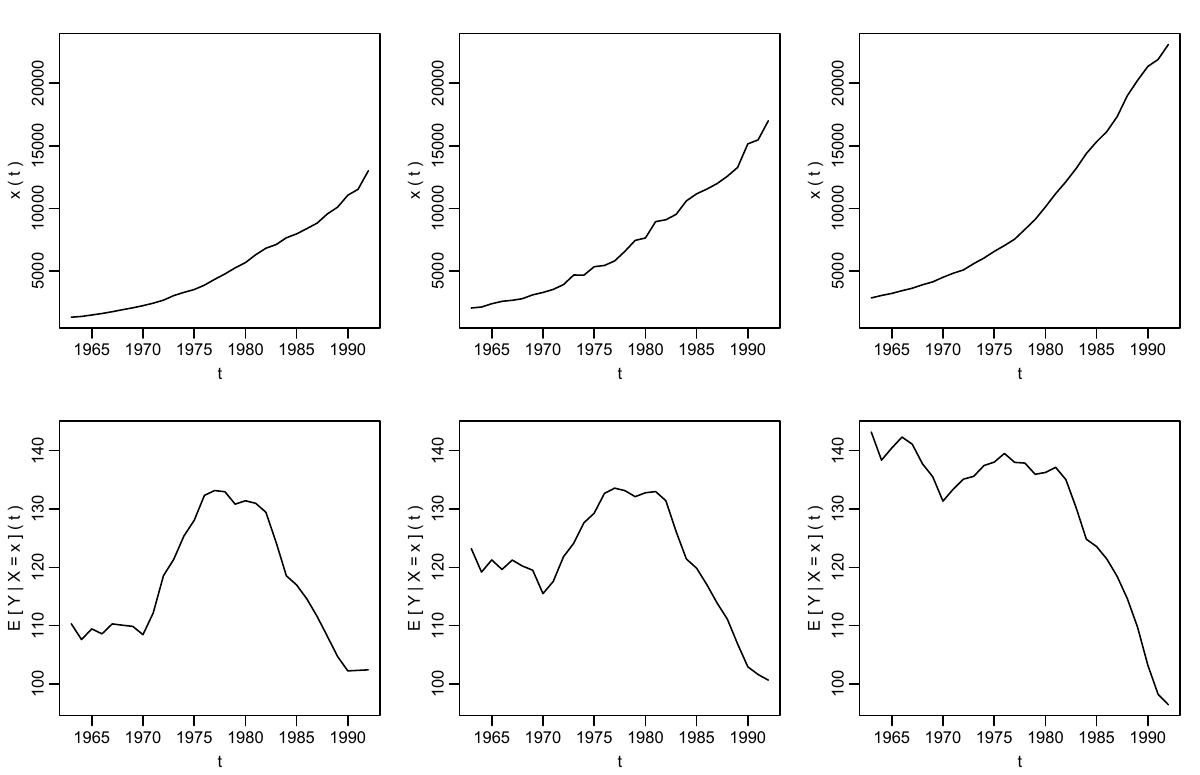}
\caption{Plots of adaptive estimates of $ \mathbb{E}[ \mathbf{Y} \midil \mathbf{X} = \mathbf{x} ] $ for 3 values of $ \mathbf{x} $ in the Cigar data.}
\label{fig:adaptcigar}
\end{figure}

Next, we demonstrate the adaptive estimates of the conditional mean in another dataset, where both the response and the covariate are random functions. The dataset we consider is the Cigar data, which is available in the `Ecdat' package in R. This dataset contains information about cigarette sales in packs per capita, per capita net disposable income (NDI) and other economic parameters in 46 states in the USA over a 30 years period from 1963 to 1992. We consider the curve of NDI over 30 years as the covariate $ \mathbf{X} $, and the curve of cigarette sales over 30 years as the response $ \mathbf{Y} $. So, both the response and the covariate in this setup are random functions, and our sample size is 46. We choose 3 sample covariate curves as values of $ \mathbf{x} $, and compute the adaptive estimates of $ \mathbb{E}[ \mathbf{Y} \midil \mathbf{X} = \mathbf{x} ] $ for these 3 values of $ \mathbf{x} $. We plot the estimated curves along with the respective covariate curves in \autoref{fig:adaptcigar}, where the first row contains the plots of the 3 curves chosen as values of $ \mathbf{x} $, and the second row contains the plots of the corresponding adaptive estimates of $ \mathbb{E}[ \mathbf{Y} \midil \mathbf{X} = \mathbf{x} ] $. The 3 estimated curves reflect the variation of $ \mathbb{E}[ \mathbf{Y} \midil \mathbf{X} = \mathbf{x} ] $ over $ \mathbf{x} $.

\section{Concluding remarks} \label{sec:5}
In this paper, we have derived the optimum convergence rate for a wide class of kernel regression estimates when the covariate as well as the response may be infinite dimensional. It is shown that the convergence rates of such estimates do not depend on the dimension of the response, but they depend critically on the dimension of the covariate.
We have seen that, for a wide class of covariates having infinite dimensional Gaussian distributions, the convergence rate is much slower than the optimum achievable rate for finite dimensional covariates. For instance, if the covariate is a real-valued continuous Gaussian Markov process in $ L_p[ 0, 1 ] $, the convergence rate is $ O( ( \log n )^{-\delta} ) $ for some $ \delta > 0 $. \autoref{thm:2} implies that if $ h_n^{2 \beta} n \phi( \mathbf{x}, h_n ) \longrightarrow 0 $ as $ n \longrightarrow \infty $, $ [ n \phi( \mathbf{x}, h_n ) ]^{1/2} c_n [ \widehat{\Theta}_n( \mathbf{x} ) - \Theta( \mathbf{x} ) ] $ converges in distribution to a Gaussian random element with zero mean as $ n \longrightarrow \infty $, where $ c_n = [ E_n^{(2)}( \mathbf{x} ) ]^{-1/2} \allowbreak E_n^{(1)}( \mathbf{x} ) $ is a sequence of positive numbers bounded and bounded away from 0. Note that this corresponds to an under-smoothed kernel estimate of $ \Theta( \mathbf{x} ) $. On the other hand, if $ h_n^{2 \beta} n \phi( \mathbf{x}, h_n ) \longrightarrow \infty $ as $ n \longrightarrow \infty $, which includes the case of our optimum bandwidth obtained in \autoref{thm:up}, we have $ h_n^{- \beta} [ \widehat{\Theta}_n( \mathbf{x} ) - \Theta( \mathbf{x} ) ] - h_n^{- \beta} \tilde{B}_n( \mathbf{x} ) \longrightarrow 0 $ \textit{in probability} as $ n \longrightarrow \infty $. Here, $ \tilde{B}_n( \mathbf{x} ) $ is a non-random deterministic object described at the beginning of \autoref{subsec:4_1}.

In \cite{ferraty2006nonparametric}, \cite{ferraty2006estimating,ferraty2010rate} and \cite{chaouch2013nonparametric,chaouch2015vector}, asymptotic properties of nonparametric regression estimates of different parameters other than the mean of the conditional distribution of the response were investigated. However, they only considered finite dimensional responses, and they did not investigate the problem of optimum convergence rates of nonparametric regression estimates.

The problem of slow convergence rate of the regression estimates with infinite dimensional covariates that has been derived in this paper may be coped with using an appropriate dimension reduction procedure on the covariate. Some procedures for such dimension reduction for infinite dimensional covariates available in the literature are the uses of functional sliced inverse regression \citep{ferre2003functional,ferre2005smoothed}, functional average derivative regression \citep{ferraty2011estimation} and distance correlation maximization \citep{vepakomma2016supervised}. If the covariate with the reduced dimension is adequate for regression analysis, the new small ball probability function in the reduced covariate space will lead to better convergence rates.

\section{Proofs and mathematical details} \label{sec:6}

\subsection{Small ball probabilities of non-Gaussian processes} \label{sec:7}

In Propositions \ref{prop:3}, \ref{prop:1} and \ref{prop:2} below, we consider two random elements $ \mathbf{T} $ and $ \mathbf{G} $, and define $ \phi_{\mathbf{T}}( \mathbf{t}, h ) = \mathbb{P}\left[ \| \mathbf{T} - \mathbf{t} \| \le h \right] $ and $ \phi_{\mathbf{G}}( \mathbf{g}, h ) = \mathbb{P}\left[ \| \mathbf{G} - \mathbf{g} \| \le h \right] $, where $ \mathbf{t} $ and $ \mathbf{g} $ are some fixed elements and $ h > 0 $.

\begin{proposition} \label{prop:3}
Let $ \mathcal{B}_1 $ and $ \mathcal{B}_2 $ be separable Banach spaces, and $ f( \cdot ) : \mathcal{B}_2 \longrightarrow \mathcal{B}_1 $ be a function such that for any $ \mathbf{u} \in \mathcal{B}_2 $, there exist constants $ r, s > 0 $, which may depend on $ \mathbf{u} $, such that for any $ \mathbf{v} \in \mathcal{B}_2 $ sufficiently close to $ \mathbf{u} $, we have
$ r \| \mathbf{v} - \mathbf{u} \|
\le \| f( \mathbf{v} ) - f( \mathbf{u} ) \| 
\le s \| \mathbf{v} - \mathbf{u} \| $.
If $ \mathbf{T} $ and $ \mathbf{G} $ are random elements with $ \mathbf{T} = f( \mathbf{G} ) $, and the small ball probability of $ \mathbf{G} $ satisfies the bounds described in \eqref{eq:1}, then similar bounds also hold for $ \mathbf{T} $.
\end{proposition}
\begin{proof}
Under the assumptions of the proposition, $ f( \cdot ) $ is a one-to-one function. Let $ \mathbf{t} $ be an element in the range of $ f( \cdot ) $. Then, $ \mathbf{t} = f( \mathbf{g} ) $ for some $ \mathbf{g} $. Consequently, for some positive constants $ r $ and $ s $, which may depend on $ \mathbf{g} $, we have for all sufficiently small $ h $,
\begin{align}
& \mathbb{P}\left[ s \| \mathbf{G} - \mathbf{g} \| \le h \right]
\le
\mathbb{P}\left[ \| f( \mathbf{G} ) - f( \mathbf{g} ) \| \le h \right]
\le
\mathbb{P}\left[ r \| \mathbf{G} - \mathbf{g} \| \le h \right]
\nonumber\\
& \iff
\phi_{\mathbf{G}}\left( \mathbf{g}, \frac{h}{s} \right)
\le
\phi_{\mathbf{T}}( \mathbf{x}, h )
\le
\phi_{\mathbf{G}}\left( \mathbf{g}, \frac{h}{r} \right) .
\label{eq:prop3:1}
\end{align}
The proof follows by applying the bounds in \eqref{eq:1} in \eqref{eq:prop3:1}.
\end{proof}

Let $ \mathbf{G} $ be a Gaussian process whose small ball probability $ \phi_{\mathbf{G}}( \mathbf{g}, h ) $ satisfies the bounds in \eqref{eq:1} for sufficiently small $ h $, so that
\begin{align*}
C_1 h^{t_1} \exp\left[ -C_2 ( 1 / h )^{t_2} ( \log ( 1 / h ) )^{t_3} \right]
& \le \phi_{\mathbf{G}}( \mathbf{g}, h ) \\
& \le C_3 h^{t_4} \exp\left[ -C_4 ( 1 / h )^{t_2} ( \log ( 1 / h ) )^{t_3} \right]
\end{align*}
as $ h \longrightarrow 0^+ $. Here, $ C_1, C_2, C_3, C_4 > 0 $ and $ t_1, t_2, t_3, t_4 \ge 0 $ are appropriate constants, all of which, except $ C_1 $, are independent of $ \mathbf{g} $. $ C_1 $ may or may not depend on $ \mathbf{g} $, but if it depends on $ \mathbf{g} $ then $ C_1 = C_1' \exp[ - ( 1 / 2 ) \| \mathbf{g} \|^2 ] $ for some positive constant $ C_1' $. Also, either $ t_2 > 0 $, or $ t_3 > 1 $ with $ C_2 = C_4 $.

In \autoref{prop:1} and \autoref{prop:2} below, we derive the bounds on the small ball probabilities of some non-Gaussian processes. There, we shall assume $ C_1 = C_1' \exp[ - ( 1 / 2 ) \| \mathbf{g} \|^2 ] $ for some positive constant $ C_1' $. Since $ C_1' \ge C_1' \exp[ - ( 1 / 2 ) \| \mathbf{g} \|^2 ] $ for all $ \mathbf{g} $, establishing the lower bound of the small ball probability, when $ C_1 = C_1' \exp[ - ( 1 / 2 ) \| \mathbf{g} \|^2 ] $, also gives an appropriate lower bound when $ C_1 $ does not depend on $ \mathbf{g} $.

\begin{proposition} \label{prop:1}
Let $ \mathbf{T} = \mathbf{G} / \mathbf{U} $, where $ \mathbf{G} $ is a Gaussian process whose small ball probability satisfies the bounds in \eqref{eq:1}, and $ \mathbf{U} $ is a bounded positive random variable independent of $ \mathbf{G} $. Then, the small ball probability of $ \mathbf{T} $ also satisfies the bounds in \eqref{eq:1}.
\end{proposition}
\begin{proof}
Note that
\begin{align}
\phi_{\mathbf{T}}( \mathbf{t}, h )
& = \mathbb{P}\left[ \left\| \mathbf{G} - \mathbf{t} \mathbf{U} \right\| \le h \mathbf{U} \right] 
= \mathbb{E}\left[ \phi_{\mathbf{G}}\left( \mathbf{t} \mathbf{U}, \; h \mathbf{U} \right) \right] .
\label{eq:prop1:1}
\end{align}
Let $ 0 \le \mathbf{U} \le u_0 $ for some $ u_0 > 0 $.
Recall from \eqref{eq:mh} that
$ m( h ) = C_2 ( 1 / h )^{t_2} ( \log ( 1 / h ) )^{t_3} $ for $ 0 < h < 1 $.
Since $ m( h u_0 ) \le m( h \mathbf{U} ) $ for all $ h > 0 $, we have
\begin{align*}
\phi_{\mathbf{G}}\left( \mathbf{t} \mathbf{U}, \; h \mathbf{U} \right)
& \le C_3 ( h \mathbf{U} )^{t_4} \exp\left[ -(C_4 / C_2) m( h \mathbf{U} ) \right] \\
& \le C_3 ( h u_0 )^{t_4} \exp\left[ -(C_4 / C_2) m( h u_0 ) \right] \\
& = C_3 u_0^{t_4} h^{t_4} \exp\left[ -C_4 \left( \frac{1}{u_0} \right)^{t_2} \left( 1 - \frac{\log u_0}{\log \frac{1}{h}} \right)^{t_3} \left( \frac{1}{h} \right)^{t_2} \left( \log \frac{1}{h} \right)^{t_3} \right] \\
& \le C_3 u_0^{t_4} h^{t_4} \exp\left[ - \frac{C_4}{2} \left( \frac{1}{u_0} \right)^{t_2} \left( \frac{1}{h} \right)^{t_2} \left( \log \frac{1}{h} \right)^{t_3} \right]
\end{align*}
for all sufficiently small $ h $. Hence, for all sufficiently small $ h $,
\begin{align}
\mathbb{E}\left[ \phi_{\mathbf{G}}\left( \mathbf{t} \mathbf{U}, \; h \mathbf{U} \right) \right]
\le
C_3 u_0^{t_4} h^{t_4} \exp\left[ - \frac{C_4}{2} \left( \frac{1}{u_0} \right)^{t_2} \left( \frac{1}{h} \right)^{t_2} \left( \log \frac{1}{h} \right)^{t_3} \right] .
\label{eq:prop1:2}
\end{align}

Now, if $ \mathbf{U} $ is a degenerate positive random variable, i.e., $ \mathbb{P}[ \mathbf{U} = u_0 ] = 1 $, then the lower bound of $ \phi_{\mathbf{G}}\left( \mathbf{t} \mathbf{U}, \; h \mathbf{U} \right) $ trivially satisfies \eqref{eq:1}. So, we assume that $ \mathbf{U} $ is non-degenerate, and $ \mathbb{P}[ l_0 \le \mathbf{U} < u_0 ] > 0 $ for some $ l_0 > 0 $.
We consider the case where the constant $ C_1 $ depends on the center of the small ball probability of $ \mathbf{G} $. The case when $ C_1 $ does not depend on the center of the small ball probability of $ \mathbf{G} $ can be covered similarly.
So, we have
\begin{align}
& \mathbb{E}\left[ \phi_{\mathbf{G}}\left( \mathbf{t} \mathbf{U}, \; h \mathbf{U} \right) \right] \nonumber\\
& \ge \mathbb{E}\left[ C_1' \exp[ - ( 1 / 2 ) \| \mathbf{t} \mathbf{U} \|^2 ] ( h \mathbf{U} )^{t_1} \exp\left[ - m( h \mathbf{U} ) \right] \right] \nonumber\\
& \ge \mathbb{E}\left[ C_1' \exp[ - ( 1 / 2 ) \mathbf{U}^2 \| \mathbf{t} \|^2 ] ( h \mathbf{U} )^{t_1} \exp\left[ - m( h \mathbf{U} ) \right] \mathbb{I}( \mathbf{U} \ge l_0 ) \right] \nonumber\\
& \ge \mathbb{E}\left[ C_1' \exp[ - ( 1 / 2 ) u_0^2 \| \mathbf{t} \|^2 ] ( h \mathbf{U} )^{t_1} \exp\left[ - m( h \mathbf{U} ) \right] \mathbb{I}( \mathbf{U} \ge l_0 ) \right] \nonumber\\
& \ge C_1' \exp[ - ( 1 / 2 ) u_0^2 \| \mathbf{t} \|^2 ] l_0^{t_1} h^{t_1} \exp\left[ - m( h l_0 ) \right] \nonumber\\
& = C_1' \exp\left[ - \frac{1}{2} u_0^2 \| \mathbf{t} \|^2 \right] l_0^{t_1} h^{t_1} \exp\left[ - \left( \frac{1}{l_0} \right)^{t_2} \left( 1 - \frac{\log l_0}{\log \frac{1}{h}} \right)^{t_3} \left( \frac{1}{h} \right)^{t_2} \left( \log \frac{1}{h} \right)^{t_3} \right] \nonumber\\
& \ge C_1' \exp[ - ( 1 / 2 ) u_0^2 \| \mathbf{t} \|^2 ] l_0^{t_1} h^{t_1} \exp\left[ - 2 \left( \frac{1}{l_0} \right)^{t_2} \left( \frac{1}{h} \right)^{t_2} \left( \log \frac{1}{h} \right)^{t_3} \right]
\label{eq:prop1:3}
\end{align}
for all sufficiently small $ h $. The proof is completed combining \eqref{eq:prop1:1}, \eqref{eq:prop1:2} and \eqref{eq:prop1:3}.
\end{proof}

Note that if $ \mathbf{T} $ is an infinite dimensional \textit{t}-process with degree $ k $, it can be expressed as $ \mathbf{T} = \mathbf{G} / \sqrt{\boldsymbol{\chi} / k} $, where $ \mathbf{G} $ is an infinite dimensional Gaussian process, $ \boldsymbol{\chi} $ follows a $ \chi^2 $ distribution with degree of freedom $ k $, and $ \boldsymbol{\chi} $ is independent of $ \mathbf{G} $. In the proposition below, we establish the bounds for the small ball probability of an infinite dimensional \textit{t}-process $ \mathbf{T} $.
\begin{proposition} \label{prop:2}
Let $ \mathbf{T} $ be an infinite dimensional \textit{t}-process in some normed vector space with corresponding Gaussian process $ \mathbf{G} $, and the small ball probability of $ \mathbf{G} $ satisfies the bounds in \eqref{eq:1} with $ t_2 > 0 $. Then, the small ball probability of $ \mathbf{T} $ also satisfies the bounds in \eqref{eq:1}.
\end{proposition}
\begin{proof}
We have
\begin{align}
\phi_{\mathbf{T}}( \mathbf{t}, h )
& = \mathbb{P}\left[ \left\| \mathbf{G} - \mathbf{t} \sqrt{\boldsymbol{\chi} / k} \right\| \le h \sqrt{\boldsymbol{\chi} / k} \right] \nonumber\\
& = \mathbb{E}\left[ \mathbb{P}\left[ \left\| \mathbf{G} - \mathbf{t} \sqrt{\boldsymbol{\chi} / k} \right\| \le h \sqrt{\boldsymbol{\chi} / k} \;\middle\arrowvert\; \boldsymbol{\chi} \right] \right] \nonumber\\
& = \frac{1}{2^{\frac{k}{2}} \Gamma\left( \frac{k}{2} \right)} \int_{0}^{\infty} \phi_{\mathbf{G}}\left( \mathbf{t} \sqrt{\frac{u}{k}}, \; h \sqrt{\frac{u}{k}} \right) e^{- \frac{u}{2}} u^{\frac{k}{2} - 1} d u .
\label{eq:prop2:1}
\end{align}
Define
$ m_1( h ) = ( 1 / h )^{t_2} ( \log ( 1 / h ) )^{t_3} $ for $ 0 < h < 1 $.
Since $ t_2 > 0 $, $ m_1( h ) \longrightarrow \infty $ as $ h \longrightarrow 0^+ $.
Let
\begin{align}
t_5 = 1 + \frac{t_2}{2} .
\label{eq:prop2:2}
\end{align}
Define
\begin{align}
U( h ) = (m_1( h ))^{\frac{1}{t_5}} .
\label{eq:prop2:3}
\end{align}
Clearly, $ U( h ) \longrightarrow \infty $ as $ h \longrightarrow 0^+ $.
Also,
\begin{align}
h \sqrt{U( h )}
= h \left[ \left( \frac{1}{h} \right)^{t_2} \left( \log \frac{1}{h} \right)^{t_3} \right]^{\frac{1}{2 t_5}}
= h^{\frac{1}{t_5}} \left( \log \frac{1}{h} \right)^{\frac{t_3}{2 t_5}}
\longrightarrow 0 \text{ as } h \longrightarrow 0^+ .
\label{eq:prop2:4}
\end{align}
So, from \eqref{eq:1} and \eqref{eq:prop2:2}, \eqref{eq:prop2:3} and \eqref{eq:prop2:4}, we have for all sufficiently small $ h $ and for any $ u \le U( h ) $,
\begin{align}
& \phi_{\mathbf{G}}\left( \mathbf{t} \sqrt{\frac{u}{k}}, \; h \sqrt{\frac{u}{k}} \right) \nonumber\\
& \le C_3 \left( h \sqrt{\frac{u}{k}} \right)^{t_4} \exp\left[ -C_4 m_1\left( h \sqrt{\frac{u}{k}} \right) \right] \nonumber\\
& = \frac{C_3}{k^{\frac{t_4}{2}}} u^{\frac{t_4}{2}} h^{t_4} \exp\left[ -C_4 k^{\frac{t_2}{2}} u^{-\frac{t_2}{2}} \left( \frac{1}{h} \right)^{t_2} \left( \log \frac{1}{h} \right)^{t_3} \left( 1 + \frac{\log \sqrt{k}}{\log \frac{1}{h}} - \frac{\log \sqrt{u}}{\log \frac{1}{h}} \right)^{t_3} \right]
\nonumber\\
& \le \frac{C_3}{k^{\frac{t_4}{2}}} u^{\frac{t_4}{2}} h^{t_4} \exp\left[ -C_4 k^{\frac{t_2}{2}} \left( \frac{1}{t_2 + 2} \right)^{t_3} \left( m_1( h ) \right)^\frac{1}{t_5} \right] ,
\label{eq:prop2:5}
\end{align}
since for all sufficiently small $ h $ and any $ u \le U( h ) $,
\begin{align*}
& 1 + \frac{\log \sqrt{k}}{\log \frac{1}{h}} - \frac{\log \sqrt{u}}{\log \frac{1}{h}} > \frac{1}{t_2 + 2} .
\end{align*}
Hence, from \eqref{eq:prop2:1} and \eqref{eq:prop2:5}, we have for all sufficiently small $ h $,
\begin{align}
& \phi_{\mathbf{T}}( \mathbf{t}, h ) \nonumber\\
& = \frac{1}{2^{\frac{k}{2}} \Gamma\left( \frac{k}{2} \right)} \int_{0}^{U( h )} \phi_{\mathbf{G}}\left( \mathbf{t} \sqrt{\frac{u}{k}}, \; h \sqrt{\frac{u}{k}} \right) e^{- \frac{u}{2}} u^{\frac{k}{2} - 1} d u \nonumber\\
& \quad + \frac{1}{2^{\frac{k}{2}} \Gamma\left( \frac{k}{2} \right)} \int_{U( h )}^{\infty} \phi_{\mathbf{G}}\left( \mathbf{t} \sqrt{\frac{u}{k}}, \; h \sqrt{\frac{u}{k}} \right) e^{- \frac{u}{2}} u^{\frac{k}{2} - 1} d u
\nonumber\\
& < \frac{1}{2^{\frac{k}{2}} \Gamma\left( \frac{k}{2} \right)} \frac{C_3}{k^{\frac{t_4}{2}}} \int_{0}^{U( h )} h^{t_4} \exp\left[ -C_4 k^{\frac{t_2}{2}} \left( \frac{1}{t_2 + 2} \right)^{t_3} \left( m_1( h ) \right)^{\frac{1}{t_5}} \right] e^{- \frac{u}{2}} u^{\frac{t_4 + k}{2} - 1} d u \nonumber\\
& \quad + \frac{1}{2^{\frac{k}{2}} \Gamma\left( \frac{k}{2} \right)} \int_{U( h )}^{\infty} \exp\left[ - \frac{1}{4} U( h ) \right] e^{- \frac{u}{4}} u^{\frac{k}{2} - 1} d u
\nonumber\\
& < \frac{1}{2^{\frac{k}{2}} \Gamma\left( \frac{k}{2} \right)} \frac{C_3}{k^{\frac{t_4}{2}}}  \left[ \int_{0}^{\infty} e^{- \frac{u}{2}} u^{\frac{t_4 + k}{2} - 1} d u \right] h^{t_4} \exp\left[ -C_4 k^{\frac{t_2}{2}} \left( \frac{1}{t_2 + 2} \right)^{t_3} \left( m_1( h ) \right)^{\frac{1}{t_5}} \right] \nonumber\\
& \quad + \frac{1}{2^{\frac{k}{2}} \Gamma\left( \frac{k}{2} \right)} \left[ \int_{0}^{\infty} e^{- \frac{u}{4}} u^{\frac{k}{2} - 1} d u \right] \exp\left[ - \frac{1}{4} ( m_1( h ) )^{\frac{1}{t_5}} \right]
\nonumber\\
& = \left( \frac{\Gamma\left( \frac{t_4 + k}{2} \right)}{\Gamma\left( \frac{k}{2} \right)} \left( \frac{2}{k} \right)^{\frac{t_4}{2}} C_3 \right) h^{t_4} \exp\left[ -C_4 k^{\frac{t_2}{2}} \left( \frac{1}{t_2 + 2} \right)^{t_3} \left( m_1( h ) \right)^{\frac{1}{t_5}} \right] \nonumber\\
& \quad
+ 2^{\frac{k}{2}} \exp\left[ - \frac{1}{4} ( m_1( h ) )^{\frac{1}{t_5}} \right]
\nonumber\\
& \le \left( \frac{\Gamma\left( \frac{t_4 + k}{2} \right)}{\Gamma\left( \frac{k}{2} \right)} \left( \frac{2}{k} \right)^{\frac{t_4}{2}} C_3 + 2^{\frac{k}{2}} \right) \times \nonumber\\
& \qquad \exp\left[ - \min\left\{ C_4 k^{\frac{t_2}{2}} \left( \frac{1}{t_2 + 2} \right)^{t_3} ,\, \frac{1}{4} \right\} \left( \frac{1}{h} \right)^{\frac{t_2}{t_5}} \left( \log \frac{1}{h} \right)^{\frac{t_3}{t_5}} \right] .
\label{eq:prop2:6}
\end{align}
We now proceed to find a lower bound for $ \phi_{\mathbf{T}}( \mathbf{t}, h ) $.
From \eqref{eq:1}, \eqref{eq:prop2:2}, \eqref{eq:prop2:3} and \eqref{eq:prop2:4}, we get that for all sufficiently small $ h $ and for any $ U( h ) \le u \le 2 U( h ) $,
\begin{align}
& \phi_{\mathbf{G}}\left( \mathbf{t} \sqrt{\frac{u}{k}}, \; h \sqrt{\frac{u}{k}} \right) \nonumber\\
& \ge C_1' \exp\left[ - \frac{1}{2} \left\| \mathbf{t} \sqrt{\frac{u}{k}} \right\|^2 \right] \left( h \sqrt{\frac{u}{k}} \right)^{t_1} \exp\left[ -C_2 m_1\left( h \sqrt{\frac{u}{k}} \right) \right] \nonumber\\
& = \frac{C_1'}{k^{\frac{t_1}{2}}} u^{\frac{t_1}{2}} h^{t_1} \times \nonumber\\
& \qquad
\exp\left[ - u \frac{\left\| \mathbf{t} \right\|^2}{2 k} -C_2 k^{\frac{t_2}{2}} u^{-\frac{t_2}{2}} \left( \frac{1}{h} \right)^{t_2} \left( \log \frac{1}{h} \right)^{t_3} \left( 1 + \frac{\log \sqrt{k}}{\log \frac{1}{h}} - \frac{\log \sqrt{u}}{\log \frac{1}{h}} \right)^{t_3} \right] 
\nonumber\\
& \ge \frac{C_1'}{k^{\frac{t_1}{2}}} u^{\frac{t_1}{2}} h^{t_1} \exp\left[ - \frac{\left\| \mathbf{t} \right\|^2}{k} (m_1( h ))^{\frac{1}{t_5}} -C_2 k^{\frac{t_2}{2}} \left( \frac{2}{t_5} \right)^{t_3} \left( m_1( h ) \right)^{1 - \frac{t_2}{2 t_5}} \right] 
\nonumber\\
& = \frac{C_1'}{k^{\frac{t_1}{2}}} u^{\frac{t_1}{2}} h^{t_1} \exp\left[ - \left( \frac{\left\| \mathbf{t} \right\|^2}{k} + C_2 k^{\frac{t_2}{2}} \left( \frac{2}{t_5} \right)^{t_3} \right) \left( m_1( h ) \right)^{\frac{1}{t_5}} \right] ,
\label{eq:prop2:7}
\end{align}
since for all sufficiently small $ h $ and any $ U( h ) \le u $,
\begin{align*}
& 1 + \frac{\log \sqrt{k}}{\log \frac{1}{h}} - \frac{\log \sqrt{u}}{\log \frac{1}{h}} < \frac{2}{t_5} .
\end{align*}
From \eqref{eq:prop2:1} and \eqref{eq:prop2:7}, we have for all sufficiently small $ h $,
\begin{align}
& \phi_{\mathbf{T}}( \mathbf{t}, h ) \nonumber\\
& \ge \frac{1}{2^{\frac{k}{2}} \Gamma\left( \frac{k}{2} \right)} \int_{U( h )}^{2 U( h )} \phi_{\mathbf{G}}\left( \mathbf{t} \sqrt{\frac{u}{k}}, \; h \sqrt{\frac{u}{k}} \right) e^{- \frac{u}{2}} u^{\frac{k}{2} - 1} d u
\nonumber\\
& \ge \frac{1}{2^{\frac{k}{2}} \Gamma\left( \frac{k}{2} \right)} \frac{C_1'}{k^{\frac{t_1}{2}}} \times \nonumber\\
& \qquad
\int_{U( h )}^{2 U( h )} h^{t_1} \exp\left[ - \left( \frac{\left\| \mathbf{t} \right\|^2}{k} + C_2 k^{\frac{t_2}{2}} \left( \frac{2}{t_5} \right)^{t_3} \right) \left( m_1( h ) \right)^{\frac{1}{t_5}} \right] e^{- \frac{u}{2}} u^{\frac{t_1 + k}{2} - 1} d u
\nonumber\\
& = \frac{1}{2^{\frac{k}{2}} \Gamma\left( \frac{k}{2} \right)} \frac{C_1'}{k^{\frac{t_1}{2}}} \left[ \int_{U( h )}^{2 U( h )} e^{- \frac{u}{2}} u^{\frac{t_1 + k}{2} - 1} d u \right] \times \nonumber\\
& \qquad
h^{t_1} \exp\left[ - \left( \frac{\left\| \mathbf{t} \right\|^2}{k} + C_2 k^{\frac{t_2}{2}} \left( \frac{2}{t_5} \right)^{t_3} \right) \left( m_1( h ) \right)^{\frac{1}{t_5}} \right]
\nonumber\\
& \ge \frac{1}{2^{\frac{k}{2}} \Gamma\left( \frac{k}{2} \right)} \frac{C_1'}{k^{\frac{t_1}{2}}} \left[ \int_{U( h )}^{2 U( h )} e^{- \frac{U( h )}{2}} \left( U( h ) \right)^{\frac{t_1 + k}{2} - 1} d u \right] \times \nonumber\\
& \qquad
h^{t_1} \exp\left[ - \left( \frac{\left\| \mathbf{t} \right\|^2}{k} + C_2 k^{\frac{t_2}{2}} \left( \frac{2}{t_5} \right)^{t_3} \right) \left( m_1( h ) \right)^{\frac{1}{t_5}} \right]
\nonumber\\
& = \frac{1}{2^{\frac{k}{2}} \Gamma\left( \frac{k}{2} \right)} \frac{C_1'}{k^{\frac{t_1}{2}}} \left( U( h ) \right)^{\frac{t_1 + k}{2}} \times \nonumber\\
& \qquad
h^{t_1} \exp\left[ - \left( \frac{1}{2} + \frac{\left\| \mathbf{t} \right\|^2}{k} + C_2 k^{\frac{t_2}{2}} \left( \frac{2}{t_5} \right)^{t_3} \right) \left( m_1( h ) \right)^{\frac{1}{t_5}} \right]
\nonumber\\
& > h^{t_1} \exp\left[ - \left( \frac{1}{2} + \frac{\left\| \mathbf{t} \right\|^2}{k} + C_2 k^{\frac{t_2}{2}} \left( \frac{2}{t_5} \right)^{t_3} \right) \left( \frac{1}{h} \right)^{\frac{t_2}{t_5}} \left( \log \frac{1}{h} \right)^{\frac{t_3}{t_5}} \right] .
\label{eq:prop2:8}
\end{align}
So, from \eqref{eq:prop2:6} and \eqref{eq:prop2:8}, we have for all sufficiently small $ h $,
\begin{align*}
& h^{t_1} \exp\left[ - u_1 \left( \frac{1}{h} \right)^{\frac{t_2}{t_5}} \left( \log \frac{1}{h} \right)^{\frac{t_3}{t_5}} \right]
< \phi_{\mathbf{T}}( \mathbf{t}, h )
< u_2 \exp\left[ - u_3 \left( \frac{1}{h} \right)^{\frac{t_2}{t_5}} \left( \log \frac{1}{h} \right)^{\frac{t_3}{t_5}} \right] ,
\\
\intertext{where}
& u_1 = \left( \frac{1}{2} + \frac{\left\| \mathbf{t} \right\|^2}{k} + C_2 k^{\frac{t_2}{2}} \left( \frac{2}{t_5} \right)^{t_3} \right) , \quad
u_2 = \left( \frac{\Gamma\left( \frac{t_4 + k}{2} \right)}{\Gamma\left( \frac{k}{2} \right)} \left( \frac{2}{k} \right)^{\frac{t_4}{2}} C_3 + 2^{\frac{k}{2}} \right) 
\\
& \text{and }
u_3 = \min\left\{ C_4 k^{\frac{t_2}{2}} \left( \frac{1}{t_2 + 2} \right)^{t_3} ,\, \frac{1}{4} \right\} .
\end{align*}
\end{proof}

\subsection{Proofs of theorems}

\begin{proof}[Proof of \autoref{ex:prop1}]
From the definitions of  $ \mathbb{L}_\mathbf{x}( \cdot ) $, $ G( \mathbf{Y} ) $ and $ F( \mathbf{z} ) $ in the statement of \autoref{ex:prop1}, and from \eqref{eq:bias} and \eqref{eq:variance}, we have
\begin{align*}
& B_n( \mathbf{x} ) \\
& = \Gamma'\left( \mathbb{E}[ \Psi( \mathbf{Y} ) \midil \mathbf{X} = \mathbf{x} ] \right)
\left( \frac{\sum_{i=1}^{n} \mathbb{E}[ \Psi( \mathbf{Y}_i ) \midil \mathbf{X}_i ] K( h_n^{-1} d( \mathbf{x}, \mathbf{X}_i ) ) }{\sum_{i=1}^{n} K( h_n^{-1} d( \mathbf{x}, \mathbf{X}_i ) )} - \mathbb{E}[ \Psi( \mathbf{Y} ) \midil \mathbf{X} = \mathbf{x} ] \right) , \\
& V_n( \mathbf{x} )
= \Gamma'\left( \mathbb{E}[ \Psi( \mathbf{Y} ) \midil \mathbf{X} = \mathbf{x} ] \right)
\left( \frac{\sum_{i=1}^{n} [ \Psi( \mathbf{Y}_i ) - \mathbb{E}[ \Psi( \mathbf{Y}_i ) \midil \mathbf{X}_i ] ] K( h_n^{-1} d( \mathbf{x}, \mathbf{X}_i ) ) }{\sum_{i=1}^{n} K( h_n^{-1} d( \mathbf{x}, \mathbf{X}_i ) )} \right) .
\end{align*}
Set $ R_n( \mathbf{x} ) = [ \widehat{\Theta}_n( \mathbf{x} ) - \Theta( \mathbf{x} ) ] - B_n( \mathbf{x} ) - V_n( \mathbf{x} ) $.
From \ref{assume:a2} and \eqref{eq:b2}, we have $ \| B_n( \mathbf{x} ) \| \longrightarrow 0 $ as $ n \longrightarrow \infty $. From \ref{assume:a2} and \autoref{thm:1}, we have $ \mathbb{E}[ \| V_n( \mathbf{x} ) \|^2 ] \longrightarrow 0 $ as $ n \longrightarrow \infty $, and consequently $ \| V_n( \mathbf{x} ) \| \longrightarrow 0 $ \textit{in probability} as $ n \longrightarrow \infty $. So, $ \| B_n( \mathbf{x} ) + V_n( \mathbf{x} ) \| \longrightarrow 0 $ \textit{in probability} as $ n \longrightarrow \infty $.
Therefore,
\begin{align}
& \| R_n( \mathbf{x} ) \| \nonumber\\
& = \left\| [ \widehat{\Theta}_n( \mathbf{x} ) - \Theta( \mathbf{x} ) ] - B_n( \mathbf{x} ) - V_n( \mathbf{x} ) \right\| \nonumber\\
& = \left\| \Gamma\left( \frac{\sum_{i=1}^{n} \Psi( \mathbf{Y}_i ) K( h_n^{-1} d( \mathbf{x}, \mathbf{X}_i ) ) }{\sum_{i=1}^{n} K( h_n^{-1} d( \mathbf{x}, \mathbf{X}_i ) )} \right)
- \Gamma( \mathbb{E}[ \Psi( \mathbf{Y} ) \midil \mathbf{X} = \mathbf{x} ] ) \right. \nonumber\\
& \left. \quad -
\Gamma'\left( \mathbb{E}[ \Psi( \mathbf{Y} ) \midil \mathbf{X} = \mathbf{x} ] \right)\left( \frac{\sum_{i=1}^{n} \Psi( \mathbf{Y}_i ) K( h_n^{-1} d( \mathbf{x}, \mathbf{X}_i ) ) }{\sum_{i=1}^{n} K( h_n^{-1} d( \mathbf{x}, \mathbf{X}_i ) )} - \mathbb{E}[ \Psi( \mathbf{Y} ) \midil \mathbf{X} = \mathbf{x} ] \right) \right\| \nonumber\\
& = o\left( \left\| B_n( \mathbf{x} ) + V_n( \mathbf{x} ) \right\| \right)
\quad \text{whenever } \left\| B_n( \mathbf{x} ) + V_n( \mathbf{x} )  \right\| \longrightarrow 0 ,
\nonumber\\
& = o\left( \| B_n( \mathbf{x} ) \| + \| V_n( \mathbf{x} ) \| \right)
\quad \text{whenever } \| B_n( \mathbf{x} ) \| + \| V_n( \mathbf{x} ) \| \longrightarrow 0 ,
\label{eq:e1}\\
& = o_\mathbb{P}\big( \max\big\{ h_n^{\beta}, [ n \phi( \mathbf{x}, h_n ) ]^{-1/2} \big\} \big)
\quad \text{as } n \longrightarrow 0 , \nonumber
\end{align}
since from \eqref{eq:b2}, we have $ \| B_n( \mathbf{x} ) \| = O( h_n^{\beta} ) $ as $ n \longrightarrow \infty $, and from \autoref{thm:1}, we have $ \| V_n( \mathbf{x} ) \| = O_\mathbb{P}( [ n \phi( \mathbf{x}, h_n ) ]^{-1/2} ) $ as $ n \longrightarrow \infty $. So, \ref{cond:4} is satisfied.
\end{proof}

\begin{proof}[Proof of \autoref{ex:prop2}]
Under the assumptions stated in \autoref{ex:2.3}, it follows that condition \ref{cond:1} holds from the Holder continuity of $ \Theta( \mathbf{z} ) $. The continuity of the linear operator $ \mathbb{L}_{\mathbf{x}}( \cdot ) $ follows from the invertibility of $ \mathbf{I}( \Theta( \mathbf{x} ) ) $, and \ref{cond:3} follows from the assumptions stated in \autoref{ex:2.3} using arguments similar to those used in the proof of Theorem 3.2 in \cite{chaudhuri1995likelihood}. We now proceed to verify condition \ref{cond:4}.

Using arguments similar to those in the proof of Theorem 3.1 in \cite{chaudhuri1995likelihood}, we get $ \widehat{\Theta}_n( \mathbf{x} ) \longrightarrow \Theta( \mathbf{x} ) $ \textit{in probability} as $ n \longrightarrow \infty $.
Using this fact, \eqref{ex2.3eq1a} and a Taylor expansion of $ \nabla g( \mathbf{Y}_i \midil \mathbf{t} ) $ at $ \mathbf{t} = \widehat{\Theta}_n( \mathbf{x} ) $, we get
\begin{align*}
& \sum_{i=1}^{n} \nabla g( \mathbf{Y}_i \midil \Theta( \mathbf{X}_i ) ) W_{i,n}( \mathbf{x} )
= \sum_{i=1}^{n} \Delta_2( g( \mathbf{Y}_i \midil \eta_i( \mathbf{x} ) ) ) \left( \Theta( \mathbf{X}_i ) - \widehat{\Theta}_n( \mathbf{x} ) \right) W_{i,n}( \mathbf{x} ) \\
& \implies \widehat{\Theta}_n( \mathbf{x} ) - \Theta( \mathbf{x} ) \\
& = \left[ \sum_{i=1}^{n} \Delta_2( g( \mathbf{Y}_i \midil \eta_i( \mathbf{x} ) ) ) W_{i,n}( \mathbf{x} ) \right]^{-1} \left( \sum_{i=1}^{n} \Delta_2( g( \mathbf{Y}_i \midil \eta_i( \mathbf{x} ) ) ) ( \Theta( \mathbf{X}_i ) - \Theta( \mathbf{x} ) ) W_{i,n}( \mathbf{x} ) \right) \\
& \quad
- \left[ \sum_{i=1}^{n} \Delta_2( g( \mathbf{Y}_i \midil \eta_i( \mathbf{x} ) ) ) W_{i,n}( \mathbf{x} ) \right]^{-1} \left( \sum_{i=1}^{n} \nabla g( \mathbf{Y}_i \midil \Theta( \mathbf{X}_i ) ) W_{i,n}( \mathbf{x} ) \right) ,
\end{align*}
where $ \eta_i( \mathbf{x} ) $ lies between $ \Theta( \mathbf{X}_i ) $ and $ \widehat{\Theta}_n( \mathbf{x} ) $.
Also, under the assumptions in \autoref{ex:2.3}, using arguments similar to those used in the proofs of Theorem 3.1 and Theorem 3.2 in \cite{chaudhuri1995likelihood}, we get that $ \| \sum_{i=1}^{n} \Delta_2( g( \mathbf{Y}_i \midil \eta_i( \mathbf{x} ) ) ) W_{i,n}( \mathbf{x} ) + \mathbf{I}( \Theta( \mathbf{x} ) ) \| \longrightarrow 0 $ \textit{in probability} as $ n \longrightarrow \infty $.
Also, since $ \Theta( \mathbf{z} ) \in \mathcal{F}( \mathbf{x}, \beta, \mathbb{R}^q ) $, we have $ \max\{ \| \Theta( \mathbf{X}_i ) - \Theta( \mathbf{x} ) \| W_{i,n}( \mathbf{x} ) \midil i = 1, \ldots, n \} \le c h_n^\beta $ for all $ n $, where $ c > 0 $ is a constant.
Consequently, it follows that
\begin{align*}
& \widehat{\Theta}_n( \mathbf{x} ) - \Theta( \mathbf{x} ) \\
& = [ \mathbf{I}( \Theta( \mathbf{x} ) ) ]^{-1} \left( \sum_{i=1}^{n} \mathbf{I}( \Theta( \mathbf{x} ) )( \Theta( \mathbf{X}_i ) - \Theta( \mathbf{x} ) ) W_{i,n}( \mathbf{x} ) \right) \\
& \quad
+ [ \mathbf{I}( \Theta( \mathbf{x} ) ) ]^{-1} \left( \sum_{i=1}^{n} \nabla g( \mathbf{Y}_i \midil \Theta( \mathbf{X}_i ) ) W_{i,n}( \mathbf{x} ) \right) \\
& \quad	+
o_P\left( h_n^\beta \right)
+ o_P\left( \left\| [ \mathbf{I}( \Theta( \mathbf{x} ) ) ]^{-1} \left( \sum_{i=1}^{n} \nabla g( \mathbf{Y}_i \midil \Theta( \mathbf{X}_i ) ) W_{i,n}( \mathbf{x} ) \right) \right\| \right) .
\end{align*}
Taking
\begin{align*}
V_n( \mathbf{x} ) &= [ \mathbf{I}( \Theta( \mathbf{x} ) ) ]^{-1} \left( \sum_{i=1}^{n} \nabla g( \mathbf{Y}_i \midil \Theta( \mathbf{X}_i ) ) W_{i,n}( \mathbf{x} ) \right) \\
\text{ and }
B_n( \mathbf{x} ) &= [ \mathbf{I}( \Theta( \mathbf{x} ) ) ]^{-1} \left( \sum_{i=1}^{n} \mathbf{I}( \Theta( \mathbf{x} ) )( \Theta( \mathbf{X}_i ) - \Theta( \mathbf{x} ) ) W_{i,n}( \mathbf{x} ) \right) ,
\end{align*}
we have
$ R_n( \mathbf{x} ) = o_\mathbb{P}\left( h_n^{\beta} + \| V_n( \mathbf{x} ) \| \right) $ as $ n \longrightarrow \infty $,
and the proof is complete using the convergence rate of $ \mathbb{E}[ \| V_n( \mathbf{x} ) \|^2 ] $ as described in the proof on \autoref{ex:prop1}.
\end{proof}

\begin{proof}[Proof of \autoref{thm:1}]
The arguments used in this proof are closely related to the arguments in the proof of Proposition 1 in \cite{chagny2016adaptive}.
Define
\begin{align*}
W_n( \mathbf{x} ) = n^{-1} \sum_{i=1}^{n} [ E_n^{(1)}( \mathbf{x} ) \phi( \mathbf{x}, h_n ) ]^{-1} K( h_n^{-1} d( \mathbf{x}, \mathbf{X}_i ) ) ,
\end{align*}
where $ E_n^{(1)}( \mathbf{x} ) $ is as defined in \eqref{eq:assume:a1}. It follows from Bernstein's inequality \citep[p.~95]{serfling2009approximation} and condition \ref{assume:a1} that
\begin{align}
\mathbb{P}[ | W_n( \mathbf{x} ) - 1 | > (1/2) ] \le 2 \exp( - c_1 n \phi( \mathbf{x}, h_n ) ) ,
\label{eq:2}
\end{align}
where $ c_1 $ is a positive constant. Note that
\begin{align}
& \mathbb{E}[ \| V_n( \mathbf{x} ) \|^2 ] \nonumber\\
&= \mathbb{E}[ \| V_n( \mathbf{x} ) \|^2 \mathbb{I}( W_n( \mathbf{x} ) < (1/2) ) ]
+ \mathbb{E}[ \| V_n( \mathbf{x} ) \|^2 \mathbb{I}( W_n( \mathbf{x} ) \ge (1/2) ) ] .
\label{eq:3}
\end{align}
For the first term on the RHS in \eqref{eq:3}, using the fact that $ \mathcal{B} $ is a type 2 Banach space and conditions \ref{assume:a1}, \ref{assume:a2} and \ref{cond:3}, we have from \eqref{eq:2},
\begin{align}
& \mathbb{E}[ \| V_n( \mathbf{x} ) \|^2 \mathbb{I}( W_n( \mathbf{x} ) < (1/2) ) ] \nonumber \\
& = \mathbb{E}[ \mathbb{E}[ \| V_n( \mathbf{x} ) \|^2 \mathbb{I}( W_n( \mathbf{x} ) < (1/2) ) \midil \mathbf{X}_1, \ldots, \mathbf{X}_n ] ] \nonumber \\
& \le c_2 \mathbb{E}\left[ \frac{\sum_{i=1}^{n} \mathbb{E}[ \| G( \mathbf{Y}_i ) - \mathbb{E}[ G( \mathbf{Y}_i ) \midil \mathbf{X}_i ] \|^2 \midil \mathbf{X}_i ] K^2( h_n^{-1} d( \mathbf{x}, \mathbf{X}_i ) )}{( \sum_{i=1}^{n} K( h_n^{-1} d( \mathbf{x}, \mathbf{X}_i ) ) )^2} \mathbb{I}\left( W_n( \mathbf{x} ) < \frac{1}{2} \right) \right] \nonumber \\
& \le c_3 \mathbb{E}\left[ \frac{\sum_{i=1}^{n} K^2( h_n^{-1} d( \mathbf{x}, \mathbf{X}_i ) )}{( \sum_{i=1}^{n} K( h_n^{-1} d( \mathbf{x}, \mathbf{X}_i ) ) )^2} \mathbb{I}\left( W_n( \mathbf{x} ) < \frac{1}{2} \right) \right] \nonumber \\
& \le c_3 \mathbb{P}[ | W_n( \mathbf{x} ) - 1 | > (1/2) ] 
\le 2 c_3 \exp( - c_1 n \phi( \mathbf{x}, h_n ) )
\label{eq:4}
\end{align}
for all sufficiently large $ n $, where $ c_2 $ and $ c_3 $ are positive constants. Since $ u e^{-u} \le e^{-1} $ for $ u > 0 $, from \eqref{eq:4}, we get that for all sufficiently large $ n $,
\begin{align}
& n \phi( \mathbf{x}, h_n ) \mathbb{E}[ \| V_n( \mathbf{x} ) \|^2 \mathbb{I}( W_n( \mathbf{x} ) < (1/2) ) ]
\le \frac{2 c_3}{c_1 e} .
\label{eq:5a}
\end{align}
Now, for the second term on the RHS in \eqref{eq:3}, again using the fact that $ \mathcal{B} $ is a type 2 Banach space, conditions \ref{assume:a1}, \ref{assume:a2}, \ref{cond:3} and inequality \eqref{eq:assume:a2}, we get that for all sufficiently large $ n $,
\begin{align}
& \mathbb{E}[ \| V_n( \mathbf{x} ) \|^2 \mathbb{I}( W_n( \mathbf{x} ) \ge (1/2) ) ] \nonumber \\
& \le \| \mathbb{L}_\mathbf{x} \|^2 \mathbb{E}\left[ \left\| \frac{1}{n} \sum_{i=1}^{n} \frac{[ G( \mathbf{Y}_i ) - \mathbb{E}[ G( \mathbf{Y}_i ) \midil \mathbf{X}_i ] ] K\left( \frac{d( \mathbf{x}, \mathbf{X}_i )}{h_n} \right)}{E_n^{(1)}( \mathbf{x} ) \phi( \mathbf{x}, h_n )} \right\|^2 \frac{\mathbb{I}\left( W_n( \mathbf{x} ) \ge \frac{1}{2} \right)}{\left( W_n( \mathbf{x} ) \right)^2} \right] \nonumber \\
& = \| \mathbb{L}_\mathbf{x} \|^2 \times \nonumber\\
& \quad \mathbb{E}\left[ \mathbb{E}\left[ \left\| \frac{1}{n} \sum_{i=1}^{n} \frac{[ G( \mathbf{Y}_i ) - \mathbb{E}[ G( \mathbf{Y}_i ) \midil \mathbf{X}_i ] ] K\left( \frac{d( \mathbf{x}, \mathbf{X}_i )}{h_n} \right)}{E_n^{(1)}( \mathbf{x} ) \phi( \mathbf{x}, h_n )} \right\|^2 \middle\arrowvert \mathbf{X}_1, \ldots, \mathbf{X}_n \right] \frac{\mathbb{I}\left( W_n( \mathbf{x} ) \ge \frac{1}{2} \right)}{\left( W_n( \mathbf{x} ) \right)^2} \right] \nonumber \\
& \le \| \mathbb{L}_\mathbf{x} \|^2 c_4 \mathbb{E}\left[ \sum_{i=1}^{n} \frac{\mathbb{E}\left[ \left\| G( \mathbf{Y}_i ) - \mathbb{E}[ G( \mathbf{Y}_i ) \midil \mathbf{X}_i ] \right\|^2 \middle\arrowvert \mathbf{X}_i \right] K^2\left( \frac{d( \mathbf{x}, \mathbf{X}_i )}{h_n} \right) \mathbb{I}\left( W_n( \mathbf{x} ) \ge \frac{1}{2} \right)}{\left( W_n( \mathbf{x} ) \right)^2 ( E_n^{(1)}( \mathbf{x} ) )^2 n^2 ( \phi( \mathbf{x}, h_n ) )^2} \right] \nonumber \\
& \le c_5 \mathbb{E}\left[ \sum_{i=1}^{n}
\frac{K^2( h_n^{-1} d( \mathbf{x}, \mathbf{X}_i ) ) \mathbb{I}( W_n( \mathbf{x} ) \ge (1/2) )}{\left( W_n( \mathbf{x} ) \right)^2 ( E_n^{(1)}( \mathbf{x} ) )^2 n^2 ( \phi( \mathbf{x}, h_n ) )^2} \right] \nonumber \\
& \le 4 c_5 \mathbb{E}\left[ \sum_{i=1}^{n}
\frac{K^2( h_n^{-1} d( \mathbf{x}, \mathbf{X}_i ) )}{( E_n^{(1)}( \mathbf{x} ) )^2 n^2 ( \phi( \mathbf{x}, h_n ) )^2} \right] \nonumber \\
& = \frac{4 c_5 E_n^{(2)}( \mathbf{x} )}{( E_n^{(1)}( \mathbf{x} ) )^2} \frac{1}{n \phi( \mathbf{x}, h_n )}
\le \frac{4 c_5 L^2}{l^2} \frac{1}{n \phi( \mathbf{x}, h_n )} \nonumber\\
& \implies
\phi( \mathbf{x}, h_n ) \mathbb{E}[ \| V_n( \mathbf{x} ) \|^2 \mathbb{I}( W_n( \mathbf{x} ) \ge (1/2) ) ]
\le \frac{4 c_5 L^2}{l^2} ,
\label{eq:5b}
\end{align}
where $ c_4 $ and $ c_5 $ are positive constants. From \eqref{eq:3}, \eqref{eq:5a} and \eqref{eq:5b}, we get $ n \phi( \mathbf{x}, h_n ) \allowbreak \mathbb{E}[ \| V_n( \mathbf{x} ) \|^2 ] = O( 1 ) $ as $ n \longrightarrow \infty $.
\end{proof}

\begin{proof}[Proof of \autoref{thm:2}]
Note that
\begin{align}
& \left[ n \phi( \mathbf{x}, h_n ) \right]^{1/2} [ E_n^{(2)}( \mathbf{x} ) ]^{-1/2} E_n^{(1)}( \mathbf{x} ) V_n( \mathbf{x} ) \nonumber\\
& = \frac{\sum_{i=1}^{n} \frac{K( h_n^{-1} d( \mathbf{x}, \mathbf{X}_i ) )}{[ E_n^{(2)}( \mathbf{x} ) ]^{1/2} \left[ n \phi( \mathbf{x}, h_n ) \right]^{1/2}} \mathbb{L}_\mathbf{x}( G( \mathbf{Y}_i ) - \mathbb{E}[ G( \mathbf{Y}_i ) \midil \mathbf{X}_i ] )}{n^{-1} \sum_{i=1}^{n} \left[ E_n^{(1)}( \mathbf{x} ) \phi( \mathbf{x}, h_n ) \right]^{-1} K( h_n^{-1} d( \mathbf{x}, \mathbf{X}_i ) )} .
\label{thm:2a}
\end{align}
Define
\begin{align*}
V^*_n( \mathbf{x} ) = \sum_{i=1}^{n} \frac{K( h_n^{-1} d( \mathbf{x}, \mathbf{X}_i ) )}{[ E_n^{(2)}( \mathbf{x} ) ]^{1/2} \left[ n \phi( \mathbf{x}, h_n ) \right]^{1/2}} \mathbb{L}_\mathbf{x}( G( \mathbf{Y}_i ) - \mathbb{E}[ G( \mathbf{Y}_i ) \midil \mathbf{X}_i ] ) .
\end{align*}
The covariance operator of $ V^*_n( \mathbf{x} ) $, denoted as $ \mathbb{D}_n( \cdot, \cdot \midil \mathbf{x} ) $, is given by
\begin{align*}
& \mathbb{D}_n( \mathbf{u}, \mathbf{v} \midil \mathbf{x} ) \\
& = \mathbb{E}\left[ 
\frac{\left\langle \mathbf{u}, \mathbb{L}_\mathbf{x}( G( \mathbf{Y} ) - \mathbb{E}[ G( \mathbf{Y} ) \midil \mathbf{X} ] ) \right\rangle
\left\langle \mathbf{v}, \mathbb{L}_\mathbf{x}( G( \mathbf{Y} ) - \mathbb{E}[ G( \mathbf{Y} ) \midil \mathbf{X} ] ) \right\rangle K^2\left( \frac{d( \mathbf{x}, \mathbf{X} )}{h_n} \right)}{E_n^{(2)}( \mathbf{x} ) \phi( \mathbf{x}, h_n )} \right]
\end{align*}
for $ \mathbf{u}, \mathbf{v} \in \mathcal{B} $.
Under conditions \ref{assume:a1}, \ref{assume:a2} and \ref{cond:5}, $ \mathbb{D}_n( \cdot, \cdot \midil \mathbf{x} ) $ converges to $ \mathbb{D}( \cdot, \cdot \midil \mathbf{x} ) $ in the trace norm as $ n \longrightarrow \infty $. Consequently, conditions (i) and (ii) in Theorem 1.1 in \cite{kundu2000central} are satisfied.
Define
\begin{align*}
\mathbf{U}_{n,i}( \mathbf{x} )
= \frac{K( h_n^{-1} d( \mathbf{x}, \mathbf{X}_i ) )}{[ E_n^{(2)}( \mathbf{x} ) ]^{1/2}} 
\mathbb{L}_\mathbf{x}( G( \mathbf{Y}_i ) - \mathbb{E}[ G( \mathbf{Y}_i ) \midil \mathbf{X}_i ] ) .
\end{align*}
Given $ \epsilon > 0 $ and $ \mathbf{b} \in \mathcal{B} $, define
\begin{align*}
L_n( \epsilon, \mathbf{b} )
= \sum_{i=1}^{n} \mathbb{E}\left[
\left\langle \frac{\mathbf{U}_{n,i}( \mathbf{x} )}{\left[ n \phi( \mathbf{x}, h_n ) \right]^{1/2}} ,
\mathbf{b} \right\rangle^2
\mathbb{I}\left[
\left| \left\langle \frac{\mathbf{U}_{n,i}( \mathbf{x} )}{\left[ n \phi( \mathbf{x}, h_n ) \right]^{1/2}} ,
\mathbf{b} \right\rangle \right|
> \epsilon \right] \right] .
\end{align*}
From \ref{assume:a1}, \ref{assume:a2} and \ref{cond:3}, we have for any $ \mathbf{b} $ with $ \| \mathbf{b} \| = 1 $,
\begin{align*}
& L_n( \epsilon, \mathbf{b} ) \\
& = \mathbb{E}\left[
\frac{\left\langle \mathbf{U}_{n,1}( \mathbf{x} ) , \mathbf{b} \right\rangle^2}{\phi( \mathbf{x}, h_n )}
\mathbb{I}\left[
\left| \left\langle \mathbf{U}_{n,1}( \mathbf{x} ) , \mathbf{b} \right\rangle \right|
> \epsilon \left[ n \phi( \mathbf{x}, h_n ) \right]^{1/2} \right] \right] \\
& \le \mathbb{E}\left[
\frac{\left\langle \mathbf{U}_{n,i}( \mathbf{x} ) , \mathbf{b} \right\rangle^2}{\phi( \mathbf{x}, h_n )}
\left[ \frac{\left| \left\langle \mathbf{U}_{n,1}( \mathbf{x} ) , \mathbf{b} \right\rangle \right|}{\epsilon \left[ n \phi( \mathbf{x}, h_n ) \right]^{1/2}} \right]^{\nu - 2}
\right] \\
& \le 
\left[ \frac{1}{\epsilon \left[ n \phi( \mathbf{x}, h_n ) \right]^{1/2}} \right]^{\nu - 2} \mathbb{E}\left[
\frac{ \| \mathbb{L}_\mathbf{x} \|^\nu \left\|
G( \mathbf{Y} ) - \mathbb{E}[ G( \mathbf{Y} ) \midil \mathbf{X} ] \right\|^\nu K^\nu( h_n^{-1} d( \mathbf{x}, \mathbf{X} ) )}{[ E_n^{(2)}( \mathbf{x} ) ]^{\nu/2} \phi( \mathbf{x}, h_n )}
\right] 
\\
& \le
c \left[ n \phi( \mathbf{x}, h_n ) \right]^{-\frac{\nu-2}{2}}
\longrightarrow
0
\end{align*}
as $ n \longrightarrow \infty $, where $ \nu > 2 $ is the constant mentioned in \ref{cond:3}. Hence, condition (iii) in Theorem 1.1 in \cite{kundu2000central} is satisfied.
Consequently,
\begin{align}
V^*_n( \mathbf{x} ) \longrightarrow \mathbf{W}
\label{thm:2b}
\end{align}
\textit{in distribution} as $ n \longrightarrow \infty $. 
Now, under conditions \ref{assume:a1}, \ref{assume:a2} and an application of the Markov inequality, we get
\begin{align}
n^{-1} \sum_{i=1}^{n} \left[ E_n^{(1)}( \mathbf{x} ) \phi( \mathbf{x}, h_n ) \right]^{-1} K( h_n^{-1} d( \mathbf{x}, \mathbf{X}_i ) )
\longrightarrow 1
\label{thm:2c}
\end{align}
\textit{in probability} as $ n \longrightarrow \infty $.
The proof is completed from \eqref{thm:2a}, \eqref{thm:2b}, \eqref{thm:2c} and an application of Slutsky's Theorem.
\end{proof}

\begin{proof}[Proof of \autoref{thm:up}]
From the upper bounds of $ \mathbb{E}\| B_n( \mathbf{x} ) \|^2 $ and $ \mathbb{E}\| V_n( \mathbf{x} ) \|^2 $ in \eqref{eq:b2} and \autoref{thm:1}, respectively, and the lower bound of $ \phi( \mathbf{x}, h_n ) $ in \eqref{eq:1}, we have
\begin{align*}
\mathbb{E}\| B_n( \mathbf{x} ) + V_n( \mathbf{x} ) \|^2
\le 2 \left[ \mathbb{E}\| B_n( \mathbf{x} ) \|^2 + \mathbb{E}\| V_n( \mathbf{x} ) \|^2 \right]
\le f_1( h_n )
\end{align*}
for all sufficiently large $ n $, where
\begin{align}
f_1( h_n ) = a h_n^{2 \beta} + \frac{b}{n C_1} ( 1 / h_n )^{t_1} \exp\left[ m( h_n ) \right] ,
\label{eq:thm:up}
\end{align}
and $ a, b > 0 $ are some constants.
We establish below that the choice of bandwidths $ \{ h_n \} $ described in the statement of \autoref{thm:up} is one which minimizes \eqref{eq:thm:up}.
Note that $ m( h ) $, which is defined in \eqref{eq:mh}, is a differentiable function of $ h $, and
\begin{align}
m'( h ) = - m( h ) ( 1 / h ) \left( t_2 + \frac{t_3}{\log ( 1 / h )} \right) .
\label{eq:thm:upa}
\end{align}
Consequently, $ f_1( h ) $ is differentiable for all $ n $, and
\begin{align}
f_1'( h )
& = 2 \beta a h^{2 \beta - 1} 
- \frac{b t_1}{n C_1} ( 1 / h )^{t_1 + 1} \exp[ m( h ) ] \nonumber\\
& \quad - \frac{b}{n C_1} ( 1 / h )^{t_1 + 1} \exp[ m( h ) ] m( h ) \left( t_2 + \frac{t_3}{\log ( 1 / h )} \right)
\label{eq:thm:upb} \\
& = \exp[ m( h ) ] \left[ \frac{2 \beta a h^{2 \beta - 1}}{\exp[ m( h ) ]}
- \frac{b t_1}{n C_1} ( 1 / h )^{t_1 + 1} \right. \nonumber\\
& \qquad\qquad\qquad\quad \left. - \frac{b}{n C_1} ( 1 / h )^{t_1 + 1} m( h ) \left( t_2 + \frac{t_3}{\log ( 1 / h )} \right) \right] .
\label{eq:thm:upb1}
\end{align}
From \eqref{eq:thm:upb1}, we get that for every fixed $ n $, $ f_1'( h ) \longrightarrow - \infty $ as $ h \longrightarrow 0^+ $, and for any $ 0 < s < 1 $, $ f_1'( s ) > 0 $ for all sufficiently large $ n $. Since $ f_1'( h ) $ is continuous in $ h $ for $ 0 < h < 1 $, given any $ 0 < s < 1 $, $ f_1'( h ) $ must have a root in $ ( 0 , s ) $ for all sufficiently large $ n $.
For any fixed $ n $, consider $ h_0 = \inf\{ h \midil f_1'( h ) = 0 \} $. Again, since $ f_1'( h ) $ is continuous in $ h $, we have $ f_1'( h_0 ) = 0 $. Further, since $ f_1'( h ) \longrightarrow - \infty $ as $ h \longrightarrow 0^+ $, from the continuity of $ f_1'( h ) $ we have $ f_1'( h ) < 0 $ for $ h < h_0 $, which implies that $ f_1( h ) $ is a decreasing function for $ h < h_0 $. Also, for any $ 0 < s < s' < 1 $, we have for all sufficiently large $ n $, $ f_1'( h ) > 0 $ for all $ s \le h \le s' $, which implies $ f_1( h ) $ is increasing in $ s \le h \le s' $. Therefore, $ f_1( h ) $ must have a minima for all sufficiently large $ n $, whose corresponding $ h $ will satisfy $ f_1'( h ) = 0 $.
Now, from \eqref{eq:thm:upb}, $ f_1'( h_n ) = 0 $ implies that
\begin{align}
& 2 \beta a h_n^{2 \beta - 1}
= n^{-1} ( 1 / h_n )^{t_1 + 1} \exp[ m( h_n ) ] \nonumber\\
& \qquad \qquad \qquad \quad \times \left[ \frac{b t_1}{C_1}
+ \frac{b}{C_1} m( h_n ) \left( t_2 + \frac{t_3}{\log ( 1 / h_n )} \right) \right]
\nonumber\\
& \iff
h_n^{2 \beta} 
= n^{-1} (1 / h_n)^{t_1} \exp[ m( h_n ) ] \nonumber\\
& \qquad \qquad \qquad \quad \times \left[ \frac{b t_1}{2 \beta a C_1}
+ \frac{b}{2 \beta a C_1} m( h_n ) \left( t_2 + \frac{t_3}{\log ( 1 / h_n )} \right) \right] 
\label{eq:proof1}\\
& \iff
n
= (1 / h_n)^{2 \beta + t_1} \exp[ m( h_n ) ] \nonumber\\
& \qquad \qquad \qquad \quad \times \left[ \frac{b t_1}{2 \beta a C_1}
+ \frac{b}{2 \beta a C_1} m( h_n ) \left( t_2 + \frac{t_3}{\log ( 1 / h_n )} \right) \right] 
\label{eq:proof2}\\
& \iff
\frac{\log n}{m( h_n )}
= 1 + ( 2 \beta + t_1 ) \frac{\log ( 1 / h_n )}{m( h_n )} + \frac{1}{m( h_n )} \nonumber\\
& \qquad \qquad \qquad
\times \log \left( \left[ \frac{b t_1}{2 \beta a C_1}
+ \frac{b}{2 \beta a C_1} m( h_n ) \left( t_2 + \frac{t_3}{\log ( 1 / h_n )} \right) \right] \right) .
\label{eq:thm:upd}
\end{align}
Let $ \{ h_n \} $ be such that $ f_1'( h_n ) = 0 $ for all $ n $.
If either $ t_2 > 0 $ or $ t_3 > 1 $, then from \eqref{eq:proof2}, we get that $ h_n \longrightarrow 0^+ $ as $ n \longrightarrow \infty $.
Consequently, from \eqref{eq:proof1}, we have
\begin{align*}
n C_1 h_n^{t_1} \exp[ - m( h_n ) ]
= 
h_n^{- 2 \beta} \left[ \frac{b t_1}{2 \beta a}
+ \frac{b}{2 \beta a} m( h_n ) \left( t_2 + \frac{t_3}{\log ( 1 / h_n )} \right) \right] 
\longrightarrow \infty
\end{align*}
as $ n \longrightarrow \infty $, which implies $ n \phi( \mathbf{x}, h_n ) \longrightarrow \infty $ as $ n \longrightarrow \infty $ from the lower bound of $ \phi( \mathbf{x}, h_n ) $ in \eqref{eq:1}. Therefore, $ \{ h_n \} $ satisfies \ref{assume:a2}.
Also, $ h_n \longrightarrow 0^+ $ as $ n \longrightarrow \infty $ implies that
\begin{align}
& \frac{\log ( 1 / h_n )}{m( h_n )} \longrightarrow 0 
\text{ as } n \longrightarrow \infty ,
\label{eq:thm:upe} \\
& \frac{1}{m( h_n )} \log \left( \left[ \frac{b t_1}{2 \beta a C_1}
+ \frac{b}{2 \beta a C_1} m( h_n ) \left( t_2 + \frac{t_3}{\log ( 1 / h_n )} \right) \right] \right) \longrightarrow 0
\label{eq:thm:upf}
\end{align}
as $ n \longrightarrow \infty $. Combining \eqref{eq:thm:upd}, \eqref{eq:thm:upe} and \eqref{eq:thm:upf}, we have
\begin{align}
\frac{\log n}{m( h_n )} \longrightarrow 1 \text{ as } n \longrightarrow \infty .
\label{eq:thm:upg}
\end{align}
Consequently,
when either $ t_2 > 0 $ or $ t_3 > 1 $, we have for all sufficiently large $ n $,
\begin{align}
a h_n^{2 \beta} < a C_2' ( m^{-1}( \log n ) )^{2 \beta} ,
\label{eq:thm:mse1}
\end{align}
where $ C_2' $ is a positive constant depending on $ C_2 $ and $ \beta $.
From \eqref{eq:thm:up} and \eqref{eq:proof1}, we get that
$ a h_n^{2 \beta} < f_1( h_n ) < 2 a h_n^{2 \beta} $
for all sufficiently large $ n $, and consequently 
$ f_1( h_n ) < 2 a C_2' ( m^{-1}( \log n ) )^{2 \beta} $ 
for all sufficiently large $ n $.
Hence, for the bandwidth sequence $ \{ h_n \} $ minimizing $ f_1( h ) $ for every fixed $ n $, we have
\begin{align}
\mathbb{E}\| B_n( \mathbf{x} ) + V_n( \mathbf{x} ) \|^2 < 2 a C_2' \left( m^{-1}\left( \log n \right) \right)^{2 \beta}
\label{eq:thm:mse2}
\end{align}
for all sufficiently large $ n $, which implies
$ \| B_n( \mathbf{x} ) + V_n( \mathbf{x} ) \| = O_\mathbb{P}\big( \left( m^{-1}\left( \log n \right) \right)^\beta \big) $
as $ n \longrightarrow \infty $.
Also, when either $ t_2 > 0 $ or $ t_3 > 1 $, from \eqref{eq:proof1} and the lower bound of $ \phi( \mathbf{x}, h ) $ in \eqref{eq:1}, we get
\begin{align}
h_n^{2 \beta} / [ [ n \phi( \mathbf{x}, h_n ) ]^{-1} ] \longrightarrow \infty
\label{eq:thm:mse3}
\end{align}
as $ n \longrightarrow \infty $. Hence, from \eqref{eq:thm:mse1} and \eqref{eq:thm:mse3}, we get that $ \left( m^{-1}\left( \log n \right) \right)^{- \beta} \| R_n( \mathbf{x} ) \| = o_\mathbb{P}( 1 ) $ as $ n \longrightarrow \infty $. Therefore, $ \big\| \widehat{\Theta}_n( \mathbf{x} ) - \Theta( \mathbf{x} ) \big\| = O_\mathbb{P}\big( \left( m^{-1}( \log n ) \right)^\beta \big) $ as $ n \longrightarrow \infty $.

Next, if $ \mathbb{E}[ \| R_n( \mathbf{x} ) \|^2 ] = o\left( \delta_n^2 \right) $ as $ n \longrightarrow \infty $, where $ \delta_n = \max\big\{ h_n^\beta, \allowbreak \big[ n \phi( \mathbf{x}, h_n ) \big]^{-1/2} \big\} $, we get from \eqref{eq:thm:mse1} and \eqref{eq:thm:mse3} that
$ \left( m^{-1}\left( \log n \right) \right)^{- 2 \beta} \mathbb{E}\| R_n( \mathbf{x} ) \|^2 \allowbreak \longrightarrow 0 $ as $ n \longrightarrow \infty $.
From \eqref{eq:main}, we have 
$ \mathbb{E}\big\| \widehat{\Theta}_n( \mathbf{x} ) - \Theta( \mathbf{x} ) \big\|^2
\le 2 \mathbb{E}\big\| B_n( \mathbf{x} ) + V_n( \mathbf{x} ) \big\|^2 + 2 \mathbb{E}\| R_n( \mathbf{x} ) \|^2 $.
Therefore, from \eqref{eq:thm:mse2}, we get
$ \mathbb{E}\big\| \widehat{\Theta}_n( \mathbf{x} ) - \Theta( \mathbf{x} ) \big\|^2 = O\big( \left( m^{-1}( \log n ) \right)^{2 \beta} \big) $ as $ n \longrightarrow \infty $.
\end{proof}

\begin{proof}[Proof of \autoref{thm:1a}]
For a sequence of bandwidths $ \{ h_n \} $ satisfying \ref{assume:a2}, it follows from the upper bound of $ \phi( \mathbf{x}, h ) $ in \eqref{eq:1} and the definition of $ m( h ) $ in \eqref{eq:mh} that
\begin{align}
& n C_3 h_n^{t_4} \exp\left[ - (C_4 / C_2) m( h_n ) \right]
\ge n \phi( \mathbf{x}, h_n ) \longrightarrow \infty \text{ as } n \longrightarrow \infty
\nonumber \\
& \implies
\log n - t_4 \log ( 1 / h_n ) - (C_4 / C_2) m( h_n ) \longrightarrow \infty
\nonumber \\
& \iff
m( h_n )
\left[ \frac{\log n}{m( h_n )} - t_4 \frac{\log ( 1 / h_n )}{m( h_n )} - \frac{C_4}{C_2} \right] \longrightarrow \infty
\label{eq:thm:1a1}
\end{align}
as $ n \longrightarrow \infty $. Now, since either $ t_2 > 0 $ or $ t_3 > 1 $, and $ h_n \longrightarrow 0 $ as $ n \longrightarrow \infty $ under assumption \ref{assume:a2}, we have
\begin{align}
m( h_n ) \longrightarrow \infty \text{ as } n \longrightarrow \infty
\quad
\text{and}
\quad
\frac{\log ( 1 / h_n )}{m( h_n )} \longrightarrow 0 \text{ as } n \longrightarrow \infty .
\label{eq:thm:1a2}
\end{align}
Hence, for \eqref{eq:thm:1a1} to be satisfied, in view of \eqref{eq:thm:1a2}, we must have, for all sufficiently large $ n $,
\begin{align}
\frac{\log n}{m( h_n )} - \frac{C_4}{C_2}
> 0 
& \iff
m^{-1}\left( \frac{\log n}{(C_4 / C_2)} \right) < h_n 
\implies
\frac{h_n}{m^{-1}\left( \log n \right)} > c_1 > 0 ,
\label{eq:thm:1a4}
\end{align}
where $ c_1 $ is a constant depending on $ C_2, C_4 $. Clearly, when $ C_2 = C_4 $, $ c_1 = 1 $.
\end{proof}

\begin{proof}[Proof of \autoref{thm:low}]
Suppose, if possible,
\begin{align}
\liminf_{n \longrightarrow \infty} \mathbb{P}\left[ \left( m^{-1}( \log n ) \right)^{-\beta} \allowbreak \left\| \widehat{\Theta}_n( \mathbf{x} ) - \Theta( \mathbf{x} ) \right\| > c \right] = 0
\label{eq:thm:low1}
\end{align}
for every $ c > 0 $. Then, given any $ c > 0 $, there is a subsequence $ \{ n' \} $ such that
\begin{align}
\lim_{n' \longrightarrow \infty} \mathbb{P}\left[ \left( m^{-1}( \log n' ) \right)^{-\beta} \left\| \widehat{\Theta}_n'( \mathbf{x} ) - \Theta( \mathbf{x} ) \right\| > c \right] = 0 .
\label{eq:thm:low2}
\end{align}
Consider the bandwidth sequence $ \{ h_{n'} \} $.
If $ \liminf_{n' \longrightarrow \infty} h_{n'}^{2 \beta} n \phi( \mathbf{x}, h_{n'} ) = 0 $, then there exists a further subsequence $ \{ n'' \} $ such that $ h_{n''}^{2 \beta} n \phi( \mathbf{x}, h_{n''} )
\longrightarrow 0 $ as $ n'' \longrightarrow \infty $. But in this case, we get a contradiction of \eqref{eq:thm:low2} from \autoref{lemma:thmlow1}.
On the other hand, if $ \limsup_{n' \longrightarrow \infty} h_{n'}^{2 \beta} n \phi( \mathbf{x}, h_{n'} ) \allowbreak = \infty $, then there exists a further subsequence $ \{ n'' \} $ such that $ h_{n''}^{2 \beta} n \phi( \mathbf{x}, h_{n''} )
\allowbreak \longrightarrow \infty $ as $ n'' \longrightarrow \infty $. But again, we get a contradiction of \eqref{eq:thm:low2} from \autoref{lemma:thmlow2}.
We consider the only remaining case, which is
$ 0
< \liminf_{n' \longrightarrow \infty} h_{n'}^{2 \beta} n \phi( \mathbf{x}, h_{n'} )
\le \limsup_{n' \longrightarrow \infty} h_{n'}^{2 \beta} n \phi( \mathbf{x}, h_{n'} ) 
< \infty $. Then,
there exist $ \epsilon_1 > 0 $, $ \epsilon_2 > 0 $ and a further subsequence $ \{ n'' \} $ such that $ 0 < \epsilon_1 <
h_{n''}^{2 \beta} n \phi( \mathbf{x}, h_{n''} )
< \epsilon_2 $ for all sufficiently large $ n'' $. But in this case also, we get a contradiction of \eqref{eq:thm:low2} from \autoref{lemma:thmlow3}.
Therefore, the assertion \eqref{eq:thm:low1} is not possible, and this completes the proof.
\end{proof}

\begin{proof}[Proof of \autoref{thm:3}]
From \eqref{eq:proof1} in the proof of \autoref{thm:up} and the lower bound of $ \phi( \mathbf{x}, h ) $ in \eqref{eq:1}, it follows that
\begin{align}
h_n^{2 \beta} n \phi( \mathbf{x}, h_n ) \longrightarrow \infty
\label{eq:thm:3a}
\end{align}
as $ n \longrightarrow \infty $. Now, choose $ \Theta( \cdot ) $ as in \autoref{thm:low} such that $ h_n^{-\beta} \| \tilde{B}_n( \mathbf{x} ) \| \ge b_1 > 0 $ for a constant $ b_1 $ and all sufficiently large $ n $. So, $ \mathbb{P}[ h_n^{- \beta} \| B_n( \mathbf{x} ) \| > b_1 / 2 ] \longrightarrow 1 $ as $ n \longrightarrow \infty $. Hence, for this choice of $ \Theta( \cdot ) $ and using \autoref{thm:2} and \eqref{eq:thm:3a}, we have
\begin{align*}
\frac{\| V_n( \mathbf{x} ) \|}{\| B_n( \mathbf{x} ) \|}
& = \frac{1}{[ n \phi( \mathbf{x}, h_n ) ]^{1/2} h_n^\beta} \frac{[ n \phi( \mathbf{x}, h_n ) ]^{1/2} \| V_n( \mathbf{x} ) \|}{h_n^{- \beta} \| B_n( \mathbf{x} ) \|}
= o_\mathbb{P}( 1 )
\text{ as }
n \longrightarrow \infty .
\end{align*}
\end{proof}

\begin{proof}[Proof of \autoref{coro:1}]
Let $ \mathbf{U} $ and $ \mathbf{V} $ be two nonnegative random variables. Then, given any $ \epsilon > 0 $ and $ \delta > 0 $, we have
\begin{align}
& \mathbb{P}\left[ \frac{\mathbf{U}}{\mathbf{V}} < \epsilon \right] 
\ge \mathbb{P}\left[ \mathbf{U} < \epsilon \delta ,\, \mathbf{V} > \delta \right] 
\ge \mathbb{P}\left[ \mathbf{U} < \epsilon \delta \right] + \mathbb{P}\left[ \mathbf{V} > \delta \right] - 1 .
\label{eq:coro1eq1}
\end{align}
We denote our optimum bandwidth minimizing \eqref{eq:thm:up} in the proof of \autoref{thm:up} as $ h_n^{(op)} $.
Given any $ \epsilon > 0 $, from \autoref{lemma:newlemma2}, we get that there is $ \delta > 0 $ such that
\begin{align}
\mathbb{P}\left[ ( h_n^{(b)} )^{-\beta} \left\| \widehat{\Theta}_n^{(b)}( \mathbf{x} ) - \Theta( \mathbf{x} ) \right\| > \delta \right]
> 1 - \epsilon 
\label{eq:coro1eq2}
\end{align}
for all sufficiently large $ n $. Further, from \autoref{lemma:newlemma1}, we get that for this constant $ \delta $,
\begin{align}
\mathbb{P}\left[ ( h_n^{(op)} )^{-\beta} \left\| \widehat{\Theta}_n^{(op)}( \mathbf{x} ) - \Theta( \mathbf{x} ) \right\| < \epsilon \delta \right]
> 1 - \epsilon
\label{eq:coro1eq3}
\end{align}
for all sufficiently large $ n $.
Therefore, from \eqref{eq:coro1eq1}, \eqref{eq:coro1eq2} and \eqref{eq:coro1eq3}, we get that
\begin{align}
\frac{( h_n^{(op)} )^{-\beta} \left\| \widehat{\Theta}_n^{(op)}( \mathbf{x} ) - \Theta( \mathbf{x} ) \right\|}{( h_n^{(b)} )^{-\beta} \left\| \widehat{\Theta}_n^{(b)}( \mathbf{x} ) - \Theta( \mathbf{x} ) \right\|}
= o_\mathbb{P}( 1 )
\quad
\text{as }
n \longrightarrow \infty .
\label{eq:coro1eq4}
\end{align}
Hence, from \eqref{eq:coro1eq4} and \autoref{lemma:coro6}, we have
\begin{align*}
\frac{\left\| \widehat{\Theta}_n^{(op)}( \mathbf{x} ) - \Theta( \mathbf{x} ) \right\|}{\left\| \widehat{\Theta}_n^{(b)}( \mathbf{x} ) - \Theta( \mathbf{x} ) \right\|}
= o_\mathbb{P}( 1 )
\quad
\text{as }
n \longrightarrow \infty .
\end{align*}
On the other hand, from \autoref{lemma:coro6}, \autoref{lemma:newlemma1} and \autoref{lemma:newlemma2}, we have
\begin{align*}
& \frac{\mathbb{E}\left\| \widehat{\Theta}_n^{(op)}( \mathbf{x} ) - \Theta( \mathbf{x} ) \right\|^2}{\mathbb{E}\left\| \widehat{\Theta}_n^{(b)}( \mathbf{x} ) - \Theta( \mathbf{x} ) \right\|^2}
= o( 1 )
\quad
\text{as }
n \longrightarrow \infty .
\end{align*}
\end{proof}

\begin{proof}[Proof of \autoref{thm:5}]
This proof is partly based on arguments used in \cite{chagny2014adaptive,chagny2016adaptive}.
For every $ h \in \mathbb{H}_n $, we have
\begin{align}
& \left\| \widehat{\Theta}_n\left( \mathbf{x}, h_n^* \right) - \Theta\left( \mathbf{x} \right) \right\|^2
\nonumber\\
& \le 
3 \left[ \left\| \widehat{\Theta}_n\left( \mathbf{x}, h_n^* \right) - \widehat{\Theta}_n\left( \mathbf{x}, \max\{ h_n^*, h \} \right) \right\|^2
+ \left\| \widehat{\Theta}_n\left( \mathbf{x}, h \right) - \widehat{\Theta}_n\left( \mathbf{x}, \max\{ h_n^*, h \} \right) \right\|^2 \right] \nonumber\\
& \quad + 3 \left\| \widehat{\Theta}_n( \mathbf{x}, h ) - \Theta( \mathbf{x} ) \right\|^2
\nonumber\\
& \le
3 \left[ \left( C_n( \mathbf{x}, h ) + D_n( \mathbf{x}, h_n^* ) \right) + \left( C_n( \mathbf{x}, h_n^* ) + D_n( \mathbf{x}, h ) \right) \right]
+ 3 \left\| \widehat{\Theta}_n( \mathbf{x}, h ) - \Theta( \mathbf{x} ) \right\|^2
\nonumber\\
& \le
6 \left[ C_n( \mathbf{x}, h ) + D_n( \mathbf{x}, h ) \right]
+ 3 \left\| \widehat{\Theta}_n( \mathbf{x}, h ) - \Theta( \mathbf{x} ) \right\|^2 .
\label{eq:adapt1}
\end{align}
From \autoref{lemma:adapt1}, \autoref{lemma:adapt2} and \autoref{lemma:adapt4}, we get
\begin{align}
& C_n( \mathbf{x}, h )
\le C_n^{(1)}( \mathbf{x}, h ) + C_n^{(2)}( \mathbf{x}, h ) .
\label{eq:adapt2}
\end{align}
Here, for all sufficiently large $ n $,
\begin{align}
& \mathbb{E}\left[ C_n^{(1)}( \mathbf{x}, h ) \right]
\le c_1 h^{2 \beta} + \frac{1}{n \log n}
\label{eq:adapt3}
\end{align}
for all $ h \in \mathbb{H}_n $ and some constant $ c_1 > 0 $ independent of $ h $. Also,
\begin{align}
& \mathbb{P}\left[ C_n^{(2)}( \mathbf{x}, h ) > n^{-2} \right] = O\left( n^{-2} \right)
\quad
\text{as }
n \longrightarrow \infty .
\label{eq:adapt4}
\end{align}
Further, $ C_n^{(2)}( \mathbf{x}, h ) = 0 $ for all $ h $ if $ R_n( \mathbf{x}, h ) = 0 $ for all $ h $.
From \autoref{lemma:adapt2}, we get that for all sufficiently large $ n $,
\begin{align}
& \mathbb{E}\left[ D_n( \mathbf{x}, h ) \right]
\le c_2 \frac{\log n}{n \phi( \mathbf{x}, h )}
\label{eq:adapt5}
\end{align}
for all $ h \in \mathbb{H}_n $, where $ c_2 > 0 $ is a constant independent of $ h $.
On the other hand, from decomposition \eqref{eq:main}, we get
\begin{align}
\left\| \widehat{\Theta}_n( \mathbf{x}, h ) - \Theta( \mathbf{x} ) \right\|^2
& \le
3 \left[ \left( \| B_n( \mathbf{x}, h ) \|^2 + M h^{2 \beta} \right) + 2 \| V_n( \mathbf{x}, h ) \|^2 \right] \nonumber\\
& \quad + 3 \left( \| R_n( \mathbf{x}, h ) \|^2 - \left( M h^{2 \beta} + \| V_n( \mathbf{x}, h ) \|^2 \right) \right)_+ ,
\end{align}
where $ M $ is the constant described in condition \ref{cond:adaptive3}.
From inequality \eqref{eq:b2} and \autoref{thm:1}, we have
\begin{align}
& \mathbb{E}\left[ \left( \| B_n( \mathbf{x}, h ) \|^2 + M h^{2 \beta} \right) + 2 \| V_n( \mathbf{x}, h ) \|^2 \right]
\le c_3 h^{2 \beta} + \frac{c_4}{n \phi( \mathbf{x}, h )}
\label{eq:adapt6}
\end{align}
for all sufficiently large $ n $ and some constants $ c_3 > 0 $ and $ c_4 > 0 $ independent of $ h $.
Also, from \autoref{lemma:adapt4}, we have
\begin{align}
& \max\limits_{h \in \mathbb{H}_n} \left( \| R_n( \mathbf{x}, h ) \|^2 - \left( M h^{2 \beta} + \| V_n( \mathbf{x}, h ) \|^2 \right) \right)_+
= o_\mathbb{P}\left( n^{-2} \right)
\label{eq:adapt7}
\end{align}
as $ n \longrightarrow \infty $.
Therefore, from \eqref{eq:adapt1}--\eqref{eq:adapt7}, we get that
\begin{align}
& \left\| \widehat{\Theta}_n\left( \mathbf{x}, h_n^* \right) - \Theta\left( \mathbf{x} \right) \right\|^2
\nonumber\\
& \le
\left[ 6 C_n^{(1)}( \mathbf{x}, h ) + 6 D_n( \mathbf{x}, h ) + 9 \left( \left( \| B_n( \mathbf{x}, h ) \|^2 + M h^{2 \beta} \right) + 2 \| V_n( \mathbf{x}, h ) \|^2 \right) \right] \nonumber\\
& \quad +
\left[ 6 C_n^{(2)}( \mathbf{x}, h ) + 9 \left( \| R_n( \mathbf{x}, h ) \|^2 - \left( M h^{2 \beta} + \| V_n( \mathbf{x}, h ) \|^2 \right) \right)_+ \right] , 
\nonumber\\
\intertext{where}
& \mathbb{E}\left[ 6 C_n^{(1)}( \mathbf{x}, h ) + 6 D_n( \mathbf{x}, h ) + 9 \left( \left( \| B_n( \mathbf{x}, h ) \|^2 + M h^{2 \beta} \right) + 2 \| V_n( \mathbf{x}, h ) \|^2 \right) \right] \nonumber\\
& \qquad
= O\left( h^{2 \beta} + \frac{\log n}{n \phi( \mathbf{x}, h )} \right)
\label{eq:adapt8}\\
\intertext{and}
& \max\limits_{h \in \mathbb{H}_n} \left[ 6 C_n^{(2)}( \mathbf{x}, h ) + 9 \left( \| R_n( \mathbf{x}, h ) \|^2 - \left( M h^{2 \beta} + \| V_n( \mathbf{x}, h ) \|^2 \right) \right)_+ \right] \nonumber\\
& \qquad
= o_\mathbb{P}\left( n^{-2} \right)
\quad
\text{as }
n \longrightarrow \infty .
\label{eq:adapt9}
\end{align}
Further, if $ R_n( \mathbf{x}, h ) = 0 $ for all $ h $, then
\begin{align}
& \max\limits_{h \in \mathbb{H}_n} \left[ 6 C_n^{(2)}( \mathbf{x}, h ) + 9 \left( \| R_n( \mathbf{x}, h ) \|^2 - \left( M h^{2 \beta} + \| V_n( \mathbf{x}, h ) \|^2 \right) \right)_+ \right]
= 0
\quad
\text{for all } h . \nonumber
\end{align}
From \eqref{eq:adapt8} and \eqref{eq:adapt9}, we get
\begin{align*}
& \left\| \widehat{\Theta}_n\left( \mathbf{x}, h_n^* \right) - \Theta\left( \mathbf{x} \right) \right\|^2
= O_\mathbb{P}\left( \lambda_n \right)
\quad
\text{as }
n \longrightarrow \infty ,
\end{align*}
and if $ R_n( \mathbf{x}, h ) = 0 $ for all $ h $, then
\begin{align*}
& \mathbb{E}\left\| \widehat{\Theta}_n\left( \mathbf{x}, h_n^* \right) - \Theta\left( \mathbf{x} \right) \right\|^2
= O\left( \lambda_n \right)
\quad
\text{as }
n \longrightarrow \infty .
\end{align*}
\end{proof}

\subsection{Results required to prove \autoref{thm:low}}

\begin{lemma} \label{lemma:thmlow4}
Let $ \{ \mathbf{U}_n \} $ be a sequence of real random variables and let $ \{ \mathbf{V}_n \} $ be another sequence of positive random variables with $ \mathbf{V}_n = o_\mathbb{P}( 1 ) $ as $ n \longrightarrow \infty $. Then, for any $ a > 0 $ and any $ \epsilon > 0 $, $ \mathbb{P}\left[ \mathbf{U}_n > a + \mathbf{V}_n \right]
> \mathbb{P}\left[ \mathbf{U}_n > 2 a \right]
- \epsilon $
for all sufficiently large $ n $.
\end{lemma}
\begin{proof}
Since $ \mathbf{V}_n = o_\mathbb{P}( 1 ) $ as $ n \longrightarrow \infty $, for any $ a > 0 $ and any $ \epsilon > 0 $,
\begin{align*}
\mathbb{P}\left[ \mathbf{U}_n > a + \mathbf{V}_n \right]
& \ge \mathbb{P}\left[ \mathbf{U}_n > 2 a \text{ and } \mathbf{V}_n < (a / 2) \right] \nonumber\\
& \ge \mathbb{P}\left[ \mathbf{U}_n > 2 a \right]
- \mathbb{P}\left[ \mathbf{V}_n > (a / 2) \right] \nonumber\\
& > \mathbb{P}\left[ \mathbf{U}_n > 2 a \right]
- \epsilon
\end{align*}
for all sufficiently large $ n $, which completes the proof.
\end{proof}

\begin{lemma} \label{lemma:thmlow1}
Suppose that in \eqref{eq:1}, we have either $ t_2 > 0 $, or $ t_3 > 1 $ with $ C_2 = C_4 $, the kernel $ K( \cdot ) $ satisfies \ref{assume:a1}, and the decomposition \eqref{eq:main} along with conditions \ref{cond:1}--\ref{cond:4}, \ref{cond:6} and \ref{cond:7} are satisfied.
Consider a bandwidth sequence $ \{ h_n \} $ that satisfies \ref{assume:a2} and $ h_n^{2 \beta} n \phi( \mathbf{x}, h_n ) \longrightarrow 0 $ as $ n \longrightarrow \infty $.
Then, there exist $ c > 0 $ and $ \delta > 0 $ such that
\begin{align*}
\mathbb{P}\left[ \left( m^{-1}( \log n ) \right)^{-\beta} \left\| \widehat{\Theta}_n( \mathbf{x} ) - \Theta( \mathbf{x} ) \right\| > c \right] > \delta
\end{align*}
for all sufficiently large $ n $.
\end{lemma}
\begin{proof}
Recall from \autoref{subsec:4_1} that $ B_n( \mathbf{x} ) = \tilde{B}_n( \mathbf{x} ) + \tilde{R}_n( \mathbf{x} ) $, where $ \tilde{R}_n( \mathbf{x} ) = o_\mathbb{P}( h_n^\beta ) $, and $ \tilde{B}_n( \mathbf{x} ) $ is a non-random quantity. So, from \eqref{eq:main} and condition \ref{cond:4}, we have
\begin{align}
\widehat{\Theta}_n( \mathbf{x} ) - \Theta( \mathbf{x} )
= \tilde{B}_n( \mathbf{x} ) + V_n( \mathbf{x} ) + Q_n( \mathbf{x} ),
\label{eq:q_n}
\end{align}
where $ Q_n( \mathbf{x} ) 
= R_n( \mathbf{x} ) + \tilde{R}_n( \mathbf{x} )
= o_\mathbb{P}\big( \max\big\{ h_n^\beta, \big[ n \phi( \mathbf{x}, h_n ) \big]^{-1/2} \big\} \big) $ as $ n \longrightarrow \infty $.

Recall the projection functional $ \tilde{\phi}_i( \cdot ) $ defined in \autoref{subsec:4_1} and the positive integer $ i_0 $ mentioned in condition \ref{cond:7}. Note that $ \| \tilde{\phi}_{i_0} \| = 1 $. So, for all $ \mathbf{v} \in \mathcal{B} $,
\begin{align}
| \tilde{\phi}_{i_0}( \mathbf{v} ) | \le \| v \| .
\label{eq:thmlow2}
\end{align}
Using \ref{assume:a1}, \ref{assume:a2}, \ref{cond:3}, \ref{cond:7} and arguments similar to those in \autoref{thm:2}, we get
\begin{align}
[ n \phi( \mathbf{x}, h_n ) ]^{1/2} [ E_n^{(2)}( \mathbf{x} ) ]^{-1/2} E_n^{(1)}( \mathbf{x} ) \tilde{\phi}_{i_0}( V_n( \mathbf{x} ) )
\longrightarrow
\mathbf{Z}
\label{eq:thmlow3}
\end{align}
\textit{in distribution} as $ n \longrightarrow \infty $, where $ \mathbf{Z} $ follows a normal distribution with mean zero and variance $ \mathbb{V}( \mathbf{x} ) > 0 $.

Next, consider $ \{ h_n \} $ that satisfies \ref{assume:a2} and
\begin{align}
h_n^{2 \beta} n \phi( \mathbf{x}, h_n ) \longrightarrow 0 \text{ as } n \longrightarrow \infty .
\label{eq:thm:1case1}
\end{align}
From \eqref{eq:thm:1a4} and \eqref{eq:thm:1case1}, we get that for all sufficiently large $ n $,
\begin{align}
& [ n \phi( \mathbf{x}, h_n ) ]^{-1/2} 
> h_n^{\beta}
> c_1^{\beta} \left( m^{-1}\left( \log n \right) \right)^\beta \nonumber\\
& \implies
\left( m^{-1}\left( \log n \right) \right)^{-\beta} [ n \phi( \mathbf{x}, h_n ) ]^{-1/2}
> c_1^{\beta} ,
\label{eq:thm:1case1a}
\end{align}
where $ c_1 > 0 $ is a constant.
Since $ Q_n( \mathbf{x} ) 
= o_\mathbb{P}\big( \max\big\{ h_n^\beta, \big[ n \phi( \mathbf{x}, h_n ) \big]^{-1/2} \big\} \big) $ as $ n \longrightarrow \infty $, from \eqref{eq:thm:1case1}, we have $ Q_n( \mathbf{x} ) 
= o_\mathbb{P}\big( \big[ n \phi( \mathbf{x}, h_n ) \big]^{-1/2} \big) $ as $ n \longrightarrow \infty $.
Further, from \ref{cond:1}, we get that $ h_n^{-\beta} \tilde{B}_n( \mathbf{x} ) $ is bounded, and hence from \eqref{eq:thm:1case1}, we have $ [ n \phi( \mathbf{x}, h_n ) ]^{1/2} \tilde{B}_n( \mathbf{x} ) \longrightarrow \mathbf{0} $ as $ n \longrightarrow \infty $. Therefore,
\begin{align}
[ n \phi( \mathbf{x}, h_n ) ]^{1/2} \left[ \big\| \tilde{B}_n( \mathbf{x} ) \big\| + \big\| Q_n( \mathbf{x} ) \big\| \right]
= o_\mathbb{P}\big( 1 \big)
\label{eq:thm:1case1b}
\end{align}
as $ n \longrightarrow 0 $.
Take
\begin{align*}
c = \frac{l c_1^{\beta}}{2 L} \quad
\text{and} \quad
\delta = \frac{1}{2} \mathbb{P}[ | \mathbf{Z} | > 1 ] ,
\end{align*}
where $ \mathbf{Z} $ is the normal random variable described in \eqref{eq:thmlow3}.
So, from \eqref{eq:assume:a2}, \autoref{lemma:thmlow4}, \eqref{eq:thmlow2}, \eqref{eq:thmlow3}, \eqref{eq:thm:1case1a}, \eqref{eq:thm:1case1b} and the triangle inequality, we have for all sufficiently large $ n $,
\begin{align*}
& \mathbb{P}\big[ \left( m^{-1}( \log n ) \right)^{-\beta} \big\| \widehat{\Theta}_n( \mathbf{x} ) - \Theta( \mathbf{x} ) \big\| > c \big] \\
& \ge \mathbb{P}\left[ \frac{ [ n \phi( \mathbf{x}, h_n ) ]^{1/2} \left[ \big\| V_n( \mathbf{x} ) \big\| - \big\| \tilde{B}_n( \mathbf{x} ) \big\| - \big\| Q_n( \mathbf{x} ) \big\| \right]}{\left( m^{-1}( \log n ) \right)^{\beta} [ n \phi( \mathbf{x}, h_n ) ]^{1/2}} > c \right] \\
& \ge \mathbb{P}\left[ [ n \phi( \mathbf{x}, h_n ) ]^{1/2} \left| \tilde{\phi}_{i_0}( V_n( \mathbf{x} ) ) \right| 
> \frac{c}{c_1^{\beta}} + [ n \phi( \mathbf{x}, h_n ) ]^{1/2} \left[ \big\| \tilde{B}_n( \mathbf{x} ) \big\| + \big\| Q_n( \mathbf{x} ) \big\| \right] \right] \\
& \ge \mathbb{P}\left[ [ n \phi( \mathbf{x}, h_n ) ]^{1/2} \left| \tilde{\phi}_{i_0}( V_n( \mathbf{x} ) ) \right| 
> 2 c c_1^{-\beta} \right] - \frac{\delta}{2} \\
& \ge \mathbb{P}\left[ \left| [ n \phi( \mathbf{x}, h_n ) ]^{1/2} [ E_n^{(2)}( \mathbf{x} ) ]^{-1/2} E_n^{(1)}( \mathbf{x} ) \tilde{\phi}_{i_0}( V_n( \mathbf{x} ) ) \right| 
> 2 c c_1^{-\beta} \frac{L}{l} \right] - \frac{\delta}{2} \\
& = \mathbb{P}\left[ \left| [ n \phi( \mathbf{x}, h_n ) ]^{1/2} [ E_n^{(2)}( \mathbf{x} ) ]^{-1/2} E_n^{(1)}( \mathbf{x} ) \tilde{\phi}_{i_0}( V_n( \mathbf{x} ) ) \right| 
> 1 \right] - \frac{\delta}{2}
> \delta .
\end{align*}
\end{proof}

\begin{lemma} \label{lemma:thmlow2}
Suppose that in \eqref{eq:1}, we have either $ t_2 > 0 $, or $ t_3 > 1 $ with $ C_2 = C_4 $, the kernel $ K( \cdot ) $ satisfies \ref{assume:a1}, and the decomposition \eqref{eq:main} along with conditions \ref{cond:1}--\ref{cond:4}, \ref{cond:6} and \ref{cond:7} are satisfied.
Consider a bandwidth sequence $ \{ h_n \} $ that satisfies \ref{assume:a2} and $ h_n^{2 \beta} n \phi( \mathbf{x}, h_n ) \longrightarrow \infty $ as $ n \longrightarrow \infty $.
Then, there exist $ c > 0 $ and $ \delta > 0 $ such that
\begin{align*}
\mathbb{P}\left[ \left( m^{-1}( \log n ) \right)^{-\beta} \left\| \widehat{\Theta}_n( \mathbf{x} ) - \Theta( \mathbf{x} ) \right\| > c \right] > \delta
\end{align*}
for all sufficiently large $ n $.
\end{lemma}
\begin{proof}
Consider $ \{ h_n \} $ that satisfies \ref{assume:a2} and
\begin{align}
h_n^{2 \beta} n \phi( \mathbf{x}, h_n ) \longrightarrow \infty
\text{ as } n \longrightarrow \infty .
\label{eq:thm:1case2}
\end{align}
Let $ Q_n( \mathbf{x} ) $ be as defined in \eqref{eq:q_n}.
Since $ Q_n( \mathbf{x} ) 
= o_\mathbb{P}\big( \max\big\{ h_n^\beta, \big[ n \phi( \mathbf{x}, h_n ) \big]^{-1/2} \big\} \big) $ as $ n \longrightarrow \infty $, from \eqref{eq:thm:1case2}, we have $ Q_n( \mathbf{x} ) 
= o_\mathbb{P}\big( h_n^\beta \big) $ as $ n \longrightarrow \infty $.
Further, from \autoref{thm:1} and \eqref{eq:thm:1case2}, we get
\begin{align*}
h_n^{- 2 \beta} \mathbb{E}[ \| V_n( \mathbf{x} ) \|^2 ]
= h_n^{- 2 \beta} \left[ n \phi( \mathbf{x}, h_n ) \right]^{-1} n \phi( \mathbf{x}, h_n ) \mathbb{E}[ \| V_n( \mathbf{x} ) \|^2 ]
\longrightarrow 0
\end{align*}
as $ n \longrightarrow \infty $, which implies $ h_n^{-\beta} V_n( \mathbf{x} ) = o_\mathbb{P}\big( 1 \big) $ as $ n \longrightarrow \infty $. Therefore,
\begin{align}
h_n^{-\beta} \left[ \big\| V_n( \mathbf{x} ) \big\| + \big\| Q_n( \mathbf{x} ) \big\| \right]
= o_\mathbb{P}\big( 1 \big)
\label{eq:thm:1case2a}
\end{align}
as $ n \longrightarrow \infty $.
Note that we have chosen $ \Theta( \mathbf{x} ) $ satisfying \ref{cond:6}, so that for any kernel $ K( \cdot ) $ satisfying \ref{assume:a1} and any sequence of bandwidths $ \{ h_n \} $ satisfying \ref{assume:a2}, we have for all sufficiently large $ n $,
\begin{align}
h_n^{-\beta} \| \tilde{B}_n( \mathbf{x} ) \| \ge b_1 > 0 ,
\label{eq:thm:1case2b}
\end{align}
where $ b_1 $ is a constant.
Take
\begin{align*}
c = \frac{b_1 c_1^{\beta}}{4} \quad
\text{and} \quad
\delta = \frac{1}{2} .
\end{align*}
Then, from \eqref{eq:thm:1a4}, \autoref{lemma:thmlow4}, \eqref{eq:thm:1case2a}, \eqref{eq:thm:1case2b} and the triangle inequality, we have for all sufficiently large $ n $,
\begin{align*}
& \mathbb{P}\big[ \left( m^{-1}( \log n ) \right)^{-\beta} \big\| \widehat{\Theta}_n( \mathbf{x} ) - \Theta( \mathbf{x} ) \big\| > c \big] \\
& \ge \mathbb{P}\left[ \frac{ h_n^{-\beta} \left[ \big\| \tilde{B}_n( \mathbf{x} ) \big\| - \big\| V_n( \mathbf{x} ) \big\| - \big\| Q_n( \mathbf{x} ) \big\| \right]}{\left( m^{-1}( \log n ) \right)^{\beta} h_n^{-\beta}} > c \right] \\
& \ge \mathbb{P}\left[ h_n^{-\beta} \big\| \tilde{B}_n( \mathbf{x} ) \big\| > c c_1^{-\beta} + h_n^{-\beta} \left[ \big\| V_n( \mathbf{x} ) \big\| + \big\| Q_n( \mathbf{x} ) \big\| \right] \right] \\
& \ge \mathbb{P}\left[ h_n^{-\beta} \big\| \tilde{B}_n( \mathbf{x} ) \big\| > 2 c c_1^{-\beta} \right] - \frac{1}{4} \\
& = \mathbb{P}\left[ h_n^{-\beta} \big\| \tilde{B}_n( \mathbf{x} ) \big\| > \frac{b_1}{2} \right]  - \frac{1}{4}
= \frac{3}{4}
> \delta .
\end{align*}
\end{proof}

\begin{lemma} \label{lemma:thmlow3}
Suppose that in \eqref{eq:1}, we have either $ t_2 > 0 $, or $ t_3 > 1 $ with $ C_2 = C_4 $, the kernel $ K( \cdot ) $ satisfies \ref{assume:a1}, and the decomposition \eqref{eq:main} along with conditions \ref{cond:1}--\ref{cond:4}, \ref{cond:6} and \ref{cond:7} are satisfied.
Consider a bandwidth sequence $ \{ h_n \} $ that satisfies \ref{assume:a2}, and $ 0 < \epsilon_1
< h_n^{2 \beta} n \phi( \mathbf{x}, h_n )
< \epsilon_2 $
for all sufficiently large $ n $ and some $ \epsilon_1 $ and $ \epsilon_2 $.
Then, there exist $ c > 0 $ and $ \delta > 0 $ such that
\begin{align*}
\mathbb{P}\left[ \left( m^{-1}( \log n ) \right)^{-\beta} \left\| \widehat{\Theta}_n( \mathbf{x} ) - \Theta( \mathbf{x} ) \right\| > c \right] > \delta
\end{align*}
for all sufficiently large $ n $.
\end{lemma}
\begin{proof}
Consider $ \{ h_n \} $ that satisfies \ref{assume:a2} and
\begin{align}
0 < \epsilon_1
< h_n^{2 \beta} n \phi( \mathbf{x}, h_n )
< \epsilon_2
\label{eq:thm:1case3}
\end{align}
for all sufficiently large $ n $ and some $ \epsilon_1 $ and $ \epsilon_2 $.
From \eqref{eq:thm:1a4} and \eqref{eq:thm:1case3}, we get that for all sufficiently large $ n $,
\begin{align}
& \left( m^{-1}\left( \log n \right) \right)^{\beta} [ n \phi( \mathbf{x}, h_n ) ]^{1/2}
< \frac{[ n \phi( \mathbf{x}, h_n ) ]^{1/2} h_n^{\beta}}{c_1^{\beta}}
< \frac{\sqrt{\epsilon_2}}{c_1^{\beta}} ,
\label{eq:thm:1case3a}
\end{align}
where $ c_1 > 0 $ is a constant.
Let $ Q_n( \mathbf{x} ) $ be as defined in \eqref{eq:q_n}.
Since $ Q_n( \mathbf{x} ) 
= o_\mathbb{P}\big( \max\big\{ h_n^\beta, \big[ n \phi( \mathbf{x}, h_n ) \big]^{-1/2} \big\} \big) $ as $ n \longrightarrow \infty $, from \eqref{eq:thm:1case3}, we have
\begin{align}
& \max\big\{ h_n^\beta, \big[ n \phi( \mathbf{x}, h_n ) \big]^{-1/2} \big\}
\le \max\{ \sqrt{\epsilon_2}, 1 \} \big[ n \phi( \mathbf{x}, h_n ) \big]^{-1/2} \nonumber\\
\implies
& [ n \phi( \mathbf{x}, h_n ) ]^{1/2} \left\| Q_n( \mathbf{x} ) \right\|
= o_\mathbb{P}\big( 1 \big)
\label{eq:thm:1case3b}
\end{align}
as $ n \longrightarrow \infty $. From \ref{assume:a2}, \ref{cond:1} and \eqref{eq:thm:1case3}, we get
\begin{align}
& [ n \phi( \mathbf{x}, h_n ) ]^{1/2} \big\| \tilde{B}_n( \mathbf{x} ) \big\|
\le [ n \phi( \mathbf{x}, h_n ) ]^{1/2} h_n^{\beta} h_n^{- \beta} \big\| \tilde{B}_n( \mathbf{x} ) \big\|
\le \sqrt{\epsilon_2} \| \mathbb{L}_\mathbf{x} \| b_F
\label{eq:thm:1case3c}
\end{align}
for all sufficiently large $ n $.
Take
\begin{align*}
c = \frac{c_1^{\beta} l}{2 \sqrt{\epsilon_2} L} \quad
\text{and} \quad
\delta = \frac{1}{2} \mathbb{P}\left[ | \mathbf{Z} | > 1 + \sqrt{\epsilon_2} \frac{L}{l} \| \mathbb{L}_\mathbf{x} \| b_F \right] ,
\end{align*}
where $ \mathbf{Z} $ is the normal random variable described in \eqref{eq:thmlow3}, and $ l $ and $ L $ are the constants described in \ref{assume:a1}.
So, from \eqref{eq:assume:a2}, \autoref{lemma:thmlow4}, \eqref{eq:thmlow2}, \eqref{eq:thmlow3}, \eqref{eq:thm:1case3a}, \eqref{eq:thm:1case3b}, \eqref{eq:thm:1case3c} and the triangle inequality, we have for all sufficiently large $ n $,
\begin{align*}
& \mathbb{P}\big[ \left( m^{-1}( \log n ) \right)^{-\beta} \big\| \widehat{\Theta}_n( \mathbf{x} ) - \Theta( \mathbf{x} ) \big\| > c \big] \\
& \ge \mathbb{P}\left[ \frac{ [ n \phi( \mathbf{x}, h_n ) ]^{1/2} \left[ \big\| V_n( \mathbf{x} ) \big\| - \big\| \tilde{B}_n( \mathbf{x} ) \big\| - \big\| Q_n( \mathbf{x} ) \big\| \right]}{\left( m^{-1}( \log n ) \right)^{\beta} [ n \phi( \mathbf{x}, h_n ) ]^{1/2}} > c \right] \\
& \ge \mathbb{P}\left[ [ n \phi( \mathbf{x}, h_n ) ]^{1/2} \left| \tilde{\phi}_{i_0}\left( V_n( \mathbf{x} ) \right) \right| > c \frac{\sqrt{\epsilon_2}}{c_1^{\beta}} 
+ \sqrt{\epsilon_2} \| \mathbb{L}_\mathbf{x} \| b_F
+ [ n \phi( \mathbf{x}, h_n ) ]^{1/2} \big\| Q_n( \mathbf{x} ) \big\| \right] \\
& \ge \mathbb{P}\left[ [ n \phi( \mathbf{x}, h_n ) ]^{1/2} \frac{E_n^{(1)}( \mathbf{x} )}{[ E_n^{(2)}( \mathbf{x} ) ]^{1/2}} \left| \tilde{\phi}_{i_0}\left( V_n( \mathbf{x} ) \right) \right| 
> \frac{2 \sqrt{\epsilon_2} L}{c_1^{\beta} l} c
+ \sqrt{\epsilon_2} \frac{L}{l} \| \mathbb{L}_\mathbf{x} \| b_F \right] - \frac{\delta}{2} \\
& \ge \mathbb{P}\left[ \left| [ n \phi( \mathbf{x}, h_n ) ]^{1/2} \frac{E_n^{(1)}( \mathbf{x} )}{[ E_n^{(2)}( \mathbf{x} ) ]^{1/2}} \tilde{\phi}_{i_0}\left( V_n( \mathbf{x} ) \right) \right| 
> 1 + \sqrt{\epsilon_2} \frac{L}{l} \| \mathbb{L}_\mathbf{x} \| b_F \right] - \frac{\delta}{2} \\
& > \delta .
\end{align*}
\end{proof}

\subsection{Results required to prove \autoref{coro:1}}

\begin{lemma} \label{lemma:coro6}
Suppose assumptions \ref{assume:a1} and \ref{assume:a2} are satisfied.
Let $ \{ h_n^{(b)} \} $ be a sequence of bandwidths that satisfies \ref{assume:a2} and balances the bias and the variance so that
\begin{align}
0 < c_1 \le ( h_n^{(b)} )^{2 \beta} n \phi( \mathbf{x}, h_n^{(b)} ) \le c_2 < \infty
\label{eq:coro1}
\end{align}
for all sufficiently large $ n $, where $ c_1, c_2 $ are some constants.
Also, let $ \{ h_n^{(op)} \} $ denote the sequence of optimum bandwidths minimizing \eqref{eq:thm:up} in the proof of \autoref{thm:up}.
Assume that $ t_2 > 0 $ in the bounds on the small ball probability of the covariate in \eqref{eq:1}.
Then,
\begin{align*}
0 < c_3 \le \frac{h_n^{(b)}}{h_n^{(op)}} \le c_4 < \infty
\end{align*}
for all sufficiently large $ n $, where $ c_3, c_4 $ are some constants.
\end{lemma}
\begin{proof}
Recall from \eqref{eq:mh} that
$ m( h ) = C_2 ( 1 / h )^{t_2} ( \log ( 1 / h ) )^{t_3} $ for $ 0 < h < 1 $.
From \eqref{eq:1} and \eqref{eq:coro1}, we have
\begin{align}
& ( h_n^{(b)} )^{2 \beta + t_1} n C_1 \exp\left[ -m( h_n^{(b)} ) \right]
\le c_2 \nonumber\\
& \text{and}\quad
c_1 \le ( h_n^{(b)} )^{2 \beta + t_4} n C_3 \exp\left[ -(C_4 / C_2) m( h_n^{(b)} ) \right] \nonumber\\
\implies
& ( h_n^{(b)} )^{2 \beta + t_1} n \exp\left[ -m( h_n^{(b)} ) \right]
\le \frac{c_2}{C_1} \nonumber\\
& \text{and}\quad
\frac{c_1}{C_3} \le ( h_n^{(b)} )^{2 \beta + t_4} n \exp\left[ -(C_4 / C_2) m( h_n^{(b)} ) \right] \nonumber\\
\implies
& \frac{- ( 2 \beta + t_1 ) \log \frac{1}{h_n^{(b)}}}{m( h_n^{(b)} )} + \frac{\log n}{m( h_n^{(b)} )} - 1
\le \frac{\log \frac{c_2}{C_1}}{m( h_n^{(b)} )} \nonumber\\
& \text{and}\quad
\frac{\log \frac{c_1}{C_3}}{m( h_n^{(b)} )} \le \frac{- ( 2 \beta + t_4 ) \log \frac{1}{h_n^{(b)}}}{m( h_n^{(b)} )} + \frac{\log n}{m( h_n^{(b)} )} -\frac{C_4}{C_2}
\label{eq:coro2}
\end{align}
for all sufficiently large $ n $. When $ t_2 > 0 $ in \eqref{eq:1}, we have
\begin{align*}
& \frac{- ( 2 \beta + t_1 ) \log \frac{1}{h_n^{(b)}}}{m( h_n^{(b)} )} \longrightarrow 0 , \quad
\frac{\log \frac{c_2}{C_1}}{m( h_n^{(b)} )} \longrightarrow 0 , \\
& \frac{\log \frac{c_1}{C_3}}{m( h_n^{(b)} )} \longrightarrow 0 \quad \text{and} \quad
\frac{- ( 2 \beta + t_4 ) \log \frac{1}{h_n^{(b)}}}{m( h_n^{(b)} )} \longrightarrow 0
\end{align*}
as $ n \longrightarrow \infty $. Therefore, given any $ \epsilon > 0 $, from \eqref{eq:coro2}, we have for all sufficiently large $ n $,
\begin{align}
& \frac{\log n}{m( h_n^{(b)} )}
\le 1 + \epsilon \quad
\text{and} \quad
\frac{C_4}{C_2} - \epsilon \le \frac{\log n}{m( h_n^{(b)} )} \nonumber\\
\implies
& \frac{\log n}{1 + \epsilon}
\le m( h_n^{(b)} ) 
\le \frac{\log n}{(C_4 / C_2) - \epsilon} \nonumber\\
\implies
& m^{-1}\left( \frac{\log n}{1 + \epsilon} \right)
\ge h_n^{(b)}
\ge m^{-1}\left( \frac{\log n}{(C_4 / C_2) - \epsilon} \right) .
\label{eq:coro3}
\end{align}
Next, we consider our optimum bandwidth $ h_n^{(op)} $. From \eqref{eq:thm:upg} in the proof of \autoref{thm:up}, we have, given any $ \epsilon > 0 $ and for all sufficiently large $ n $,
\begin{align}
m^{-1}\left( \frac{\log n}{1 + \epsilon} \right)
\ge h_n^{(op)}
\ge m^{-1}\left( \frac{\log n}{1 - \epsilon} \right) .
\label{eq:coro4}
\end{align}
Since $ m( h ) $ is strictly monotone decreasing function for $ h \in ( 0, 1 ) $ and $ m( h ) \longrightarrow \infty $ as $ h \longrightarrow 0^+ $, $ m^{-1}( u ) $ is well-defined for all $ u > 1 $ and $ m^{-1}( u ) \longrightarrow 0^+ $ as $ u \longrightarrow \infty $.
Given $ \epsilon > 0 $, we have
\begin{align*}
& m\left( c^{- \frac{1}{t_2}} ( 1 + \epsilon ) h \right)
= c m( h ) \frac{1}{( 1 + \epsilon )^{t_2}} \left( 1 - \frac{\log( 1 + \epsilon )}{\log \frac{1}{h}} + \frac{1}{t_2} \frac{\log c}{\log \frac{1}{h}} \right)^{t_3} , 
\\
& m\left( c^{- \frac{1}{t_2}} ( 1 - \epsilon ) h \right)
= c m( h ) \frac{1}{( 1 - \epsilon )^{t_2}} \left( 1 - \frac{\log( 1 - \epsilon )}{\log \frac{1}{h}} + \frac{1}{t_2} \frac{\log c}{\log \frac{1}{h}} \right)^{t_3} .
\end{align*}
For sufficiently small $ h > 0 $, we have
\begin{align*}
& \frac{1}{( 1 + \epsilon )^{t_2}} \left( 1 - \frac{\log( 1 + \epsilon )}{\log \frac{1}{h}} + \frac{1}{t_2} \frac{\log c}{\log \frac{1}{h}} \right)^{t_3} \\
& < 1 \\
& < \frac{1}{( 1 - \epsilon )^{t_2}} \left( 1 - \frac{\log( 1 - \epsilon )}{\log \frac{1}{h}} + \frac{1}{t_2} \frac{\log c}{\log \frac{1}{h}} \right)^{t_3} ,
\end{align*}
which implies
\begin{align*}
& m\left( c^{- \frac{1}{t_2}} ( 1 + \epsilon ) h \right)
< c m( h )
< m\left( c^{- \frac{1}{t_2}} ( 1 - \epsilon ) h \right)
\end{align*}
for all sufficiently small $ h > 0 $. Hence, for all sufficiently large $ u $, we have
\begin{align}
& m\left( c^{- \frac{1}{t_2}} ( 1 + \epsilon ) m^{-1}( u ) \right)
< c u
< m\left( c^{- \frac{1}{t_2}} ( 1 - \epsilon ) m^{-1}( u ) \right)
\nonumber\\
& \implies
c^{- \frac{1}{t_2}} ( 1 - \epsilon )
< \frac{m^{-1}( c u )}{m^{-1}( u )}
< c^{- \frac{1}{t_2}} ( 1 + \epsilon ) .
\label{eq:coro5}
\end{align}
From \eqref{eq:coro5}, we get that for any $ c > 0 $,
\begin{align}
\frac{m^{-1}( c u )}{m^{-1}( u )} \longrightarrow c^{- \frac{1}{t_2}}
\label{eq:coro6}
\end{align}
as $ u \longrightarrow \infty $.
Therefore, using \eqref{eq:coro3}, \eqref{eq:coro4} and \eqref{eq:coro6}, we have
\begin{align*}
0 < c_3 \le \frac{h_n^{(b)}}{h_n^{(op)}} \le c_4 < \infty
\end{align*}
for all sufficiently large $ n $, where $ c_3, c_4 $ are some constants.
\end{proof}

\begin{lemma} \label{lemma:newlemma1}
We denote our optimum bandwidth minimizing \eqref{eq:thm:up} in the proof of \autoref{thm:up} as $ h_n^{(op)} $.
Let $ \widehat{\Theta}_n^{(op)}( \mathbf{x} ) $ be as defined in \autoref{coro:1}.
Then, under the conditions in \autoref{coro:1},
\begin{align*}
& ( h_n^{(op)} )^{-\beta} \left\| \widehat{\Theta}_n^{(op)}( \mathbf{x} ) - \Theta( \mathbf{x} ) \right\|
= o_\mathbb{P}( 1 )
\quad
\text{as } n \longrightarrow \infty ,
\\
& \text{and}\quad
( h_n^{(op)} )^{- 2 \beta} \mathbb{E}\left\| \widehat{\Theta}_n^{(op)}( \mathbf{x} ) - \Theta( \mathbf{x} ) \right\|^2
\longrightarrow 0
\quad
\text{as } n \longrightarrow \infty .
\end{align*}
\end{lemma}
\begin{proof}
From \eqref{eq:proof1} in the proof of \autoref{thm:up} and the lower bound of $ \phi( \mathbf{x}, h ) $ in \eqref{eq:1}, we get
\begin{align}
(h_n^{(op)})^{2 \beta} n \phi( \mathbf{x}, h_n^{(op)} ) \longrightarrow \infty
\text{ as }
n \longrightarrow \infty .
\label{eq:coro7}
\end{align}
Since $ F( \cdot ) \in \mathcal{F}( \mathbf{x}, \beta_1, \mathcal{G} ) $ for some $ \beta_1 > \beta $, we have
\begin{align*}
\left( d( \mathbf{x}, \mathbf{z} ) \right)^{- \beta} \| F( \mathbf{z} ) - F( \mathbf{x} ) \| \longrightarrow 0
\quad
\text{as }
d( \mathbf{x}, \mathbf{z} ) \longrightarrow 0 .
\end{align*}
Consequently,
\begin{align}
& ( h_n^{(op)} )^{- \beta} \left\| B_n^{(op)}( \mathbf{x} ) \right\|
= o_\mathbb{P}( 1 )
\quad
\text{as }
n \longrightarrow \infty , 
\label{eq:coro9a}\\
& \text{and } ( h_n^{(op)} )^{- 2 \beta} \mathbb{E}\left\| B_n^{(op)}( \mathbf{x} ) \right\|^2
\longrightarrow 0 
\quad
\text{as }
n \longrightarrow \infty .
\label{eq:coro11}
\end{align}
From \autoref{thm:1} and \eqref{eq:coro7}, we have
\begin{align}
& ( h_n^{(op)} )^{-2 \beta} \mathbb{E}\left\| V_n^{(op)}( \mathbf{x} ) \right\|^2 \nonumber\\
& = \left( ( h_n^{(op)} )^{-2 \beta} \left( n \phi( \mathbf{x}, h_n^{(op)} ) \right)^{- 1} \right)
n \phi( \mathbf{x}, h_n^{(op)} ) \mathbb{E}\left\| V_n^{(op)}( \mathbf{x} ) \right\|^2
\longrightarrow 0
\label{eq:coro15}
\end{align}
as $ n \longrightarrow \infty $, and from \eqref{eq:coro15} and the Markov inequality, we get
\begin{align}
& ( h_n^{(op)} )^{-\beta} \left\| V_n^{(op)}( \mathbf{x} ) \right\| 
= o_\mathbb{P}( 1 )
\text{ as }
n \longrightarrow \infty .
\label{eq:coro16}
\end{align}
From condition \ref{cond:4} and \eqref{eq:coro7}, we have
\begin{align}
& ( h_n^{(op)} )^{-\beta} \left\| R_n^{(op)}( \mathbf{x} ) \right\|
= o_\mathbb{P}( 1 )
\quad
\text{as }
n \longrightarrow \infty .
\label{eq:coro17a}
\end{align}
Since $ \| \tilde{\phi}_{i_0} \| = 1 $,
when $ \mathbb{E}[ \| R_n( \mathbf{x} ) \|^2 ] = o( \delta_n^2 ) $ as $ n \longrightarrow \infty $, from \eqref{eq:coro7}, we have
\begin{align}
& ( h_n^{(op)} )^{-2 \beta} \mathbb{E}\left\| R_n^{(op)}( \mathbf{x} ) \right\|^2
\longrightarrow 0
\quad
\text{as }
n \longrightarrow \infty .
\label{eq:coro19}
\end{align}
Therefore, from \eqref{eq:coro9a}, \eqref{eq:coro16} and \eqref{eq:coro17a}, we have
\begin{align*}
& ( h_n^{(op)} )^{-\beta} \left\| \widehat{\Theta}_n^{(op)}( \mathbf{x} ) - \Theta( \mathbf{x} ) \right\|
\\
& \le ( h_n^{(op)} )^{-\beta} \left\| B_n^{(op)}( \mathbf{x} ) \right\|
+ ( h_n^{(op)} )^{-\beta} \left\| V_n^{(op)}( \mathbf{x} ) \right\|
+ ( h_n^{(op)} )^{-\beta} \left\| R_n^{(op)}( \mathbf{x} ) \right\|
\\
& = o_\mathbb{P}( 1 )
\quad
\text{as } n \longrightarrow \infty .
\end{align*}
Further, from \eqref{eq:coro11}, \eqref{eq:coro15} and \eqref{eq:coro19}, we have
\begin{align*}
& ( h_n^{(op)} )^{- 2 \beta} \mathbb{E}\left\| \widehat{\Theta}_n^{(op)}( \mathbf{x} ) - \Theta( \mathbf{x} ) \right\|^2
\\
& \le 3 ( h_n^{(op)} )^{-2 \beta} \mathbb{E}\left\| B_n^{(op)}( \mathbf{x} ) \right\|^2
+ 3 ( h_n^{(op)} )^{-2 \beta} \mathbb{E}\left\| V_n^{(op)}( \mathbf{x} ) \right\|^2 \\
& \quad
+ 3 ( h_n^{(op)} )^{-2 \beta} \mathbb{E}\left\| R_n^{(op)}( \mathbf{x} ) \right\|^2
\\
& \longrightarrow 0
\quad
\text{as } n \longrightarrow \infty .
\end{align*}
\end{proof}

\begin{lemma} \label{lemma:newlemma2}
Let $ h_n^{(b)} $ and $ \widehat{\Theta}_n^{(b)}( \mathbf{x} ) $ be as defined in \autoref{coro:1}.
Then, under the conditions in \autoref{coro:1}, given any $ \epsilon > 0 $, there is $ \delta > 0 $ such that
\begin{align*}
\mathbb{P}\left[ ( h_n^{(b)} )^{-\beta} \left\| \widehat{\Theta}_n^{(b)}( \mathbf{x} ) - \Theta( \mathbf{x} ) \right\| > \delta \right]
> 1 - \epsilon
\end{align*}
for all sufficiently large $ n $.
Further,
\begin{align*}
( h_n^{(b)} )^{- 2 \beta} \mathbb{E}\left\| \widehat{\Theta}_n^{(b)}( \mathbf{x} ) - \Theta( \mathbf{x} ) \right\|^2
\text{ is bounded away from 0 as } n \longrightarrow \infty .
\end{align*}
\end{lemma}
\begin{proof}
Let $ h_n^{(b)} $ satisfy \eqref{eq:coro1}.
Since $ F( \cdot ) \in \mathcal{F}( \mathbf{x}, \beta_1, \mathcal{G} ) $ for some $ \beta_1 > \beta $, we have
\begin{align}
\left( d( \mathbf{x}, \mathbf{z} ) \right)^{- \beta} \| F( \mathbf{z} ) - F( \mathbf{x} ) \| \longrightarrow 0
\quad
\text{as }
d( \mathbf{x}, \mathbf{z} ) \longrightarrow 0 .
\label{eq:neweq1}
\end{align}
Consequently,
\begin{align}
& ( h_n^{(b)} )^{- \beta} \left\| B_n^{(b)}( \mathbf{x} ) \right\|
= o_\mathbb{P}( 1 )
\quad
\text{as }
n \longrightarrow \infty .
\label{eq:coro9}
\end{align}
Let $ \mathbf{Z} $ be the normal random variable described in \eqref{eq:thmlow3}. Given any $ \epsilon > 0 $, there exists $ \delta > 0 $ such that
\begin{align}
\mathbb{P}\left[ \left| \mathbf{Z} \right| > 2 \delta \sqrt{c_2} l^{-1} L \right]
> 1 - \epsilon ,
\label{eq:coro11a}
\end{align}
where $ c_2 $ is a constant described in \eqref{eq:coro1}, and $ l, L $ are constants described in assumption \ref{assume:a1}.
Hence, from \eqref{eq:assume:a2}, \eqref{eq:thmlow3}, \eqref{eq:coro1} and \eqref{eq:coro11a}, we have
\begin{align}
& \mathbb{P}\left[ \left\| ( h_n^{(b)} )^{-\beta} V_n^{(b)}( \mathbf{x} ) \right\| > 2 \delta \right] \nonumber\\
& = \mathbb{P}\left[ \left( ( h_n^{(b)} )^{-\beta} \left( n \phi( \mathbf{x}, h_n^{(b)} ) \right)^{- 1/2} \right) \left\| \left( n \phi( \mathbf{x}, h_n^{(b)} ) \right)^{1/2} V_n^{(b)}( \mathbf{x} ) \right\| > 2 \delta \right] \nonumber\\
& \ge \mathbb{P}\left[ \left\| \left( n \phi( \mathbf{x}, h_n^{(b)} ) \right)^{1/2} V_n^{(b)}( \mathbf{x} ) \right\| > 2 \delta \sqrt{c_2} \right] \nonumber\\
& \ge \mathbb{P}\left[ \left| \left( n \phi( \mathbf{x}, h_n^{(b)} ) \right)^{1/2} \tilde{\phi}_{i_0}( V_n^{(b)}( \mathbf{x} ) )  \right| > 2 \delta \sqrt{c_2} \right] \nonumber\\
& \ge \mathbb{P}\left[ \left| \left( n \phi( \mathbf{x}, h_n^{(b)} ) \right)^{1/2} [ E_n^{(2)}( \mathbf{x} ) ]^{-1/2} E_n^{(1)}( \mathbf{x} ) \tilde{\phi}_{i_0}( V_n^{(b)}( \mathbf{x} ) )  \right| > 2 \delta \sqrt{c_2} l^{-1} L \right] \nonumber\\
& > 1 - \epsilon
\label{eq:coro12}
\end{align}
for all sufficiently large $ n $.
From condition \ref{cond:4} and \eqref{eq:coro1}, we have
\begin{align}
& ( h_n^{(b)} )^{-\beta} \left\| R_n^{(b)}( \mathbf{x} ) \right\|
= o_\mathbb{P}( 1 )
\quad
\text{as }
n \longrightarrow \infty .
\label{eq:coro17}
\end{align}
Therefore, from \autoref{lemma:thmlow4}, \eqref{eq:coro9}, \eqref{eq:coro12} and \eqref{eq:coro17}, we have
\begin{align*}
& \mathbb{P}\left[ ( h_n^{(b)} )^{-\beta} \left\| \widehat{\Theta}_n^{(b)}( \mathbf{x} ) - \Theta( \mathbf{x} ) \right\| > \delta \right]
\\
& \ge
\mathbb{P}\left[ ( h_n^{(b)} )^{-\beta} \left\| V_n^{(b)}( \mathbf{x} ) \right\|
- ( h_n^{(b)} )^{-\beta} \left\| B_n^{(b)}( \mathbf{x} ) \right\|
- ( h_n^{(b)} )^{-\beta} \left\| R_n^{(b)}( \mathbf{x} ) \right\|
> \delta \right]
\\
& =
\mathbb{P}\left[ ( h_n^{(b)} )^{-\beta} \left\| V_n^{(b)}( \mathbf{x} ) \right\|
> \delta + ( h_n^{(b)} )^{-\beta} \left\| B_n^{(b)}( \mathbf{x} ) \right\|
+ ( h_n^{(b)} )^{-\beta} \left\| R_n^{(b)}( \mathbf{x} ) \right\| \right]
\\
& >
\mathbb{P}\left[ ( h_n^{(b)} )^{-\beta} \left\| V_n^{(b)}( \mathbf{x} ) \right\|
> 2 \delta \right]
> 1 - \epsilon
\end{align*}
for all sufficiently large $ n $.

We proceed to prove the second part of the lemma.
Since $ | \tilde{\phi}_{i_0}( \mathbf{v} ) | \le \| \mathbf{v} \| $ for any $ \mathbf{v} $, from an application of the Cauchy-Schwarz inequality, we have
\begin{align}
& \mathbb{E}\left\| \widehat{\Theta}_n^{(b)}( \mathbf{x} ) - \Theta( \mathbf{x} ) \right\|^2 \nonumber\\
& = \mathbb{E}\left\| B_n^{(b)}( \mathbf{x} ) + V_n^{(b)}( \mathbf{x} ) + R_n^{(b)}( \mathbf{x} ) \right\|^2 \nonumber\\
& \ge \mathbb{E}\left[ \tilde{\phi}_{i_0}( B_n^{(b)}( \mathbf{x} ) ) + \tilde{\phi}_{i_0}( V_n^{(b)}( \mathbf{x} ) ) + \tilde{\phi}_{i_0}( R_n^{(b)}( \mathbf{x} ) ) \right]^2 \nonumber\\
& = \mathbb{E}\left[ \left( \tilde{\phi}_{i_0}( B_n^{(b)}( \mathbf{x} ) ) \right)^2 \right]
+ \mathbb{E}\left[ \left( \tilde{\phi}_{i_0}( V_n^{(b)}( \mathbf{x} ) ) \right)^2 \right]
+ \mathbb{E}\left[ \left( \tilde{\phi}_{i_0}( R_n^{(b)}( \mathbf{x} ) ) \right)^2 \right] \nonumber\\
& \quad
+ 2 \mathbb{E}\left[ \tilde{\phi}_{i_0}( R_n^{(b)}( \mathbf{x} ) ) \left( \tilde{\phi}_{i_0}( B_n^{(b)}( \mathbf{x} ) ) + \tilde{\phi}_{i_0}( V_n^{(b)}( \mathbf{x} ) ) \right) \right] \nonumber\\
& \ge \mathbb{E}\left[ \left( \tilde{\phi}_{i_0}( B_n^{(b)}( \mathbf{x} ) ) \right)^2 \right]
+ \mathbb{E}\left[ \left( \tilde{\phi}_{i_0}( V_n^{(b)}( \mathbf{x} ) ) \right)^2 \right]
+ \mathbb{E}\left[ \left( \tilde{\phi}_{i_0}( R_n^{(b)}( \mathbf{x} ) ) \right)^2 \right] \nonumber\\
& \quad
- 2 \left[ \mathbb{E}\left[ \left( \tilde{\phi}_{i_0}( R_n^{(b)}( \mathbf{x} ) ) \right)^2 \right] \right]^{1/2}
\left[ \mathbb{E}\left[ \left( \tilde{\phi}_{i_0}( B_n^{(b)}( \mathbf{x} ) ) \right)^2 \right]
+ \mathbb{E}\left[ \left( \tilde{\phi}_{i_0}( V_n^{(b)}( \mathbf{x} ) ) \right)^2 \right] \right]^{1/2} .
\label{eq:coro8}
\end{align}
From \eqref{eq:neweq1}, we have
\begin{align}
& ( h_n^{(b)} )^{- 2 \beta} \mathbb{E}\left[ \left( \tilde{\phi}_{i_0}( B_n^{(b)}( \mathbf{x} ) ) \right)^2 \right]
\le
( h_n^{(b)} )^{- 2 \beta} \mathbb{E}\left\| B_n^{(b)}( \mathbf{x} ) \right\|^2
\longrightarrow 0 
\label{eq:coro10}
\end{align}
as $ n \longrightarrow \infty $.
From \eqref{eq:assume:a2}, \eqref{eq:thmlow3} and \eqref{eq:coro1}, we have
\begin{align}
& ( h_n^{(b)} )^{-2 \beta} \mathbb{E}\left[ \left( \tilde{\phi}_{i_0}( V_n^{(b)}( \mathbf{x} ) ) \right)^2 \right] \nonumber\\
& = \left( ( h_n^{(b)} )^{-2 \beta} \left( n \phi( \mathbf{x}, h_n^{(b)} ) \right)^{- 1} \right)
n \phi( \mathbf{x}, h_n^{(b)} ) \mathbb{E}\left[ \left( \tilde{\phi}_{i_0}( V_n^{(b)}( \mathbf{x} ) ) \right)^2 \right] \nonumber\\
& \ge \frac{1}{c_2} \, \mathbb{P}\left[ \left| \left( n \phi( \mathbf{x}, h_n^{(b)} ) \right)^{1/2} [ E_n^{(2)}( \mathbf{x} ) ]^{-1/2} E_n^{(1)}( \mathbf{x} ) \tilde{\phi}_{i_0}( V_n^{(b)}( \mathbf{x} ) )  \right| > l^{-1} L \right]
\nonumber\\
& > c_6 > 0
\label{eq:coro13}
\end{align}
for all sufficiently large $ n $ and for some constant $ c_6 $.
Further, since $ \| \tilde{\phi}_{i_0} \| = 1 $, from \autoref{thm:1} and \eqref{eq:coro1}, we have
\begin{align}
& ( h_n^{(b)} )^{-2 \beta} \mathbb{E}\left[ \left( \tilde{\phi}_{i_0}( V_n^{(b)}( \mathbf{x} ) ) \right)^2 \right] \nonumber\\
& \le \left( ( h_n^{(b)} )^{-2 \beta} \left( n \phi( \mathbf{x}, h_n^{(b)} ) \right)^{- 1} \right)
n \phi( \mathbf{x}, h_n^{(b)} ) \mathbb{E}\left\| V_n^{(b)}( \mathbf{x} ) ) \right\|^2 
\le \frac{c_7}{c_1}
\label{eq:coro14}
\end{align}
for some constant $ c_7 > 0 $ and for all sufficiently large $ n $.
Since $ \| \tilde{\phi}_{i_0} \| = 1 $,
when $ \mathbb{E}[ \| R_n( \mathbf{x} ) \|^2 ] = o( \delta_n^2 ) $ as $ n \longrightarrow \infty $, from \eqref{eq:coro1} and \eqref{eq:coro7}, we have
\begin{align}
& ( h_n^{(b)} )^{-2 \beta} \mathbb{E}\left[ \left( \tilde{\phi}_{i_0}( R_n^{(b)}( \mathbf{x} ) ) \right)^2 \right]
\le ( h_n^{(b)} )^{-2 \beta} \mathbb{E}\left\| R_n^{(b)}( \mathbf{x} ) \right\|^2
\longrightarrow 0 
\label{eq:coro18}
\end{align}
as $ n \longrightarrow \infty $.

Therefore, from \eqref{eq:coro8}, \eqref{eq:coro10}, \eqref{eq:coro13}, \eqref{eq:coro14} and \eqref{eq:coro18}, we have
\begin{align*}
( h_n^{(b)} )^{-2 \beta} \mathbb{E}\left\| \widehat{\Theta}_n^{(b)}( \mathbf{x} ) - \Theta( \mathbf{x} ) \right\|^2
\ge \frac{c_6}{2} > 0
\end{align*}
for all sufficiently large $ n $.
\end{proof}

\subsection{Results required to prove \autoref{thm:5}}

\begin{lemma} \label{lemma:adapt1}
Let $ 0 < \epsilon_0 < 0.5 $ be fixed.
For $ h \in \mathbb{H}_n $, define
\begin{align*}
& \tilde{D}_n( \mathbf{x}, h ) 
= \frac{1}{(1 + \epsilon_0)} \sigma^2 \zeta_n \frac{\log n}{n \phi( \mathbf{x}, h )} , \\
& \tilde{C}_n( \mathbf{x}, h ) 
= \max_{ h' \in \mathbb{H}_n } \left( \left\| \widehat{\Theta}_n( \mathbf{x}, h' ) - \widehat{\Theta}_n( \mathbf{x}, \max\{ h, h' \} ) \right\|^2 - \tilde{D}_n( \mathbf{x}, h' ) \right)_+ .
\end{align*}
Then,
\begin{align*}
& C_n( \mathbf{x}, h ) 
\le \tilde{C}_n( \mathbf{x}, h )
+ \max_{ h' \in \mathbb{H}_n } \left( \tilde{D}_n( \mathbf{x}, h' ) - D_n( \mathbf{x}, h' ) \right)_+ .
\end{align*}
\end{lemma}
\begin{proof}
The proof is straight forward from the definitions of $ C_n( \mathbf{x}, h ) $, $ D_n( \mathbf{x}, h ) $, $ \tilde{C}_n( \mathbf{x}, h ) $ and $ \tilde{D}_n( \mathbf{x}, h ) $.
\end{proof}

\begin{lemma} \label{lemma:adapt2}
Let $ \tilde{D}_n( \mathbf{x}, h ) $ be as defined in \autoref{lemma:adapt1}, where $ h \in \mathbb{H}_n $. Then, there exists a positive integer $ N_1 $ such that for all $ n \ge N_1 $,
\begin{align*}
& \mathbb{E}\left[ \max_{ h' \in \mathbb{H}_n } \left( \tilde{D}_n( \mathbf{x}, h' ) - D_n( \mathbf{x}, h' ) \right)_+ \right]
< \frac{1}{n^2} , \\
& \text{and}\quad \mathbb{E}\left[ D_n( \mathbf{x}, h ) \right]
\le 
3 \tilde{D}_n( \mathbf{x}, h )
+ \frac{3 \zeta_0 \sigma^2}{n^2} .
\end{align*}
\end{lemma}
\begin{proof}
Define the event
\begin{align*}
\mathbb{U}( \mathbf{x} ) = \bigcap\limits_{h' \in \mathbb{H}_n} \left\{ \left| \frac{\widehat{\phi}( \mathbf{x}, h' )}{\phi( \mathbf{x}, h' )} - 1 \right| < \epsilon_0 \right\} ,
\end{align*}
where $ \epsilon_0 $ is as in \autoref{lemma:adapt1}.
Since the cardinality of $ \mathbb{H}_n $ is at most $ n $, from an application of the Bernstein inequality, we get that there exists an integer $ n_1 $ such that for all $ n \ge n_1 $,
\begin{align}
\mathbb{P}\left[ ( \mathbb{U}( \mathbf{x} ) )^c \right]
& = \mathbb{P}\left[ \bigcup\limits_{h' \in \mathbb{H}_n} \left\{ \left| \widehat{\phi}( \mathbf{x}, h' ) - \phi( \mathbf{x}, h' ) \right| \ge \epsilon_0 \phi( \mathbf{x}, h' ) \right\} \right]
\nonumber\\
& \le
\sum_{h' \in \mathbb{H}_n} \mathbb{P}\left[ \left| \sum_{i=1}^{n} \left[ \mathbb{I}\left( d( \mathbf{x}, \mathbf{X}_i ) \le h' \right) - \phi( \mathbf{x}, h' ) \right] \right| \ge \epsilon_0 n \phi( \mathbf{x}, h' ) \right]
\nonumber\\
& <
2 \sum_{h' \in \mathbb{H}_n} \exp\left[ - 4 \log n \right]
\le \frac{2}{n^3} .
\label{lemmaadapt2:eq1}
\end{align}
Note that
\begin{align}
& \mathbb{E}\left[ \max_{ h' \in \mathbb{H}_n } \left( \tilde{D}_n( \mathbf{x}, h' ) - D_n( \mathbf{x}, h' ) \right)_+ \right]
\nonumber\\
& =
\mathbb{E}\left[ \max_{ h' \in \mathbb{H}_n } \left( \tilde{D}_n( \mathbf{x}, h' ) - D_n( \mathbf{x}, h' ) \right)_+ \mathbb{I}( \mathbb{U}( \mathbf{x} ) ) \right] \nonumber\\
& \quad +
\mathbb{E}\left[ \max_{ h' \in \mathbb{H}_n } \left( \tilde{D}_n( \mathbf{x}, h' ) - D_n( \mathbf{x}, h' ) \right)_+ \mathbb{I}\left( ( \mathbb{U}( \mathbf{x} ) )^c \right) \right] .
\label{lemmaadapt2:eq2}
\end{align}
When $ \mathbb{I}( \mathbb{U}( \mathbf{x} ) ) = 1 $, we have
\begin{align}
& (1 - \epsilon_0) \phi( \mathbf{x}, h' )
< \widehat{\phi}( \mathbf{x}, h' )
< (1 + \epsilon_0) \phi( \mathbf{x}, h' )
\quad
\text{for all } h' \in \mathbb{H}_n \nonumber\\
\iff
& \frac{1}{(1 + \epsilon_0)} \frac{1}{\phi( \mathbf{x}, h' )}
< \frac{1}{\widehat{\phi}( \mathbf{x}, h' )}
< \frac{1}{(1 - \epsilon_0)} \frac{1}{\phi( \mathbf{x}, h' )}
\quad
\text{for all } h' \in \mathbb{H}_n
\label{lemmaadapt2:eq3}\\
\implies
& \max_{ h' \in \mathbb{H}_n } \left( \tilde{D}_n( \mathbf{x}, h' ) - D_n( \mathbf{x}, h' ) \right)_+ \mathbb{I}( \mathbb{U}( \mathbf{x} ) ) = 0
\nonumber\\
\implies
& \mathbb{E}\left[ \max_{ h' \in \mathbb{H}_n } \left( \tilde{D}_n( \mathbf{x}, h' ) - D_n( \mathbf{x}, h' ) \right)_+ \mathbb{I}( \mathbb{U}( \mathbf{x} ) ) \right] = 0 .
\label{lemmaadapt2:eq4}
\end{align}
Let $ n_2 $ be a positive integer such that for all $ n \ge n_2 $, $ \zeta_n \le (1 + \epsilon_0) \zeta_0 $.
So, from \eqref{lemmaadapt2:eq1}, we get that for all $ n \ge \max\{ n_1, n_2 \} $,
\begin{align}
& \mathbb{E}\left[ \max_{ h' \in \mathbb{H}_n } \left( \tilde{D}_n( \mathbf{x}, h' ) - D_n( \mathbf{x}, h' ) \right)_+ \mathbb{I}\left( ( \mathbb{U}( \mathbf{x} ) )^c \right) \right]
\nonumber\\
& \le
\sum_{ h' \in \mathbb{H}_n } \mathbb{E}\left[ \left( \tilde{D}_n( \mathbf{x}, h' ) - D_n( \mathbf{x}, h' ) \right)_+ \mathbb{I}\left( ( \mathbb{U}( \mathbf{x} ) )^c \right) \right]
\nonumber\\
& \le
\sum_{ h' \in \mathbb{H}_n } \tilde{D}_n( \mathbf{x}, h' ) \mathbb{P}\left[ ( \mathbb{U}( \mathbf{x} ) )^c \right]
\nonumber\\
& =
\sum_{ h' \in \mathbb{H}_n } \frac{1}{(1 + \epsilon_0)} \sigma^2 \zeta_n \frac{\log n}{n \phi( \mathbf{x}, h' )} \mathbb{P}\left[ ( \mathbb{U}( \mathbf{x} ) )^c \right]
<
2 \zeta_0 \sigma^2 \frac{1}{\log n} \frac{1}{n^2} .
\label{lemmaadapt2:eq5}
\end{align}
Let $ n_3 = \min\{ n \midil \log n > (2 / (1 + \epsilon_0)) \sigma^2 \zeta_0 \} $. Then, from \eqref{lemmaadapt2:eq2}, \eqref{lemmaadapt2:eq4} and \eqref{lemmaadapt2:eq5}, we get that for all $ n \ge \max\{ n_1, n_2, n_3 \} $,
\begin{align}
& \mathbb{E}\left[ \max_{ h' \in \mathbb{H}_n } \left( D_n( \mathbf{x}, h' ) - \tilde{D}_n( \mathbf{x}, h' ) \right)_+ \right]
< \frac{1}{n^2} .
\label{lemmaadapt2:eq6}
\end{align}
Next, from \eqref{lemmaadapt2:eq1} and \eqref{lemmaadapt2:eq3}, we have for all $ n \ge n_1 $,
\begin{align}
\mathbb{E}\left[ D_n( \mathbf{x}, h ) \right]
& =
\mathbb{E}\left[ D_n( \mathbf{x}, h ) \mathbb{I}\left( \mathbb{U}( \mathbf{x} ) \right) \right]
+ \mathbb{E}\left[ D_n( \mathbf{x}, h ) \mathbb{I}\left( ( \mathbb{U}( \mathbf{x} ) )^c \right) \right]
\nonumber\\
& \le 
\frac{(1 + \epsilon_0)}{(1 - \epsilon_0)} \tilde{D}_n( \mathbf{x}, h )
+ \sigma^2 \zeta_n n \mathbb{P}\left[ ( \mathbb{U}( \mathbf{x} ) )^c \right] \nonumber\\
& < 
3 \tilde{D}_n( \mathbf{x}, h )
+ \frac{3 \zeta_0 \sigma^2}{n^2} .
\label{lemmaadapt2:eq7}
\end{align}
Taking $ N_1 = \max\{ n_1, n_2, n_3 \} $, the proof is complete from \eqref{lemmaadapt2:eq6} and \eqref{lemmaadapt2:eq7}.
\end{proof}

\begin{lemma} \label{lemma:adapt3}
Let the assumptions of \autoref{thm:5} be satisfied. Let $ y > 0 $. We have for all sufficiently large $ n $,
\begin{align*}
& \mathbb{P}\left[ \left\| \sum_{i=1}^{n} \mathbb{L}_\mathbf{x} \left( G( \mathbf{Y}_i ) - \mathbb{E}[ G( \mathbf{Y}_i ) \midil \mathbf{X}_i ] \right) \frac{K( h'^{-1} d( \mathbf{x}, \mathbf{X}_i ) )}{n \mathbb{E}[ K( h'^{-1} d( \mathbf{x}, \mathbf{X} ) ) ]} \right\| 
> y \right] 
\le
n^{-3}
\end{align*}
for all $ h' \in \mathbb{H}_n $. Further, given any $ c_1 > 0 $, $ c_2 > 0 $ and any $ 0 < \epsilon < 1 $, we have, for all sufficiently large $ n $,
\begin{align*}
& \mathbb{P}\left[ \left\| \sum_{i=1}^{n} \mathbb{L}_\mathbf{x} \left( G( \mathbf{Y}_i ) - \mathbb{E}[ G( \mathbf{Y}_i ) \midil \mathbf{X}_i ] \right) \frac{K\left( \frac{d( \mathbf{x}, \mathbf{X}_i )}{h'} \right)}{n \mathbb{E}\left[ K\left( \frac{d( \mathbf{x}, \mathbf{X} )}{h'} \right) \right]} \right\| 
> c_2 \sqrt{c_1 D_n( \mathbf{x}, h' ) + t} \right]
\nonumber\\
& \le
\exp\left[ - \frac{( 1 - \epsilon )^2 l^2 n \phi( \mathbf{x}, h' ) c_2^2 ( c_1 D_n( \mathbf{x}, h' ) + t )}{16 \sigma^2 L^2} \right]
\nonumber\\
& \quad +
\exp\left[ - \frac{( 1 - \epsilon )^2 l^2 n \phi( \mathbf{x}, h' ) c_2 \sqrt{c_1 D_n( \mathbf{x}, h' ) + t}}{16 \sigma L^2} \right] .
\end{align*}
for all $ h' \in \mathbb{H}_n $ and all $ t \ge 0 $.
\end{lemma}
\begin{proof}
We use the following result from \cite{yurinskiui1976exponential}:
Let $ \xi_1, \ldots, \xi_n \in \mathcal{B} $ be independent random elements with
\begin{align*}
\mathbb{E}\| \xi_j \|^m \le (m! / 2) b_j^2 H^{m-2}
\end{align*}
for all integers $ m \ge 2 $. Let
\begin{align*}
\beta_n \ge \mathbb{E}\| \xi_1 + \cdots + \xi_n \| ,
\qquad
U_n^2 = b_1^2 + \cdots + b_n^2 .
\end{align*}
If $ \bar{u} = u - (\beta_n / U_n) > 0 $, then
\begin{align}
\mathbb{P}[ \| \xi_1 + \cdots + \xi_n \| \ge u U_n ]
\le
\exp\left[ - \frac{\bar{u}^2}{8 ( 1 + ( \bar{u} H / 2 U_n ) )} \right] .
\label{eq:expon1}
\end{align}
Now, we choose
\begin{align*}
\xi_i
= \mathbb{L}_\mathbf{x} \left( G( \mathbf{Y}_i ) - \mathbb{E}[ G( \mathbf{Y}_i ) \midil \mathbf{X}_i ] \right) \frac{K( h'^{-1} d( \mathbf{x}, \mathbf{X}_i ) )}{n \mathbb{E}[ K( h'^{-1} d( \mathbf{x}, \mathbf{X} ) ) ]}
\end{align*}
for $ i = 1, \ldots, n $. Since $ \mathcal{B} $ is a type 2 Banach space, from \ref{cond:adaptive2}, we have
\begin{align*}
& \mathbb{E}\| \xi_1 + \cdots + \xi_n \| \nonumber\\
& = \mathbb{E}\left\| \sum_{i=1}^{n} \mathbb{L}_\mathbf{x} \left( G( \mathbf{Y}_i ) - \mathbb{E}[ G( \mathbf{Y}_i ) \midil \mathbf{X}_i ] \right) \frac{K( h'^{-1} d( \mathbf{x}, \mathbf{X}_i ) )}{n \mathbb{E}[ K( h'^{-1} d( \mathbf{x}, \mathbf{X} ) ) ]} \right\|
\nonumber\\
& \le \left[ \mathbb{E}\left\| \sum_{i=1}^{n} \mathbb{L}_\mathbf{x} \left( G( \mathbf{Y}_i ) - \mathbb{E}[ G( \mathbf{Y}_i ) \midil \mathbf{X}_i ] \right) \frac{K( h'^{-1} d( \mathbf{x}, \mathbf{X}_i ) )}{n \mathbb{E}[ K( h'^{-1} d( \mathbf{x}, \mathbf{X} ) ) ]} \right\|^2 \right]^{\frac{1}{2}}
\nonumber\\
& \le \left[ c \sum_{i=1}^{n} \mathbb{E}\left\| \mathbb{L}_\mathbf{x} \left( G( \mathbf{Y}_i ) - \mathbb{E}[ G( \mathbf{Y}_i ) \midil \mathbf{X}_i ] \right) \frac{K( h'^{-1} d( \mathbf{x}, \mathbf{X}_i ) )}{n \mathbb{E}[ K( h'^{-1} d( \mathbf{x}, \mathbf{X} ) ) ]} \right\|^2 \right]^{\frac{1}{2}}
\nonumber\\
& = \sqrt{c} \left[ \sum_{i=1}^{n} \mathbb{E}\left[ \mathbb{E}\left[ \left\| \mathbb{L}_\mathbf{x} \left( G( \mathbf{Y}_i ) - \mathbb{E}[ G( \mathbf{Y}_i ) \midil \mathbf{X}_i ] \right) \right\|^2 \middle\arrowvert \mathbf{X}_i \right] \frac{K^2\left( \frac{d( \mathbf{x}, \mathbf{X}_i )}{h'} \right)}{\left( n \mathbb{E}\left[ K\left( \frac{d( \mathbf{x}, \mathbf{X} )}{h'} \right) \right] \right)^2} \right] \right]^{\frac{1}{2}}
\nonumber\\
& \le \sqrt{c} \frac{\sigma  L}{l \sqrt{n \phi( \mathbf{x}, h' )}} = \beta_n ,
\end{align*}
where $ c $ is a positive constant.
Also, again using \ref{cond:adaptive2}, we get
\begin{align*}
\mathbb{E}\| \xi_i \|^m
& = \mathbb{E}\left\| \mathbb{L}_\mathbf{x} \left( G( \mathbf{Y}_i ) - \mathbb{E}[ G( \mathbf{Y}_i ) \midil \mathbf{X}_i ] \right) \frac{K( h'^{-1} d( \mathbf{x}, \mathbf{X}_i ) )}{n \mathbb{E}[ K( h'^{-1} d( \mathbf{x}, \mathbf{X} ) ) ]} \right\|^m
\nonumber\\
& \le \frac{m!}{2} \left( \frac{\sigma L}{l n \phi( \mathbf{x}, h' )} \right)^{m-2} \frac{\sigma^2 L^2}{l^2 n^2 \phi( \mathbf{x}, h' )} ,
\end{align*}
and we can take
\begin{align*}
U_n^2 = \frac{\sigma^2 L^2}{l^2 n \phi( \mathbf{x}, h' )}
\quad
\text{and}\quad
H = \frac{\sigma L}{l n \phi( \mathbf{x}, h' )} .
\end{align*}
So, $ ( \beta_n / U_n ) = \sqrt{c} $.
Now,
\begin{align*}
\frac{y}{U_n} - \frac{\beta_n}{U_n}
= \frac{y l \sqrt{n \phi( \mathbf{x}, h' )}}{\sigma L} - \sqrt{c}
\ge \frac{y l \log n}{\sigma L} - \sqrt{c} > 0
\end{align*}
for all sufficiently large $ n $ and for all $ h' \in \mathbb{H}_n $.
Also,
\begin{align*}
\left( \frac{y}{U_n} - \frac{\beta_n}{U_n} \right) \frac{H}{2 U_n}
= \left( \frac{y l \sqrt{n \phi( \mathbf{x}, h' )}}{\sigma L} - \sqrt{c} \right) \frac{1}{2 \sqrt{n \phi( \mathbf{x}, h' )}}
< \frac{y l}{2 \sigma L} .
\end{align*} So, from \eqref{eq:expon1}, we get that for all sufficiently large $ n $ (depending on $ y $),
\begin{align*}
& \mathbb{P}\left[ \left\| \sum_{i=1}^{n} \mathbb{L}_\mathbf{x} \left( G( \mathbf{Y}_i ) - \mathbb{E}[ G( \mathbf{Y}_i ) \midil \mathbf{X}_i ] \right) \frac{K( h'^{-1} d( \mathbf{x}, \mathbf{X}_i ) )}{n \mathbb{E}[ K( h'^{-1} d( \mathbf{x}, \mathbf{X} ) ) ]} \right\| 
> y \right]
\nonumber\\
& <
\exp\left[ - \frac{\left( y l \log n - \sqrt{c} \sigma L \right)^2}{8 \sigma^2 L^2 + 4 y l \sigma L} \right]
< \exp\left[ - 3 \log n \right]
= n^{-3} .
\end{align*}
For the next part in the statement of this lemma, we have
\begin{align*}
& \min_{t \ge 0} \frac{c_2 \sqrt{c_1 D_n( \mathbf{x}, h' ) + t}}{U_n}
\ge \sqrt{\log n} \frac{l c_2 \sqrt{c_1 \frac{2}{3} \sigma^2 \zeta_n}}{\sigma L}
> \sqrt{c} = \frac{\beta_n}{U_n}
\end{align*}
for all sufficiently large $ n $ and all $ h' \in \mathbb{H}_n $. Also, given any $ 0 < \epsilon < 1 $, we have, for all sufficiently large $ n $,
\begin{align*}
& \epsilon \frac{c_2 \sqrt{c_1 D_n( \mathbf{x}, h' ) + t}}{U_n}
\ge \epsilon \sqrt{\log n} \left( \frac{l c_2 \sqrt{c_1 \frac{2}{3} \sigma^2 \zeta_n}}{\sigma L} \right)
> \sqrt{c}
\nonumber\\
\implies
& \left( \frac{c_2 \sqrt{c_1 D_n( \mathbf{x}, h' ) + t}}{U_n} - \sqrt{c} \right)^2
> ( 1 - \epsilon )^2 c_2^2 \frac{c_1 D_n( \mathbf{x}, h' ) + t}{U_n^2}
\end{align*}
for all $ h' \in \mathbb{H}_n $ and all $ t \ge 0 $.
Now,
\begin{align*}
\left( \frac{c_2 \sqrt{c_1 D_n( \mathbf{x}, h' ) + t}}{U_n} - \frac{\beta_n}{U_n} \right) \frac{H}{2 U_n}
& \le c_2 \sqrt{c_1 D_n( \mathbf{x}, h' ) + t} \frac{H}{2 U_n^2} \\
& < c_2 \sqrt{c_1 D_n( \mathbf{x}, h' ) + t} \frac{l}{\sigma L} 
\end{align*}
for all $ h' \in \mathbb{H}_n $ and all $ t \ge 0 $.
So, from \eqref{eq:expon1}, we get that for all sufficiently large $ n $,
\begin{align*}
& \mathbb{P}\left[ \left\| \sum_{i=1}^{n} \mathbb{L}_\mathbf{x} \left( G( \mathbf{Y}_i ) - \mathbb{E}[ G( \mathbf{Y}_i ) \midil \mathbf{X}_i ] \right) \frac{K( h'^{-1} d( \mathbf{x}, \mathbf{X}_i ) )}{n \mathbb{E}[ K( h'^{-1} d( \mathbf{x}, \mathbf{X} ) ) ]} \right\| 
> c_2 \sqrt{c_1 D_n( \mathbf{x}, h' ) + t} \right]
\nonumber\\
& \le
\exp\left[ - \frac{( 1 - \epsilon )^2 c_2^2 l^2 n \phi( \mathbf{x}, h' ) ( c_1 D_n( \mathbf{x}, h' ) + t )}{8 \sigma L^2 \left( \sigma + c_2 \sqrt{c_1 D_n( \mathbf{x}, h' ) + t} \right)} \right]
\nonumber\\
& \le
\exp\left[ - \frac{( 1 - \epsilon )^2 l^2 n \phi( \mathbf{x}, h' ) c_2^2 ( c_1 D_n( \mathbf{x}, h' ) + t )}{16 \sigma^2 L^2} \right]
\nonumber\\
& \quad +
\exp\left[ - \frac{( 1 - \epsilon )^2 l^2 n \phi( \mathbf{x}, h' ) c_2 \sqrt{c_1 D_n( \mathbf{x}, h' ) + t}}{16 \sigma L^2} \right] 
\end{align*}
for all $ h' \in \mathbb{H}_n $ and for all $ t \ge 0 $.
\end{proof}

\begin{lemma} \label{lemma:adapt4}
Let $ \tilde{C}_n( \mathbf{x}, h ) $ be as defined in \autoref{lemma:adapt1}, where $ h \in \mathbb{H}_n $.
Let the assumptions in \autoref{thm:5} be satisfied. Then, there exists an integer $ N_2 $ such that for all $ n \ge N_2 $,
\begin{align*}
\tilde{C}_n( \mathbf{x}, h )
& \le
M_1 h^{2 \beta} 
+ 24 \max_{ h' \in \mathbb{H}_n,\, h' \le h } \left( \left\| V_n( \mathbf{x}, h' ) \right\|^2 - \frac{\tilde{D}_n( \mathbf{x}, h' )}{24} \right)_+ \nonumber\\
& \quad + 12 \max_{ h' \in \mathbb{H}_n,\, h' \le h } \left( \left\| R_n( \mathbf{x}, h' ) \right\|^2 - \left( M h'^{2 \beta} + \left\| V_n( \mathbf{x}, h' ) \right\|^2 \right) \right)_+ 
\end{align*}
for all $ h \in \mathbb{H}_n $, where $ M_1 > 0 $ is some constant. Further, for all $ n \ge N_2 $ and all $ h \in \mathbb{H}_n $, we have
\begin{align*}
& \mathbb{P}\left[ \max_{ h' \in \mathbb{H}_n,\, h' \le h } \left( \left\| R_n( \mathbf{x}, h' ) \right\|^2 - \left( M h'^{2 \beta} + \left\| V_n( \mathbf{x}, h' ) \right\|^2 \right) \right)_+
> \frac{1}{n^2} \right]
\le 2 n^{-2} ,
\\
& \text{and} \quad
\mathbb{E}\left[ \max_{ h' \in \mathbb{H}_n,\, h' \le h } \left( \left\| V_n( \mathbf{x}, h' ) \right\|^2 - \frac{\tilde{D}_n( \mathbf{x}, h' )}{24} \right)_+ \right]
< \frac{1}{n} .
\end{align*}
\end{lemma}
\begin{proof}
Note that
\begin{align}
& \tilde{C}_n( \mathbf{x}, h )
\nonumber\\
& = \max_{ h' \in \mathbb{H}_n,\, h' \le h } \left( \left\| \widehat{\Theta}_n( \mathbf{x}, h' ) - \widehat{\Theta}_n( \mathbf{x}, h ) \right\|^2 - \tilde{D}_n( \mathbf{x}, h' ) \right)_+
\nonumber\\
& \le \max_{ h' \in \mathbb{H}_n,\, h' \le h } \left( 2 \left\| \widehat{\Theta}_n( \mathbf{x}, h' ) - \Theta( \mathbf{x} ) \right\|^2
+ 2 \left\| \widehat{\Theta}_n( \mathbf{x}, h ) - \Theta( \mathbf{x} ) \right\|^2
- \tilde{D}_n( \mathbf{x}, h' ) \right)_+
\nonumber\\
& \le 2 \max_{ h' \in \mathbb{H}_n,\, h' \le h } \left( \left\| \widehat{\Theta}_n( \mathbf{x}, h' ) - \Theta( \mathbf{x} ) \right\|^2 - \frac{\tilde{D}_n( \mathbf{x}, h' )}{4} \right)_+ \nonumber\\
& \quad +
2 \max_{ h' \in \mathbb{H}_n,\, h' \le h } \left( \left\| \widehat{\Theta}_n( \mathbf{x}, h ) - \Theta( \mathbf{x} ) \right\|^2
- \frac{\tilde{D}_n( \mathbf{x}, h' )}{4} \right)_+
\nonumber\\
& \le 4 \max_{ h' \in \mathbb{H}_n,\, h' \le h } \left( \left\| \widehat{\Theta}_n( \mathbf{x}, h' ) - \Theta( \mathbf{x} ) \right\|^2 - \frac{\tilde{D}_n( \mathbf{x}, h' )}{4} \right)_+
\label{eq:thmadaptive1}
\end{align}
since $ \tilde{D}_n( \mathbf{x}, h' ) \ge \tilde{D}_n( \mathbf{x}, h ) $ for $ h' \le h $.
From \eqref{eq:main}, we have
\begin{align}
& \max_{ h' \in \mathbb{H}_n,\, h' \le h } \left( \left\| \widehat{\Theta}_n( \mathbf{x}, h' ) - \Theta( \mathbf{x} ) \right\|^2 - \frac{\tilde{D}_n( \mathbf{x}, h' )}{4} \right)_+
\nonumber\\
& \le
3 \max_{ h' \in \mathbb{H}_n,\, h' \le h } \left( \left\| B_n( \mathbf{x}, h' ) \right\|^2 + \left\| V_n( \mathbf{x}, h' ) \right\|^2 + \left\| R_n( \mathbf{x}, h' ) \right\|^2 - \frac{\tilde{D}_n( \mathbf{x}, h' )}{12} \right)_+
\nonumber\\
& \le
3 \max_{ h' \in \mathbb{H}_n,\, h' \le h } \left( \left\| B_n( \mathbf{x}, h' ) \right\|^2 + M h'^{2 \beta} \right) 
\nonumber\\
& \quad + 6 \max_{ h' \in \mathbb{H}_n,\, h' \le h } \left( \left\| V_n( \mathbf{x}, h' ) \right\|^2 - \frac{\tilde{D}_n( \mathbf{x}, h' )}{24} \right)_+ 
\nonumber\\
& \quad + 3 \max_{ h' \in \mathbb{H}_n,\, h' \le h } \left( \left\| R_n( \mathbf{x}, h' ) \right\|^2 - \left( M h'^{2 \beta} + \left\| V_n( \mathbf{x}, h' ) \right\|^2 \right) \right)_+ .
\label{eq:thmadaptive1a}
\end{align}
From assumption \ref{cond:1} and the fact that $ \max\{ h' \midil h' \in \mathbb{H}_n \} \longrightarrow 0 $ as $ n \longrightarrow \infty $, we get that for all sufficiently large $ n $,
\begin{align}
& \max_{ h' \in \mathbb{H}_n,\, h' \le h } \left( \left\| B_n( \mathbf{x}, h' ) \right\|^2 + M h'^{2 \beta} \right)
\le M_1 h^{2 \beta}
\label{eq:thmadaptive1b}
\end{align}
for all $ h \in \mathbb{H}_n $, where $ M_1 > 0 $ is a constant.

Next, define the event
\begin{align*}
\mathbb{S}( \mathbf{x}, h' )
= \left\{ \frac{1}{n} \sum_{i=1}^{n} \frac{K( h'^{-1} d( \mathbf{x}, \mathbf{X}_i ) )}{\mathbb{E}\left[ K( h'^{-1} d( \mathbf{x}, \mathbf{X} ) ) \right]}
> (1 - \epsilon_0) \right\} ,
\end{align*}
where $ \epsilon_0 $ is the number described in \autoref{lemma:adapt1}.
From assumption \ref{cond:adaptive3} and the fact that $ \max\{ h' \midil h' \in \mathbb{H}_n \} \longrightarrow 0 $ as $ n \longrightarrow \infty $, we have for all sufficiently large $ n $,
\begin{align}
& \mathbb{P}\left[ \max_{ h' \in \mathbb{H}_n,\, h' \le h } \left( \left\| R_n( \mathbf{x}, h' ) \right\|^2 - \left( M h'^{2 \beta} + \left\| V_n( \mathbf{x}, h' ) \right\|^2 \right) \right)_+
> \frac{1}{n^2} \right]
\nonumber\\
& \le
\sum_{ h' \in \mathbb{H}_n }\mathbb{P}\left[ \left( \left\| R_n( \mathbf{x}, h' ) \right\|^2 - \left( M h'^{2 \beta} + \left\| V_n( \mathbf{x}, h' ) \right\|^2 \right) \right)_+
> 0 \right]
\nonumber\\
& \le
\sum_{ h' \in \mathbb{H}_n }\mathbb{P}\left[ \left\| R_n( \mathbf{x}, h' ) \right\|^2
> M h'^{2 \beta} + \left\| V_n( \mathbf{x}, h' ) \right\|^2 \right]
\nonumber\\
& \le
\sum_{ h' \in \mathbb{H}_n }\mathbb{P}\left[ \left\| V_n( \mathbf{x}, h' ) \right\|
> \epsilon_2 \right]
\nonumber\\
& \le
\sum_{ h' \in \mathbb{H}_n }\mathbb{P}\left[ \left\| V_n( \mathbf{x}, h' ) \right\| > \epsilon_2 \text{ and } \mathbb{I}( \mathbb{S}( \mathbf{x}, h' ) ) = 1 \right]
+ \sum_{ h' \in \mathbb{H}_n }\mathbb{P}\left[ ( \mathbb{S}( \mathbf{x}, h' ) )^c \right] .
\label{eq:thmadaptive3}
\end{align}
Now, using assumption \ref{assume:a1}, the fact that $ n \phi( \mathbf{x}, h' ) \ge ( \log n )^2 $ for all $ h' \in \mathbb{H}_n $ and the Bernstein inequality, we get that for all sufficiently large $ n $,
\begin{align}
\sum_{ h' \in \mathbb{H}_n }\mathbb{P}\left[ ( \mathbb{S}( \mathbf{x}, h' ) )^c \right] 
& = \sum_{ h' \in \mathbb{H}_n }\mathbb{P}\left[ \frac{1}{n} \sum_{i=1}^{n} \left[ 1 - \frac{K( h'^{-1} d( \mathbf{x}, \mathbf{X}_i ) )}{\mathbb{E}\left[ K( h'^{-1} d( \mathbf{x}, \mathbf{X} ) ) \right]} \right]
\ge \epsilon_0 \right]
\nonumber\\
& \le \sum_{ h' \in \mathbb{H}_n } \exp\left[ - 3 \log n \right]
\le n^{-2} .
\label{eq:thmadaptive4}
\end{align}
Also, from \autoref{lemma:adapt3}, we get
\begin{align}
& \sum_{ h' \in \mathbb{H}_n } \mathbb{P}\left[ \left\| V_n( \mathbf{x}, h' ) \right\| > \epsilon_2 \text{ and } \mathbb{I}( \mathbb{S}( \mathbf{x}, h' ) ) = 1 \right]
\nonumber\\
& \le
\sum_{ h' \in \mathbb{H}_n } \mathbb{P}\left[ \left\| \sum_{i=1}^{n} \mathbb{L}_\mathbf{x} \left( G( \mathbf{Y}_i ) - \mathbb{E}[ G( \mathbf{Y}_i ) \midil \mathbf{X}_i ] \right) \frac{K( h'^{-1} d( \mathbf{x}, \mathbf{X}_i ) )}{n \mathbb{E}[ K( h'^{-1} d( \mathbf{x}, \mathbf{X} ) ) ]} \right\| 
> (1 - \epsilon_0) \epsilon_2 \right]
\nonumber\\
& \le n^{-2}
\label{eq:thmadaptive5}
\end{align}
for all sufficiently large $ n $. Hence, from \eqref{eq:thmadaptive3}, \eqref{eq:thmadaptive4} and \eqref{eq:thmadaptive5}, we have
\begin{align}
& \mathbb{P}\left[ \max_{ h' \in \mathbb{H}_n,\, h' \le h } \left( \left\| R_n( \mathbf{x}, h' ) \right\|^2 - \left( M h'^{2 \beta} + \left\| V_n( \mathbf{x}, h' ) \right\|^2 \right) \right)_+
> \frac{1}{n^2} \right]
\le 2 n^{-2}
\label{eq:thmadaptive5a}
\end{align}
for all sufficiently large $ n $ and all $ h \in \mathbb{H}_n $.
Next,
\begin{align}
& \mathbb{E}\left[ \max_{ h' \in \mathbb{H}_n,\, h' \le h } \left( \left\| V_n( \mathbf{x}, h' ) \right\|^2 - \frac{\tilde{D}_n( \mathbf{x}, h' )}{24} \right)_+ \right]
\nonumber\\
& \le
\sum_{ h' \in \mathbb{H}_n } \mathbb{E}\left[ \left( \left\| V_n( \mathbf{x}, h' ) \right\|^2 - \frac{\tilde{D}_n( \mathbf{x}, h' )}{24} \right)_+ \right]
\nonumber\\
& \le
\sum_{ h' \in \mathbb{H}_n } \mathbb{E}\left[ \left( \left\| V_n( \mathbf{x}, h' ) \right\|^2 \mathbb{I}( \mathbb{S}( \mathbf{x}, h' ) ) - \frac{\tilde{D}_n( \mathbf{x}, h' )}{24} \right)_+ \right] \nonumber\\
& \quad + 
\sum_{ h' \in \mathbb{H}_n } \mathbb{E}\left[ \left\| V_n( \mathbf{x}, h' ) \right\|^2 \mathbb{I}\left( ( \mathbb{S}( \mathbf{x}, h' ) )^c \right) \right] .
\label{eq:thmadaptive6}
\end{align}
Since $ \mathcal{B} $ is a type 2 Banach space, from \ref{cond:adaptive2} and \eqref{eq:thmadaptive4}, we have
\begin{align}
& \sum_{ h' \in \mathbb{H}_n } \mathbb{E}\left[ \left\| V_n( \mathbf{x}, h' ) \right\|^2 \mathbb{I}\left( ( \mathbb{S}( \mathbf{x}, h' ) )^c \right) \right]
\nonumber\\
& =
\sum_{ h' \in \mathbb{H}_n } \mathbb{E}\left[ \frac{\mathbb{E}\left[ \left\| \sum_{i=1}^{n} \mathbb{L}_\mathbf{x}\left( G( \mathbf{Y}_i ) - \mathbb{E}[ G( \mathbf{Y}_i ) \midil \mathbf{X}_i ] \right) K\left( \frac{d( \mathbf{x}, \mathbf{X}_i )}{h'} \right) \right\|^2 \middle\arrowvert \mathbf{X}_1, \ldots, \mathbf{X}_n \right]}{\left( \sum_{i=1}^{n} K\left( \frac{d( \mathbf{x}, \mathbf{X}_i )}{h'} \right) \right)^2} \mathbb{I}\left( ( \mathbb{S}( \mathbf{x}, h' ) )^c \right) \right]
\nonumber\\
& \le
\sum_{ h' \in \mathbb{H}_n } \mathbb{E}\left[ \frac{c \sum_{i=1}^{n} \mathbb{E}\left[ \left\| \mathbb{L}_\mathbf{x}\left( G( \mathbf{Y}_i ) - \mathbb{E}[ G( \mathbf{Y}_i ) \midil \mathbf{X}_i ] \right) \right\|^2 \middle\arrowvert \mathbf{X}_i \right] K^2\left( \frac{d( \mathbf{x}, \mathbf{X}_i )}{h'} \right)}{\left( \sum_{i=1}^{n} K( h'^{-1} d( \mathbf{x}, \mathbf{X}_i ) ) \right)^2} \mathbb{I}\left( ( \mathbb{S}( \mathbf{x}, h' ) )^c \right) \right]
\nonumber\\
& \le
c \sigma^2 \sum_{ h' \in \mathbb{H}_n } \mathbb{E}\left[ \frac{\sum_{i=1}^{n} K^2( h'^{-1} d( \mathbf{x}, \mathbf{X}_i ) )}{\left( \sum_{i=1}^{n} K( h'^{-1} d( \mathbf{x}, \mathbf{X}_i ) ) \right)^2} \mathbb{I}\left( ( \mathbb{S}( \mathbf{x}, h' ) )^c \right) \right]
\nonumber\\
& \le
c \sigma^2 \sum_{ h' \in \mathbb{H}_n } \mathbb{P}\left[ ( \mathbb{S}( \mathbf{x}, h' ) )^c \right]
\le c \sigma^2 n^{-2} 
\label{eq:thmadaptive7}
\end{align}
for all sufficiently large $ n $, where $ c > 0 $ is a constant.
On the other hand, taking $ \epsilon = \epsilon_0 $ in \autoref{lemma:adapt3}, we have for all sufficiently large $ n $,
\begin{align}
& \sum_{ h' \in \mathbb{H}_n } \mathbb{E}\left[ \left( \left\| V_n( \mathbf{x}, h' ) \right\|^2 \mathbb{I}( \mathbb{S}( \mathbf{x}, h' ) ) - \frac{\tilde{D}_n( \mathbf{x}, h' )}{24} \right)_+ \right]
\nonumber\\
& =
\sum_{ h' \in \mathbb{H}_n } \int_{0}^{\infty} \mathbb{P}\left[ \left( \left\| V_n( \mathbf{x}, h' ) \right\|^2 \mathbb{I}( \mathbb{S}( \mathbf{x}, h' ) ) - \frac{\tilde{D}_n( \mathbf{x}, h' )}{24} \right)_+ \ge t \right] dt
\nonumber\\
& =
\sum_{ h' \in \mathbb{H}_n } \int_{0}^{\infty} \mathbb{P}\left[ \left\| V_n( \mathbf{x}, h' ) \right\| \mathbb{I}( \mathbb{S}( \mathbf{x}, h' ) ) \ge \sqrt{\frac{\tilde{D}_n( \mathbf{x}, h' )}{24} + t} \right] dt
\nonumber\\
& \le
\sum_{ h' \in \mathbb{H}_n } \int_{0}^{\infty} \mathbb{P}\left[ \left\| \sum_{i=1}^{n} \frac{\mathbb{L}_\mathbf{x}\left( G( \mathbf{Y}_i ) - \mathbb{E}[ G( \mathbf{Y}_i ) \midil \mathbf{X}_i ] \right) K\left( \frac{d( \mathbf{x}, \mathbf{X}_i )}{h'} \right)}{n \mathbb{E}\left[ K\left( \frac{d( \mathbf{x}, \mathbf{X}_i )}{h'} \right) \right]} \right\| 
\ge (1 - \epsilon_0) \sqrt{\frac{\tilde{D}_n( \mathbf{x}, h' )}{24} + t} \right] dt
\nonumber\\
& \le
\sum_{ h' \in \mathbb{H}_n } \int_{0}^{\infty} \exp\left[ - \frac{( 1 - \epsilon_0 )^4 l^2 n \phi( \mathbf{x}, h' )}{16 \sigma^2 L^2} \left( \frac{1}{24} D_n( \mathbf{x}, h' ) + t \right) \right] dt
\nonumber\\
& \quad +
\sum_{ h' \in \mathbb{H}_n } \int_{0}^{\infty} \exp\left[ - \frac{( 1 - \epsilon_0 )^3 l^2 n \phi( \mathbf{x}, h' )}{16 \sigma L^2} \sqrt{\frac{1}{24} D_n( \mathbf{x}, h' ) + t} \right] dt .
\label{eq:thmadaptive8}
\end{align}
Now, for the second term on the right hand side of \eqref{eq:thmadaptive8}, we have
\begin{align}
& \sum_{ h' \in \mathbb{H}_n } \int_{0}^{\infty} \exp\left[ - \frac{( 1 - \epsilon_0 )^3 l^2 n \phi( \mathbf{x}, h' )}{16 \sigma L^2} \sqrt{\frac{1}{24} D_n( \mathbf{x}, h' ) + t} \right] dt
\nonumber\\
& = 2 \sum_{ h' \in \mathbb{H}_n } \int_{\sqrt{\frac{1}{24} D_n( \mathbf{x}, h' )}}^{\infty} \exp\left[ - \frac{( 1 - \epsilon_0 )^3 l^2 n \phi( \mathbf{x}, h' )}{16 \sigma L^2} s \right] s ds
\nonumber\\
& < \frac{1}{n \log n}
\label{eq:thmadaptive9}
\end{align}
for all sufficiently large $ n $.
Next, we take
\begin{align}
\zeta_0 \ge 768 \frac{( 1 + \epsilon_0 )^2}{( 1 - \epsilon_0 )^4} \frac{L^2}{l^2} .
\label{eq:adaptivezetabound}
\end{align}
Since $ \zeta_n \longrightarrow \zeta_0 $ as $ n \longrightarrow \infty $, we have
\begin{align}
\zeta_n > 768 \frac{( 1 + \epsilon_0 )}{( 1 - \epsilon_0 )^4} \frac{L^2}{l^2}
\label{eq:adaptivezeta_nbound}
\end{align}
for all sufficiently large $ n $.
Consequently, for the first term on the right hand side of \eqref{eq:thmadaptive8}, we have from \eqref{eq:adaptivezeta_nbound},
\begin{align}
& \sum_{ h' \in \mathbb{H}_n } \int_{0}^{\infty} \exp\left[ - \frac{( 1 - \epsilon_0 )^4 l^2 n \phi( \mathbf{x}, h' )}{16 \sigma^2 L^2} \left( \frac{1}{24} D_n( \mathbf{x}, h' ) + t \right) \right] dt
\nonumber\\
& = \sum_{ h' \in \mathbb{H}_n } \int_{\frac{1}{24} D_n( \mathbf{x}, h' )}^{\infty} \exp\left[ - \frac{( 1 - \epsilon_0 )^4 l^2 n \phi( \mathbf{x}, h' )}{16 \sigma^2 L^2} s \right] ds
\nonumber\\
& = \sum_{ h' \in \mathbb{H}_n } \frac{16 \sigma^2 L^2}{( 1 - \epsilon_0 )^4 l^2 n \phi( \mathbf{x}, h' )} \exp\left[ - \frac{( 1 - \epsilon_0 )^4 l^2 n \phi( \mathbf{x}, h' )}{16 \sigma^2 L^2} \frac{1}{24} D_n( \mathbf{x}, h' ) \right]
\nonumber\\
& = \sum_{ h' \in \mathbb{H}_n } \frac{16 \sigma^2 L^2}{( 1 - \epsilon_0 )^4 l^2 n \phi( \mathbf{x}, h' )} \exp\left[ - \frac{1}{768} \frac{( 1 - \epsilon_0 )^4}{( 1 + \epsilon_0 )} \frac{l^2}{L^2} \zeta_n ( 2 \log n ) \right]
\nonumber\\
& \le \frac{16 \sigma^2 L^2}{( 1 - \epsilon )^4 l^2 ( \log n )^2} n^{-1}
< \frac{1}{n \log n}
\label{eq:thmadaptive10}
\end{align}
for all sufficiently large $ n $.
Hence, from \eqref{eq:thmadaptive8}, \eqref{eq:thmadaptive9} and \eqref{eq:thmadaptive10}, we have
\begin{align}
& \sum_{ h' \in \mathbb{H}_n } \mathbb{E}\left[ \left( \left\| V_n( \mathbf{x}, h' ) \right\|^2 \mathbb{I}( \mathbb{S}( \mathbf{x}, h' ) ) - \frac{\tilde{D}_n( \mathbf{x}, h' )}{24} \right)_+ \right]
< \frac{2}{n \log n}
\label{eq:thmadaptive11}
\end{align}
for all sufficiently large $ n $.
Therefore, from \eqref{eq:thmadaptive6}, \eqref{eq:thmadaptive7} and \eqref{eq:thmadaptive11}, we get that for all sufficiently large $ n $ and all $ h \in \mathbb{H}_n $,
\begin{align}
& \mathbb{E}\left[ \max_{ h' \in \mathbb{H}_n,\, h' \le h } \left( \left\| V_n( \mathbf{x}, h' ) \right\|^2 - \frac{\tilde{D}_n( \mathbf{x}, h' )}{24} \right)_+ \right]
< \frac{1}{n} .
\label{eq:thmadaptive12}
\end{align}
We choose an integer $ N_2 $ large enough such that the assertions in \eqref{eq:thmadaptive1b}, \eqref{eq:thmadaptive5a} and \eqref{eq:thmadaptive12} are satisfied for all $ n \ge N_2 $ and all $ h \in \mathbb{H}_n $. Hence, the proof is complete from \eqref{eq:thmadaptive1}, \eqref{eq:thmadaptive1a}, \eqref{eq:thmadaptive1b}, \eqref{eq:thmadaptive5a} and \eqref{eq:thmadaptive12}.
\end{proof}

From \eqref{eq:adaptivezetabound}, we see that $ \zeta_0 $ depends on the choice of $ \epsilon_0 $, and it increases with an increase in the value of $ \epsilon_0 $.
Taking $ \epsilon_0 = 0.1 $ we see that
\begin{align}
\zeta_0 = 1500 \frac{L^2}{l^2}
\label{eq:adaptivezetaboundfinal}
\end{align}
satisfies \eqref{eq:adaptivezetabound}. Taking smaller values of $ \epsilon_0 $, we can further decrease the value of $ \zeta_0 $, but it cannot be less than 768 in view of \eqref{eq:adaptivezetabound}.

\bibliographystyle{apa}
\bibliography{bibliography_database}

\begin{thebibliography}{}

\bibitem[\protect\astroncite{Aerts and Claeskens}{1997}]{aerts1997local}
Aerts, M. and Claeskens, G. (1997).
\newblock Local polynomial estimation in multiparameter likelihood models.
\newblock {\em Journal of the American Statistical Association},
  92(440):1536--1545.

\bibitem[\protect\astroncite{Araujo and Gin{\'e}}{1980}]{araujo1980central}
Araujo, A. and Gin{\'e}, E. (1980).
\newblock {\em The central limit theorem for real and {Banach} valued random
  variables}.
\newblock New York: Wiley.

\bibitem[\protect\astroncite{Bhatia}{2009}]{bhatia2009notes}
Bhatia, R. (2009).
\newblock {\em Notes on functional analysis}.
\newblock Hindustan Book Agency, New Delhi.

\bibitem[\protect\astroncite{Burba et~al.}{2009}]{burba2009k}
Burba, F., Ferraty, F., and Vieu, P. (2009).
\newblock k-nearest neighbour method in functional nonparametric regression.
\newblock {\em Journal of Nonparametric Statistics}, 21(4):453--469.

\bibitem[\protect\astroncite{Cameron and
  Martin}{1944}]{cameron1944transformations}
Cameron, R.~H. and Martin, W.~T. (1944).
\newblock Transformations of {Weiner} integrals under translations.
\newblock {\em Annals of Mathematics}, 45(2):386--396.

\bibitem[\protect\astroncite{Chagny and Roche}{2014}]{chagny2014adaptive}
Chagny, G. and Roche, A. (2014).
\newblock Adaptive and minimax estimation of the cumulative distribution
  function given a functional covariate.
\newblock {\em Electronic Journal of Statistics}, 8(2):2352--2404.

\bibitem[\protect\astroncite{Chagny and Roche}{2016}]{chagny2016adaptive}
Chagny, G. and Roche, A. (2016).
\newblock Adaptive estimation in the functional nonparametric regression model.
\newblock {\em Journal of Multivariate Analysis}, 146:105--118.

\bibitem[\protect\astroncite{Chaouch and
  La{\"\i}b}{2013}]{chaouch2013nonparametric}
Chaouch, M. and La{\"\i}b, N. (2013).
\newblock Nonparametric multivariate $ \mbox{L}_1 $-median regression
  estimation with functional covariates.
\newblock {\em Electronic Journal of Statistics}, 7:1553--1586.

\bibitem[\protect\astroncite{Chaouch and La{\"\i}b}{2015}]{chaouch2015vector}
Chaouch, M. and La{\"\i}b, N. (2015).
\newblock Vector-on-function quantile regression for stationary ergodic
  processes.
\newblock {\em Journal of the Korean Statistical Society}, 44(2):161--178.

\bibitem[\protect\astroncite{Chaudhuri and
  Dewanji}{1995}]{chaudhuri1995likelihood}
Chaudhuri, P. and Dewanji, A. (1995).
\newblock On a likelihood-based approach in nonparametric smoothing and
  cross-validation.
\newblock {\em Statistics \& Probability Letters}, 22(1):7--15.

\bibitem[\protect\astroncite{Dereich and
  Lifshits}{2005}]{dereich2005probabilities}
Dereich, S. and Lifshits, M. (2005).
\newblock Probabilities of randomly centered small balls and quantization in
  {Banach} spaces.
\newblock {\em The Annals of Probability}, 33(4):1397--1421.

\bibitem[\protect\astroncite{Dette et~al.}{2012}]{dette2012testing}
Dette, H., Marchlewski, M., and Wagener, J. (2012).
\newblock Testing for a constant coefficient of variation in nonparametric
  regression by empirical processes.
\newblock {\em Annals of the Institute of Statistical Mathematics},
  64(5):1045--1070.

\bibitem[\protect\astroncite{Dette and Wieczorek}{2009}]{dette2009testing}
Dette, H. and Wieczorek, G. (2009).
\newblock Testing for a constant coefficient of variation in nonparametric
  regression.
\newblock {\em Journal of Statistical Theory and Practice}, 3(3):587--612.

\bibitem[\protect\astroncite{Donoho and Liu}{1991a}]{donoho1991geometrizingII}
Donoho, D.~L. and Liu, R.~C. (1991a).
\newblock Geometrizing rates of convergence, {II}.
\newblock {\em The Annals of Statistics}, 19(2):633--667.

\bibitem[\protect\astroncite{Donoho and Liu}{1991b}]{donoho1991geometrizingIII}
Donoho, D.~L. and Liu, R.~C. (1991b).
\newblock Geometrizing rates of convergence, {III}.
\newblock {\em The Annals of Statistics}, 19(2):668--701.

\bibitem[\protect\astroncite{Ferraty et~al.}{2010}]{ferraty2010rate}
Ferraty, F., Laksaci, A., Tadj, A., and Vieu, P. (2010).
\newblock Rate of uniform consistency for nonparametric estimates with
  functional variables.
\newblock {\em Journal of Statistical Planning and Inference}, 140(2):335--352.

\bibitem[\protect\astroncite{Ferraty et~al.}{2006}]{ferraty2006estimating}
Ferraty, F., Laksaci, A., and Vieu, P. (2006).
\newblock Estimating some characteristics of the conditional distribution in
  nonparametric functional models.
\newblock {\em Statistical Inference for Stochastic Processes}, 9(1):47--76.

\bibitem[\protect\astroncite{Ferraty et~al.}{2007}]{ferraty2007nonparametric}
Ferraty, F., Mas, A., and Vieu, P. (2007).
\newblock Nonparametric regression on functional data: Inference and practical
  aspects.
\newblock {\em Australian \& New Zealand Journal of Statistics},
  49(3):267--286.

\bibitem[\protect\astroncite{Ferraty et~al.}{2011}]{ferraty2011estimation}
Ferraty, F., Park, J., and Vieu, P. (2011).
\newblock Estimation of a functional single index model.
\newblock In Ferraty, F., editor, {\em Recent Advances in Functional Data
  Analysis and Related Topics}, chapter~17, pages 111--116. New York: Springer.

\bibitem[\protect\astroncite{Ferraty et~al.}{2012}]{ferraty2012regression}
Ferraty, F., Van~Keilegom, I., and Vieu, P. (2012).
\newblock Regression when both response and predictor are functions.
\newblock {\em Journal of Multivariate Analysis}, 109:10--28.

\bibitem[\protect\astroncite{Ferraty and Vieu}{2006}]{ferraty2006nonparametric}
Ferraty, F. and Vieu, P. (2006).
\newblock {\em Nonparametric functional data analysis: Theory and practice}.
\newblock New York: Springer.

\bibitem[\protect\astroncite{Ferr{\'e} and Yao}{2005}]{ferre2005smoothed}
Ferr{\'e}, L. and Yao, A. (2005).
\newblock Smoothed functional inverse regression.
\newblock {\em Statistica Sinica}, 15(3):665.

\bibitem[\protect\astroncite{Ferr{\'e} and Yao}{2003}]{ferre2003functional}
Ferr{\'e}, L. and Yao, A.-F. (2003).
\newblock Functional sliced inverse regression analysis.
\newblock {\em Statistics}, 37(6):475--488.

\bibitem[\protect\astroncite{Hardle}{1990}]{hardle1990applied}
Hardle, W. (1990).
\newblock Applied nonparametric regression.
\newblock {\em Cambridge University Press, UK}.

\bibitem[\protect\astroncite{Hoffmann-Jorgensen
  et~al.}{1979}]{hoffmann1979lower}
Hoffmann-Jorgensen, J., Shepp, L.~A., and Dudley, R.~M. (1979).
\newblock On the lower tail of {Gaussian} seminorms.
\newblock {\em The Annals of Probability}, 7(2):319--342.

\bibitem[\protect\astroncite{Ibragimov and
  Ha\'{s}minskii}{1980}]{ibragimov1980nonparametric}
Ibragimov, I.~A. and Ha\'{s}minskii, R.~Z. (1980).
\newblock On nonparametric estimation of regression.
\newblock {\em Soviet Mathematics Doklady}, 21:810--814.

\bibitem[\protect\astroncite{Klemel{\"a}}{2014}]{klemela2014multivariate}
Klemel{\"a}, J.~S. (2014).
\newblock {\em Multivariate nonparametric regression and visualization: With R
  and applications to finance}.
\newblock Hoboken: Wiley.

\bibitem[\protect\astroncite{Kundu et~al.}{2000}]{kundu2000central}
Kundu, S., Majumdar, S., and Mukherjee, K. (2000).
\newblock Central limit theorems revisited.
\newblock {\em Statistics \& Probability Letters}, 47(3):265--275.

\bibitem[\protect\astroncite{Li}{2001}]{li2001small}
Li, W.~V. (2001).
\newblock Small ball probabilities for {Gaussian} {Markov} processes under the
  $ \mbox{L}_p $-norm.
\newblock {\em Stochastic Processes and Their Applications}, 92(1):87--102.

\bibitem[\protect\astroncite{Li and Shao}{2001}]{li2001gaussian}
Li, W.~V. and Shao, Q.-M. (2001).
\newblock Gaussian processes: Inequalities, small ball probabilities and
  applications.
\newblock {\em Stochastic Processes: Theory and Methods}, 19:533--597.

\bibitem[\protect\astroncite{Lian}{2012}]{lian2012convergence}
Lian, H. (2012).
\newblock Convergence of nonparametric functional regression estimates with
  functional responses.
\newblock {\em Electronic Journal of Statistics}, 6:1373--1391.

\bibitem[\protect\astroncite{Lifshits}{2013}]{lifshits2013gaussian}
Lifshits, M.~A. (2013).
\newblock {\em Gaussian random functions}, volume 322.
\newblock Dordrecht: Springer.

\bibitem[\protect\astroncite{Luki{\'c} and Beder}{2001}]{lukic2001stochastic}
Luki{\'c}, M. and Beder, J. (2001).
\newblock Stochastic processes with sample paths in reproducing kernel
  {Hilbert} spaces.
\newblock {\em Transactions of the American Mathematical Society},
  353(10):3945--3969.

\bibitem[\protect\astroncite{Mas}{2012}]{mas2012lower}
Mas, A. (2012).
\newblock Lower bound in regression for functional data by representation of
  small ball probabilities.
\newblock {\em Electronic Journal of Statistics}, 6:1745--1778.

\bibitem[\protect\astroncite{Masry}{2005}]{masry2005nonparametric}
Masry, E. (2005).
\newblock Nonparametric regression estimation for dependent functional data:
  Asymptotic normality.
\newblock {\em Stochastic Processes and Their Applications}, 115(1):155--177.

\bibitem[\protect\astroncite{Nadaraya}{1964}]{nadaraya1964estimating}
Nadaraya, E.~A. (1964).
\newblock On estimating regression.
\newblock {\em Theory of Probability \& Its Applications}, 9(1):141--142.

\bibitem[\protect\astroncite{{\O}ksendal}{2003}]{oksendal2003stochastic}
{\O}ksendal, B. (2003).
\newblock {\em Stochastic differential equations: An introduction with
  applications}.
\newblock New York: Springer.

\bibitem[\protect\astroncite{Rachdi and Vieu}{2007}]{rachdi2007nonparametric}
Rachdi, M. and Vieu, P. (2007).
\newblock Nonparametric regression for functional data: Automatic smoothing
  parameter selection.
\newblock {\em Journal of Statistical Planning and Inference},
  137(9):2784--2801.

\bibitem[\protect\astroncite{Serfling}{2009}]{serfling2009approximation}
Serfling, R.~J. (2009).
\newblock {\em Approximation theorems of mathematical Statistics}, volume 162.
\newblock Hoboken: Wiley.

\bibitem[\protect\astroncite{Staniswalis}{1989}]{staniswalis1989kernel}
Staniswalis, J.~G. (1989).
\newblock The kernel estimate of a regression function in likelihood-based
  models.
\newblock {\em Journal of the American Statistical Association},
  84(405):276--283.

\bibitem[\protect\astroncite{Stone}{1980}]{stone1980optimal}
Stone, C.~J. (1980).
\newblock Optimal rates of convergence for nonparametric estimators.
\newblock {\em The Annals of Statistics}, 8(6):1348--1360.

\bibitem[\protect\astroncite{Stone}{1982}]{stone1982optimal}
Stone, C.~J. (1982).
\newblock Optimal global rates of convergence for nonparametric regression.
\newblock {\em The Annals of Statistics}, 10(4):1040--1053.

\bibitem[\protect\astroncite{Vepakomma et~al.}{2016}]{vepakomma2016supervised}
Vepakomma, P., Tonde, C., and Elgammal, A. (2016).
\newblock Supervised dimensionality reduction via distance correlation
  maximization.
\newblock {\em arXiv preprint arXiv:1601.00236}.

\bibitem[\protect\astroncite{Watson}{1964}]{watson1964smooth}
Watson, G.~S. (1964).
\newblock Smooth regression analysis.
\newblock {\em Sankhy{\=a}: The Indian Journal of Statistics, Series A},
  26(4):359--372.

\bibitem[\protect\astroncite{Yatracos}{1988}]{yatracos1988lower}
Yatracos, Y.~G. (1988).
\newblock A lower bound on the error in nonparametric regression type problems.
\newblock {\em The Annals of Statistics}, 16(3):1180--1187.

\bibitem[\protect\astroncite{Yurinski{\u\i}}{1976}]{yurinskiui1976exponential}
Yurinski{\u\i}, V. (1976).
\newblock Exponential inequalities for sums of random vectors.
\newblock {\em Journal of Multivariate Analysis}, 6(4):473--499.

\end{thebibliography}

\end{document}